\documentclass[twocolumn, headings=normal]{scrartcl}

\usepackage[english]{babel}
\usepackage[utf8x]{inputenc}
\usepackage[T1]{fontenc}

\usepackage{amsmath}
\usepackage{amsthm}
\usepackage{amsfonts}
\usepackage{graphicx}
\usepackage[colorinlistoftodos]{todonotes}
\usepackage[colorlinks=true, allcolors=blue]{hyperref}

\usepackage[left=1.5cm, right=2cm]{geometry}
\usepackage{graphicx,color,psfrag}
\newcommand{\OPV}[1]{\boldsymbol{#1}}

\newcommand{\OPstatevector}[0]{\OPV{x}}

\DeclareMathAlphabet{\mathitbf}{OML}{cmm}{b}{it}

\def\min{\mbox{min}}

\newcommand{\mbeq}{\overset{!}{=}}

\newcommand{\mat}[1]{\mathchoice{\displaystyle\mathbf#1}
{\textstyle\mathbf#1}{\scriptstyle\mathbf#1}
{\scriptscriptstyle\mathbf#1}}
\newcommand{\vek}[1]{\mathchoice{\displaystyle\boldsymbol#1}
{\textstyle\boldsymbol#1}{\scriptstyle\boldsymbol#1}
{\scriptscriptstyle\boldsymbol#1}}

\newcommand{\ignore}[1]{}

\newtheorem{thm}{Theorem}[section]

\newcommand{\refeq}[1]{Eq.~(\ref{#1})}
\newcommand{\refalg}[1]{Alg.~(\ref{#1})}
\newcommand{\refig}[1]{Fig.~(\ref{#1})}

\newcommand{\vDelta}{\varDelta}

\newcommand{\vOmega}{\varOmega}

\usepackage[percent]{overpic}

\usepackage{algorithm}
\usepackage{algpseudocode}

\newtheorem{definition}[thm]{Definition}
\newtheorem{theorem}[thm]{Theorem}

\usepackage{epstopdf}
\usepackage{pdfpages}

\title{Stochastic state estimation via incremental iterative sparse polynomial chaos based Bayesian-Gauss-Newton-Markov-Kalman filter}
\author{Bojana Rosi\'{c}\\
Applied Mechanics and Data Analysis\\
University of Twente\\
Netherlands}

\begin{document}
\maketitle

\begin{abstract}
\textit{
In this paper is proposed a novel incremental iterative Gauss-Newton-Markov-Kalman filter method for state 
estimation of dynamic models given noisy measurements. The filter is constructed by projecting the random variable representing the unknown state onto 
the subspace generated by data.
The approximation of projection, i.e. the conditional expectation of the state given data, is evaluated by minimising 
the expected Bregman's loss.
The mathematical
formulation of the proposed filter is based on the construction of an optimal nonlinear map between the observable and 
parameter (state) spaces via
a convergent sequence of linear maps obtained by successive linearisation of the observation operator in a Gauss-Newton-like form. 
To allow automatic linearisation of the dynamical system in a sparse form, 
the smoother is designed in a hierarchical setting such that the forward map and its linearised counterpart are estimated in a Bayesian manner
given a forecasted data set. For this purpose the relevance vector machine approach is used. To improve the algorithm convergence,
the smoother is further reformulated in its incremental form in which the current and intermediate states are assimilated before the initial one,
and the corresponding posterior estimates 
are taken as pseudo-measurements. As the latter ones are random variables, 
and not deterministic any more, the novel stochastic iterative filter is designed 
to take this into account. To correct the bias in the posterior outcome, the procedure is built in a predictor-corrector form 
in which the predictor phase is used to assimilate noisy measurement data, whereas the corrector phase is constructed to correct the mean bias. 
The resulting filter is further discretised via 
time-adapting sparse polynomial chaos expansions obtained
either via modified Gram-Schmidt orthogonalisation 
or by a carefully chosen nonlinear mapping, both of which are estimated in a Bayesian manner by promoting 
the sparsity of the outcomes. The time adaptive basis with non-Gaussian arguments is further mapped to the polynomial chaos one
by a suitably chosen isoprobabilistic transformation. 
Finally, the proposed method is tested on a chaotic nonlinear Lorenz 1984 system.
%
%
%
%
}

\end{abstract}
\section{Introduction}

Probabilistic inverse estimation is gaining momentum 
in computational practice today. Bayes's rule as given in its classical form often cannot be used in
practice because the evaluation of the posterior 
distribution requires the use of slowly convergent random walk strategies such as Markov chain Monte 
Carlo-like algorithms \cite{gamerman1997,Bazargan2013, Smith93}. On 
the other hand,  its linear approximation 
in the form of a Kalman filter \cite{kalman} became a very important industrial tool for the prediction/forecast of the system state describing 
various types of dynamical systems.
 However, Kalman filters are not good at coping with highly nonlinear system responses, 
and many attempts have been made to resolve this issue. The vast majority of studies on this subject can be broadly classified into two groups: 
stochastic strategies based on the sequential Monte Carlo algorithm also known as particle/ensemble filters 
(e.g.~\cite{ristic,reich2016}), and deterministic methods based on the linearisation of the measurement operator such as extended \cite{einicke1999robust,kai2010robust} 
and unscented \cite{wan2000unscented,van2001square} Kalman filters. The former theories are based on the approximation of
the posterior distribution via a convex combination of the Diract delta measure such that the
corresponding filter requires only few simulation calls. But, it is well known that the ensemble in the particle form
may collapse, which is especially evident for small ensembles. On the other hand, the deterministic 
filters based on the first order Taylor expansion of the measurement operator may become inaccurate 
when used in a highly nonlinear setting. 


It is well known that the Bayesian update is theoretically based
on the notion of conditional expectation \cite{ABobrowski}. Here the conditional expectation is not only used as
a theoretical basis, but also as a basic computational tool for the identification of the initial state of the dynamical system. 
Being a unique optimal projector for all Bregman's loss functions, the conditional expectation allows the estimation of the
posterior moments by finding an optimal map between the measurement and the parameter/state space 
that minimises the expected Bregman's loss.  Therefore,
being able to numerically approximate conditional expectations, one can build various filtering techniques
for the state assimilation. To accommodate the nonlinearities present in the estimation problem,
in this paper an iterative version of the filter in the Gauss-Newton form is suggested for the backpropagation of information 
on the state in the current time moment to the initial one. Several previous studies have investigated the linearisation
idea by building the filter either 
as an iterative version of ensemble Kalman filters as presented in 
 \cite{sakov2012iterative,bell}, or procedures coming from the randomised likelihood (e.g.~\cite{chen2012ensemble}) and maximum
 a posteriori error estimate 
(e.g.~\cite{wang2016randomized}). In this paper the iterative filtering technique is based on the approximation of the conditional 
expectation of the state given observation, as well as its inverse map, via a sequence of linearised maps obtained by 
minimising the corresponding expected quadratic Bregman's loss functions, or by using Bayesian estimation.
In this manner the Gauss-Newton filtering procedure obtains its hierarchical structure and does not require special differentiation 
techniques as the estimation of the Jacobian comes as the by-product. 
To improve the local convergence, the Gauss-Newton estimation is here improved by substituting the direct state estimation with the incremental one
based on the pseudo-time discretisations. 
 The idea is to build the optimal map between the observation and the initial state as a composition of linearised maps displaying the
intermediate state posteriors chacareterised by pseudo-time discretisations. In contrast to the direct estimation this approach takes the 
estimated intermediate states as pseudo-measurements for the preceding ones. Hence, the dynamic of the
filter's incremental form is driven by pseudo-time stepping in which the global optimal linear map of one update step is
substituted by few optimal local maps obtained by splitting the update step into smaller increments (pseudo-update steps). As the pseudo-measurements
are random variables
and not deterministic ones, here is suggested a novel stochastic Gauss-Newton filter for the state estimation in a predictor-corrector form. 

In contrast to
most sampling approaches to Bayesian updating that typically start from the classical formulation involving conditional measures and densities,
the conditional expectation as the computationally prime object allows a direct estimation of the posterior random variable 
in a functional approximation form. As a stochastic Gauss-Newton filter operates on random variables, not densities, its numerical implementation 
is achieved by discretising the random variables of consideration 
 via time dependent polynomial chaos expansion (PCEs). 
 The time adaptive nature of discretisation is used to prevent an overestimation of the 
 measurement prediction after 
 long-time integration, which is known to be a side-effect of the classical polynomial chaos representations. Therefore, 
 the observation random variables are first discretised in a non-Gaussian basis, 
 which is further transformed to the Gaussian one by a nonlinear isoprobabilistic transformation. 
 The non-Gaussian basis is chosen either as an orthogonal one 
 by employing the
 stochastic modified Gram-Schmidt orthogonalisation as already discussed by \cite{Gerritsma:2010} for purely uncertainty
 quantification purposes, or as a non-orthogonal one 
 taking the form of a 
 nonlinear polynomial map between two consecutive states. 
%
 To promote for sparsity,  the functional representations are estimated in a data-driven Bayesian way 
 by using the relevance vector machine approach \cite{Tipping}. By using the sparse time dependent PCE approximations, the filter is finally designed in 
 its minimal form that is estimated by using a minimal number of model evaluations.

The paper is organised as follows: Section \ref{ss:model} gives a concise introduction to the Bayesian 
state estimation of the abstract dynamical system. Section \ref{ss:cond} considers the approximate Bayesian estimation from a conditional 
expectation point of view. Numerical approximations of conditional expectation are shortly studied in 
Section \ref{opmap}, and hence the Gauss-Newton filtering procedure is introduced. The Bayesian point of view
on the Gauss-Newton filter is further studied in
Section \ref{ss:unbiased_filter}, whereas its incremental version in predictor-corrector form is discussed in Section \ref{pcp}. The filter
discretisation
and its computational form are given in Section \ref{itpce}. Here the filter is studied from the perspective of time adaptive sparse 
random variable discretisations. The paper is concluded with Section \ref{concl}. 
%
%

\section{Model problem}\label{ss:model}
Let the state of the dynamical system ${x} \in \mathbb{R}^d$ satisfying the nonlinear initial value problem 
\begin{equation}
\label{initvp}
 \dot{x}=f(x,t), \quad x_0=x(0)
\end{equation}
be observed in time moments $0\leq t_k \leq T, \textrm{ } t_k=k\Delta t, k\in \mathbb{N}_0$ given time increment $\Delta t$ via 
\begin{equation}
\label{det_observ}
 {y}_{k}=Y(x_n)
\end{equation}
in which $Y$ is a nonlinear observation operator, whereas $x_n, \textrm{ }n\in \mathbb{N}_0$ either denotes the current state when $n=k$,
or an unknown previous state
when $x_n=x_{k-r}$ for $r\in \mathbb{N}$, respectively. Assuming that $y_k$ is possibly not measured in 
its full component form, i.e.~${y}_k\in \mathbb{R}^m, \textrm{ } m\leq d$, the goal is to 
estimate the state ${x}_n$ given noisy measurements 
\begin{equation}
\label{inv}
 {y}^{mes}=Y(x^{tru})+\hat{{\varepsilon}}
\end{equation}
in which $x^{tru}$ denotes the so-called truth, whereas $\hat{{\varepsilon}}$ stands for the corresponding
realisation of the measurement noise. 

Formally, in a Bayesian setting the unknown state ${x}_n$ in \refeq{inv} is modelled as a random variable (a priori knowledge or forecast)
\begin{equation}
\label{unstate}
 {x}_{nf}(\omega):\vOmega \rightarrow \mathbb{R}^d
\end{equation}
on a probability space $(\vOmega,\mathcal{F},\mathbb{P})$ endowed with the set of elementary events 
$\vOmega$, a $\sigma$-algebra of measurable events $\mathcal{F}$, and a probability measure $\mathbb{P}$. The common choice is to assume that $x_{nf}\in\mathcal{X}:=L_2(\vOmega,\mathcal{F},\mathbb{P};\mathbb{R}^d)$,
the space of real valued random variables with finite variance. 
As $x_n$ is a random variable, so is the observation in \refeq{det_observ}, here obtaining the form of
\begin{equation}
\label{pred_eq}
 \mathcal{Y}\ni {y}_{kf}(\omega)=Y({x}_{nf}(\omega))+\varepsilon_k(\omega)
\end{equation}
in which $\varepsilon_k(\omega)\sim \mathcal{N}(0,{C}_{\varepsilon_k})$ forecasts the measurement error 
usually taking the form of
zero-mean Gaussian noise with covariance ${C}_{\varepsilon_k}$.

Assuming that $x_n$ and $y_k$ have a joint probability density function $\pi(x_n,y_k)$, one may use 
Bayes's theorem in its density form
\begin{equation}
\label{bayes_thm}
 \pi_{x|y}(x_n|y_k)=\frac{\pi(x_n,y_k)}{P(y_k)}=\frac{\pi_{y|x}(y_k|x_n)\pi_x(x_n)}{P(y_k)}
\end{equation}
to incorporate (assimilate) new information $y^{mes}$ into the probabilistic description 
given in Eqs.~(\ref{unstate})-(\ref{pred_eq}).
Here, $\pi_x(x_n)$ denotes the prior density function, $\pi_{y|x}(y_k|x_n)$ is the likelihood,
the form of which depends on the measurement error, and $P(y_k)=\int_\vOmega \pi(x_n,y_k)dx_n $
is the normalisation factor or evidence. If both the prior and the likelihood are conjugate, i.e.~belong to the exponential family of
distributions with predefined statistics, the posterior $\pi_{x|y}(x_n|y_k)$ in \refeq{bayes_thm}
can be analytically evaluated.
 Otherwise, 
the estimation
boils down to computationally intense random walk algorithms of the Markov chain Monte Carlo type. However, both
computations essentially lead to the extraction of neccessary 
information from the posterior by evaluating some form of expectation w.r.t.~the posterior, an
example of which is the conditional mean
\begin{equation}
\label{mean}
 \mathbb{E}(x_n|y_k)=\int_{\vOmega} x_n \pi_{x|y}(x_n|y_k) \textrm{d}x_n.
\end{equation}
 Having done so, one may avoid expensive evaluation of the full posterior by targeting a direct calculation
 of desired estimates such as the one given in \refeq{mean}. 
 To achieve this, one may design filtering procedures based on conditional expectation as further described.
 
%

\section{Conditional expectation}\label{ss:cond}

The conditional expectation is defined as the unique optimal projector for all Bregman's loss functions (BLFs) \cite{bregman}
\begin{equation}
\label{eq:optimality_blf}
 x^*:=\mathbb{E}(x|\mathfrak{B})=\underset{\hat{x}\in L_2(\varOmega,\mathcal{B},\mathbb{P};\mathbb{R}^d)}{\arg \normalfont{\min}}
 \textrm{} \mathbb{E}(\mathcal{D}_\phi(x,\hat{x}))
\end{equation}
over all $\mathcal{B}$-measurable random variables $\hat{x}$ in which $\mathcal{B}:=\sigma(y)$ is the
sub-$\sigma$-algebra generated by measurement $y$.
The Bregman's loss function is defined as
\begin{definition}
 Let $\phi: \mathbb{R}^d \mapsto \mathbb{R}$ be a strictly convex, differentiable function. Then the Bregman loss function 
 $\mathcal{D}_\phi: \mathbb{R}^d\times \mathbb{R}
 \mapsto \mathbb{R}_+:=[0,+\infty)$ is defined as
 \begin{equation}\label{BLF}
  \mathcal{D}_\phi(x,y)=\mathcal{H}(x)-\mathcal{H}(y)=\phi(x)-\phi(y)-\langle x-y,\nabla \phi(y)\rangle
 \end{equation}
 in which $\mathcal{H}(x)=\phi(y)+\langle x-y,\nabla \phi(y)\rangle$ is hyperplane tangent to $\phi$ at point $y$. 
\end{definition}
The optimality in \refeq{eq:optimality_blf} then follows from \cite{banerjee}
 \begin{theorem}
  Let $\phi: \mathbb{R}^d \mapsto \mathbb{R}$ be a strictly convex, differentiable function and let $\mathcal{D}_\phi$
  be the corresponding BLF. Let $(\varOmega,\mathfrak{F},\mathbb{P})$ be an arbitrary probability space and let $\mathfrak{B}$
  be a sub-$\sigma$-algebra of $\mathfrak{F}$. Let $x$ be any $\mathfrak{F}$-measurable random variable taking values in $\mathbb{R}^d$
  for which both $\mathbb{E}(x)$ and $\mathbb{E}(\phi(x))$ are finite. Then, among all $\mathfrak{B}$-measurable random variables, the conditional
  expectation is the unique minimiser (up to a.s.~equivalence) of the expected Bregman loss, i.e.
  \begin{equation}
  x^*:=\mathbb{E}(x|\mathfrak{B})=\underset{\hat{x}\in L_2(\varOmega,\mathcal{B},\mathbb{P};\mathbb{R}^d)}{\arg \normalfont{\min}}
  \textrm{} \mathbb{E}(\mathcal{D}_\phi(x,\hat{x})).
\end{equation}
 \end{theorem}
The proof of the theorem can be shortly sketched as follows: 
\begin{proof}
 Let $\hat{x}$ be any $\mathfrak{B}$-measurable random variable, and $x^*=\mathbb{E}(x|\mathfrak{B})$, then one has 
 \begin{eqnarray}
 \label{proof_eq_1}
  && \mathbb{E}(\mathcal{D}_\phi(x,\hat{x}))-\mathbb{E}(\mathcal{D}_\phi(x,x^*))=\mathbb{E}(\phi(x^*)-\phi(\hat{x})\nonumber\\
  &&-
  \langle x-\hat{x},\nabla \phi(\hat{x})\rangle
  +\langle x-x^*,\nabla \phi(x^*)\rangle. 
 \end{eqnarray}
Using the law of total expectation, e.g. $\mathbb{E}(x)=\mathbb{E}(\mathbb{E}(x|\mathfrak{B}))$, one may further state
\begin{eqnarray}
 \mathbb{E}(\langle x-\hat{x},\nabla \phi(\hat{x})\rangle)&=&\mathbb{E}\left(\mathbb{E}(\langle x-\hat{x},\nabla \phi(\hat{x})\rangle
 |\mathfrak{B}\right))\nonumber \\
 &=&
\mathbb{E}(\langle \mathbb{E}(x|\mathfrak{B})-\hat{x},\nabla \phi(\hat{x})\rangle)\nonumber\\
&=&\mathbb{E}(x^*-\hat{x},\nabla \phi(\hat{x}))
\end{eqnarray}
Similarly,
\begin{eqnarray}
 \mathbb{E}(\langle x-x^*,\nabla \phi(\hat{x})\rangle)
 &=&\mathbb{E}(\mathbb{E}(\langle x-x^*,\nabla \phi(\hat{x})\rangle|\mathfrak{B}))\nonumber\\
 &=&
\mathbb{E}(x^*-x^*,\nabla \phi(\hat{x}))\nonumber\\
&\equiv&0.
\end{eqnarray}
Following this, the relation in \refeq{proof_eq_1} reduces to 
\begin{eqnarray}\label{breg_in}
&& \mathbb{E}(\mathcal{D}_\phi(x,\hat{x}))-\mathbb{E}(\mathcal{D}_\phi(x,x^*))=\nonumber \\
&&\mathbb{E}(\phi(x^*)-
 \phi(\hat{x})-\langle x^*-\hat{x},\nabla \phi(\hat{x})\rangle)\nonumber\\
 &&=\mathbb{E}(\mathcal{D}_\phi(x^*,\hat{x})).
\end{eqnarray}
\end{proof}
%
The last relation in \refeq{breg_in} defines the Bregman Pythagorean inequality 
\begin{equation}
\label{in_breg}
  \mathbb{E}(\mathcal{D}_\phi(x,\hat{x}))\geq \mathbb{E}(\mathcal{D}_\phi(x,x^*))+\mathbb{E}(\mathcal{D}_\phi(x^*,\hat{x}))
\end{equation}
such that one may state
\begin{theorem}
Let $\phi: \mathbb{R}^d \mapsto \mathbb{R}$ be a strictly convex, differentiable function and let $\mathcal{D}_\phi$
  be the corresponding BLF. Let $(\varOmega,\mathfrak{F},\mathbb{P})$ be an arbitrary probability space and let $\mathfrak{B}$
  be a sub-$\sigma$-algebra of $\mathfrak{F}$. Let $\hat{x}$ and $x$ be any $\mathfrak{F}$-measurable random variable taking values in $\mathbb{R}^d$
  for which both pairs $(\mathbb{E}(\hat{x}),\mathbb{E}(x))$ and $(\mathbb{E}(\phi(\hat{x})),\mathbb{E}(\phi(x)))$ are finite. Then, we have
  \begin{equation}
  \mathbb{E}(\mathcal{D}_\phi(x,\hat{x}))\geq \mathbb{E}(\mathcal{D}_\phi(x,x^*))+\mathbb{E}(\mathcal{D}_\phi(x^*,\hat{x}))
\end{equation}
in which the unique point $x^*$ is called the Bayesian projection of $x$ onto $\mathfrak{B}$ and is defined as following
\begin{equation}
  x^*:=\mathbb{E}(x|\mathfrak{B})=P_{\mathfrak{B}} x=
  \underset{\hat{x}\in L_2(\varOmega,\mathcal{B},\mathbb{P};\mathbb{R}^d)}{\arg \normalfont{\min}} \textrm{} \mathbb{E}(\mathcal{D}_\phi(x,\hat{x}))
\end{equation}
\end{theorem}
 
%
%
Note that if we took $x^*=\mathbb{E}(x)$ then the term $ \mathbb{E}(\mathcal{D}_\phi(x,x^*))$ is known as the Bregman's variance 
\begin{equation}
   \textrm{var}_\phi(x)=\mathbb{E}(\mathcal{D}_\phi(x||\mathbb{E}(x)))=\mathbb{E}(\phi(x))-\phi(\mathbb{E}(x))\geq 0
  \end{equation}
  for which holds (see \cite{banerjee})
\begin{theorem}\label{thmber}
 Let $x$ be a random variable with mean $\mathbb{E}(x)$ and variance $\textrm{var}(x)$. The Bregman variance 
 $\textrm{var}_\phi(x)\neq \textrm{var}(x)$
 is then defined as follows
 \begin{eqnarray}
 \label{eqbt1}
  \textrm{var}_\phi(x)&=&\mathbb{E}(\mathcal{D}_\phi(x||\mathbb{E}(x)))\nonumber\\
  &=&\mathbb{E}(\phi(x))-\phi(\mathbb{E}(x))\geq 0. 
 \end{eqnarray}
From inequality \refeq{breg_in} one may further state
\begin{eqnarray}
  \textrm{var}_\phi(x)&=&\mathbb{E}(\mathcal{D}_\phi(x||\mathbb{E}(x)))\nonumber\\
  &=&\mathbb{E}(\mathcal{D}_\phi(x,\hat{x}))
  -\mathbb{E}(\mathcal{D}_\phi(\mathbb{E}(x),\hat{x}))\nonumber \\
  &\geq & 0
  \end{eqnarray} for any random variable $\hat{x}$. This then leads to
  \begin{equation}
   \mathbb{E}(x)=\underset{\hat{x}\in L_2(\varOmega,\mathcal{F},\mathbb{P};\mathbb{R}^d)}\arg \normalfont{\min }\textrm{ } \mathbb{E}(\mathcal{D}_\phi(x,\hat{x}))
  \end{equation}
  which is the same minimum point for any expected Bregman's divergence.
\end{theorem}
The key result of the previous theorems justifies using a mean as
a representative of a random variable, particularly in a Bayesian estimation.

In a special case when $\phi$ takes the quadratic form, i.e. $\phi(x)=\frac{1}{2}\|x\|_{L_2}^2$, 
the Bregman's divergence in \refeq{BLF} modifies to the squared-Euclidean distance
\begin{equation}
 \mathcal{D}_\phi(x||y)=\|x-y\|^2.
\end{equation}
In such a case the Bregman Pythagorean theorem \refeq{in_breg} reduces to the classical 
Pythagorean theorem as already discussed by the author and co-workers in \cite{Matthies2016}. 

Following the authors previous works, the conditional expectation ${E}({x}|{y})$ of a random variable $x$ given the measurement $y$
 in terms of Bregman's quadratic loss functions is an orthogonal projection $P_\mathcal{B}({x})$
of ${x}$ onto the subspace $L_2(\vOmega,\mathcal{B},\mathbb{P};\mathbb{R}^d)$ of all random variables consistent with 
the data ${y}$, i.e.~generated by 
the sub-sigma algebra $\mathcal{B}:=\sigma({y})$. This further means that ${x}$ can be orthogonally 
decomposed into two components $x_p$ and $x_o$: 
\begin{equation}
\label{ff}
 {x}={x}_p+{x}_o
\end{equation}
in which the projected part reads ${x}_p:=P_\mathcal{B}({x})$, whereas the orthogonal component $x_o$
equals $(I-P_\mathcal{B}){x}$. 

As an observation $y^{mes}$ arrives, the first term in \refeq{ff}, $x_p$, 
 is altered by the data $y^{mes}$, whereas the latter one, ${x}_o$,
embodies the remaining (residuals) of the prior information $x_f$. This idea leads to the analogy of 
${x}_p$ with $\mathbb{E}({x}_f|{y}^{mes})$ and of ${x}_o$ with ${x}_f(\omega)-\mathbb{E}({x}_f|{y}_f)$ in which $y_f$ takes the form given in
\refeq{pred_eq} such that 
\begin{equation}
\label{ffa}
 {x}_a=\mathbb{E}({x}_f|{y}^{mes})+({x}_f-\mathbb{E}({x}_f|{y}_f))
\end{equation}
holds. This is the filtering form of the decomposition given in \refeq{ff}, in which the indices $a$ and $f$ are used to denote the assimilated (posterior)
state and forecast (prior) state, respectively.   
 Following the Doob-Dynkin lemma, the previous equation can be rewritten as 
\begin{equation}
\label{form1}
 {x}_a=\varphi({y}^{mes})+({x}_f-\varphi({y}_f)),
\end{equation}
in which the conditional expectation $\mathbb{E}({x}_f|{y}^{mes})$ is represented by a measurable map $\varphi({y}^{mes})$, and similarly $\mathbb{E}({x}_f|{y}_f)$
is expressed as $\varphi({y}_f)$.
By rearranging the terms in \refeq{form1} one obtains
\begin{equation}
\label{fil}
 {x}_a={x}_f+\varphi({y}^{mes})-\varphi({y}_f),
\end{equation}
the general form that is further used to construct the nonlinear filtering procedure. The advantage of
\refeq{fil} compared to \refeq{bayes_thm} is that all quantities of consideration 
are given in terms of random variables, and not probability measures. Hence, it is easier to functionally approximate and computationally manipulate 
\refeq{fil} than \refeq{bayes_thm}, as further discussed.

%
%
%
%
%
%
%
\section{Optimal map}\label{opmap}
To obtain the maximal information gain in \refeq{fil}, 
the task is to find the optimal map $\varphi$ among all measurable maps $\mathcal{Y}\rightarrow \mathcal{X}$. However, this step is not computationally tractable, and thus 
additional approximations are required. The simplest possible choice is to consider a linear approximation
\begin{equation}
\label{cond_map}
\mathbb{E}(x_f|y_f)\approx K{y_f}+{b}
\end{equation}
in which the map coefficients $(K,b)$ are obtained 
by minimising the orthogonal component in \refeq{ffa}, i.e.~
\begin{eqnarray}
\label{inver_map_k}
 &&\underset{K,b}{\arg \min}\textrm{ } \mathbb{E}(\|x_f-\mathbb{E}(x_f|y_f)\|_2^2)\nonumber\\
 &&=\underset{K,b}{\arg \min}\textrm{ } \mathbb{E}(\|x_f-(Ky_f+b)\|_2^2).
\end{eqnarray}
From the optimality condition
\begin{equation}
 \forall \chi:\quad \mathbb{E}(\langle x_f-(Ky_f+b),\chi \rangle) =0.
\end{equation}
 one obtains
\begin{eqnarray}
 &&\mathbb{E}(\langle x_f-K{y}_f-b,y_f\rangle)=0\nonumber\\
 &&\mathbb{E}\left(x_f-K y_f-b\right)=0
\end{eqnarray}
which further results in a linear Gauss-Markov-Kalman filter equation
\begin{equation}
\label{line_gaus_mar}
{x}_a(\omega)={x}_f(\omega)+{K}({y}^{mes}-{y}_f(\omega)),
\end{equation}
 specified by the well-known Kalman gain 
\begin{equation}
\label{kalmna_gain_eq}
K={C}_{{x}_f,{y}_f}({C}_{{y}_f})^{\dagger}.
\end{equation}
Here, $\dagger$ denotes the pseudo-inverse, ${C}_{{x}_f,{y}_f}$ is the covariance
between the prior $x_f$ and the observation forecast $y_f$, and 
${C}_{{y}_f}={C}_{Y({x}_f)}+C_\varepsilon$ is the auto-covariance of $y_f$ consisting of
forecast covariance ${C}_{Y({x}_f)}$ and the measurement covariance 
$C_\varepsilon$. 

Even though computationally cheap, the previous formula uses only pieces of provided information in ${y}^{mes}$ and may lead to 
over-- or under-- estimation 
in highly nonlinear systems. Namely, the term $y_f$ in \refeq{line_gaus_mar} is essentially nonlinear and 
does not comply with the linear approximation of the map $\varphi$. To resolve nonlinearity,
let the measurement operator $Y$ be the Fr\'{e}chet differentiable
with Lipschitz continuous derivative $H:=\partial Y/\partial x$ such that
 \begin{equation}
\label{mesoperequiv}
Y({x}) \approx Y(\check{{x}})+H({x}-\check{{x}})=:Y_\ell({x})
\end{equation}
holds. Following this assumption, one may further state
\begin{equation}
\label{mesr_equiv}
 y\approx Y_\ell(x)+\varepsilon=:y_\ell(x)
 \end{equation}
 in which $y_\ell(x)$ represents the linearised measurement around the point $\check{x}$. 
 As $Y_\ell$ is linear, the new Gauss-Markov-Kalman formula obtains a similar form to the one given in \refeq{line_gaus_mar} and reads
 \begin{eqnarray}
 \label{linearised_linear_filter}
 {x}_a&=&{x}_f+{K}_\ell({y}^{mes}-{y}_\ell(x_f))\\
 &=&x_f+K_\ell(y^{mes}-Y(\check{x})-H(x_f-\check{x})-\varepsilon)\nonumber.
\end{eqnarray}
 Here, $x_f$ is the forecast parameter, $y_{\ell}(x_f)$ is the forecast value of linearised 
 measurement around the point $\check{x}$ given the prior $x_f$, $\varepsilon$ is the model of
 the measurement error,
 and $K_\ell$ is the corresponding Kalman gain calculated via
 \begin{eqnarray}
  \label{kalman_gain_linearised}
 K_\ell&=&C_{x_fy_{\ell}}(C_{y_{\ell}})^{\dagger}\nonumber\\
 &=&C_{x_f}H^T(HC_{x_f}H^T+C_\varepsilon)^{\dagger}.
 \end{eqnarray}
 Note that in a special case when all distributions of consideration are known to be Gaussian, the last formula obtains
 a similar form to the extended Kalman filter \cite{einicke2012smoothing,simon2006optimal}.

The map in \refeq{linearised_linear_filter} is not optimal as it highly depends on the choice
of the point $\check{x}$. Obviously, $\check{x}$ taken as $\mathbb{E}(x_f)$ is not always the best choice.
To find an optimal linearisation point, one may introduce 
the sequence of the first order approximants
\begin{equation}
 \label{mesoperequiv_seq}
 Y_\ell^{(i)}(x):= Y(\check{x}^{(i)})+H^{(i)}(x-\check{x}^{(i)})
 \end{equation}
with
\begin{equation}
 \label{mesr_equiv_seq}
 H^{(i)}:=\frac{\partial Y}{\partial x}\Big |_{\check{x}^{(i)}},
\end{equation}
and
\begin{equation}
  y_\ell^{(i)}(x)=Y_\ell^{(i)}(x)+\varepsilon.
  \end{equation}
 As a result, the optimal map $Y$ is iteratively found via the sequence of Kalman gains
 \begin{eqnarray}
  \label{kalman_gain_seq}
 K_\ell^{(i)}&=&C_{x_fy_\ell^{(i)}}C_{y_\ell^{(i)}}^{\dagger}\\
 &=& C_{x_f}(H^{(i)})^T(H^{(i)}C_{x_f}(H^{(i)})^T+C_\varepsilon)^{\dagger},\nonumber
 \end{eqnarray}
 and subsequently the posterior state is estimated via an iterative procedure
 \begin{equation}
 \label{newton_gauss_markov}
  x_a^{(i+1)}=x_f+K_\ell^{(i)}(y^{mes}-y_\ell^{(i)}(x_f)),
 \end{equation}
 \begin{equation}
  \check{x}^{(i)}=\mathbb{E}(x_a^{(i+1)}),
 \end{equation}
here called the Gauss-Newton-Markov-Kalman filter. Under Gaussianity assumptions 
one may show that the previous equation represents the Gauss-Newton procedure for the maximum aposteriori
estimate (MAP) as shown in
 \cite{bell}. Note that no such assumption is made here.
 
The convergence properties of the algorithm can be studied via fixed point theorem \cite{gratton},
according to which the algorithm has local convergence characterised by a spectral radius of $\rho(K_\ell^{(i)}{H}^{(i)})$.

%

 \section{Bayesian estimation of optimal map}\label{ss:unbiased_filter}
 
In the form given in \refeq{newton_gauss_markov} the Gauss-Newton-Kalman filter has two drawbacks: first the 
filter requires the time consuming evaluation of
the Jacobian ${H}^{(i)}$, and second
the filter is biased 
as it assumes that
 \begin{equation}
 \label{eqal_mean_mome}
  \mathbb{E}\left[\left(Y(x_f)\right)^k\right]=\mathbb{E}\left[\left(Y(\check{x})+H(\check{x})(x_f-\check{x})\right)^k\right]
 \end{equation}
 holds for $k=1,..,n$. Therefore, the straightforward linearisation is not the best possible choice.
 Instead, one may search for the optimal linear map in a similar setting as given in Section \ref{opmap}.  
 
 In numerical practice the measurement operator $Y(x_f)$ is encoded in the 
 corresponding computer software/simulator of the physical model, and hence is not explicitly known. But, using the classical 
 uncertainty quantification procedures (e.g.~the pseudo-spectral method or similar) one may obtain $z_f:=Y(x_f)$
 given $x_f$ in a non-intrusive way. In such a case both $z_f$ and $x_f$ are known, and hence the estimation 
 of the measurement operator $Y(x_f)$ in a linearised form becomes simple. It only requires an estimation of the map $\varphi_y: x_f \mapsto z_f$, i.e.
 the conditional expectation $\mathbb{E}(z_f|x_f)$. By taking the Bregman's squared loss function, as already 
 discussed, the parameterised map $\varphi_y(\beta)$ 
 can be estimated 
 by minimising 
 \begin{equation}
 \label{mmse_inverse1}
\beta^*= \underset{\beta}{\arg \min}\textrm{ }\mathbb{E}(\|z_f-\varphi_y(x_f,\beta)\|_2^2).
\end{equation}
%
%
In a special affine case
\begin{equation}
  \mathbb{E}(z_f|x_f)\approx \check{H}(x_f-\check{x})+h=:\varphi_y(x_f,\beta),
\end{equation}
with $\beta:=(\check{H},h)$ the previous optimisation problem reduces to
 \begin{equation}
 \label{mmse_inverse}
 \underset{\check{H},h}{\arg \min}\textrm{ }\mathbb{E}(\|z_f-(\check{H}(x_f-\check{x})+h)\|_2^2),
\end{equation}
%
%
%
%
%
 the solution of which 
 \begin{equation}
 \label{chechkhapp}
 \check{H}=C_{Y(x_f),x_f}C_{x_f}^\dagger
 \end{equation}
 represents the approximation of the Jacobian, and
 \begin{equation}
  h:=\mathbb{E}(Y(x_f))-\mathbb{E}(x_f-\check{x})
 \end{equation}
 is the linear constant. Note that if $Y(x_f)$ is originally linear described by the true Jacobian $H$, then using the formula in \refeq{chechkhapp} 
 one has that 
 \begin{equation}
 \label{eq1top}
  \check{H}=C_{Y(x_f),x_f}C_{x_f}^\dagger=HC_{x_f}C_{x_f}^\dagger\equiv H.
 \end{equation}
 Similarly, for the inverse map $z_f \mapsto x_f$ holds
 \begin{equation}
 \label{eq2top}
  C_{x_f,z_f}C_{z_f}^\dagger\equiv H^\dagger.
 \end{equation}
 Employing the previous two relations one may conclude that the Jacobian of the forward map is equal to the inverse Kalman gain
 when the observation $y_f=z_f+\varepsilon$ does not contain the measurement/modelling/approximation error $\varepsilon$.

 However, note that \refeq{mmse_inverse} holds only if linearisation is done once as in the extended Kalman filter procedure. Otherwise, 
 given $z_a^{(i)}:=Y(x_a^{(i)})$ one solves the following problem
 \begin{equation}
 \label{mmse_inverse_iter}
 \underset{\check{H}^{(i)},h^{(i)}}{\arg \min}\textrm{ }\mathbb{E}(\|z_a^{(i)}-(\check{H}^{(i)}(x_a^{(i)}-\check{x}^{(i)})+h^{(i)})\|_2^2),
\end{equation}
%
%
such that the
 filter in \refeq{newton_gauss_markov} obtains its unbiased form
 \begin{eqnarray}
 \label{newton_gauss_markov1}
  x_a^{(i+1)}&=&x_f+K_\ell^{(i)}(y^{mes}-y_h^{(i)}(x_f))
 \end{eqnarray}
 in which
 \begin{eqnarray}
  y_h^{(i)}(x_f):&=&\check{H}^{(i)}(x_f-\check{x}^{(i)})+h^{(i)}+\varepsilon.
 \end{eqnarray}
 
 Note that previously we have assumed that we know random variables $z_f$ and $x_f$ resp. $z_a^{(i)},x_a^{(i)}$, which is 
 often not the case. Instead, in numerical simulations we may only know their samples. Let us denote the set of 
  samples of the variable $x_a^{(i)}$ by $x^{sim}:=(x_a^{(i)}(\omega_i)_{i=1}^N)$. Similarly, let us denote the set of forecasted samples
 by $z^{sim}:=(z_a^{(i)}(\omega_i)_{i=1}^N)$
 such that 
 \begin{equation}
 \label{eqeq}
  z_a^{(i)}(\omega_i)=\varphi_y(x_a^{(i)}(\omega_i))+\varepsilon_y(\omega_i)
 \end{equation}
holds. 
%
In such a case, the approximation of the Jacobian 
can be estimated from
 \begin{equation}
\label{optmaplin}
  z_a^{(i)}(\omega_i)=\check{H}^{(i)}(x_a^{(i)}\omega_i)-\check{x}^{(i)})+h^{(i)}+\varepsilon_y(\omega_i)
 \end{equation}
in a Bayesian framework given measurement data $d^{sim}:=(x^{sim},z^{sim})$ by assuming that the pair $h,\check{H}$ and
the approximation error $\varepsilon_y$ 
are unknown, and hence modelled as uncertain. In a Bayesian setting the map parameters $\beta:=(\check{H},h,\varepsilon_y)$ can be estimated as:
\begin{equation}
\label{bayes_mapp}
 \pi_{\beta|d^{sim}}(\beta|d^{sim})\propto\pi_{d^{sim}|\beta}
 (d^{sim}|\beta)
 \pi_{ \beta}(\beta)
\end{equation}
in which $\pi_{\beta}(\beta)$ is a joint prior distribution on $\beta$ 
 here factorised according to
$\pi_{\beta}=\pi_{\check{H}}(\check{H})\pi_h(h)\pi_{\varepsilon_y}(\varepsilon_y)$.
The prior information can be imposed further such that each element of the prior is of Gaussian type. 
As \refeq{optmaplin} is of linear type, the Bayesian estimation in such a case reduces to the Kalman filter estimate. 
For this purpose one may assume that the prior mean for the Jacobian is close to the inverse of the 
previously estimated Kalman gain, see \refeq{eq1top} and \refeq{eq2top}. 
To include more information into the prior such as sparsity of the matrix, the prior has to be carefully designed, as  
discussed in Section \ref{spaopt}. 

 
Note that same type of approach can be also used for the estimation of the Kalman gain in \refeq{kalmna_gain_eq}. Following \refeq{eqeq} one may pose 
the following problem: given samples $(x_f(\omega_i),y_h^{(i)}(x_f(\omega_i))$ estimate $\mathbb{E}(x_f|y_f)=\varphi(y_h^{(i)})$  such that
\begin{equation}
 x_f=\mathbb{E}(x_f|y_f)+\varepsilon_x=\varphi(y_h^{(i)})+\varepsilon_x^{(i)}
\end{equation}
holds.
Assuming linear map
\begin{equation}
 \varphi_y(y_h^{(i)})=K^{(i)} y_h^{(i)}+b
\end{equation}
and given the data set $d^{sim}:=(x_f(\omega_i),y_h^{(i)}(\omega_i))$ one may use Bayes's rule to estimate 
$\beta:=(K,b,\varepsilon_x)$ in a similar manner as in \refeq{bayes_mapp}. The numerical advantage of
Bayes's rule compared to \refeq{inver_map_k} lies in the prior knowledge which can be 
imposed on the Kalman gain, e.g.~the sparsity information on the mapping coefficients as discussed in Section \ref{spaopt}.
This further allow us
to use the previously described filter in a "hierarchical sense"
for both solving the inverse problem, as well as for estimating the optimal linear map. 
In particular, the hierarchical approach is interesting 
when one would like to estimate the approximation/modelling/linearisation error $\epsilon$ as further discussed in Section \ref{spaopt}.
However, note that by using Bayes's rule to obtain a Kalman gain we do not satisfy the orthogonality condition, and hence we do not have 
a Kalman filter estimate as understood in the classical sense.

 \section{Predictor-corrector Bayesian-Gauss-Newton-Markov-Kalman filter for backpropagation}\label{pcp}
 
\begin{figure}
\begin{center}
\includegraphics[width=0.35\textwidth]{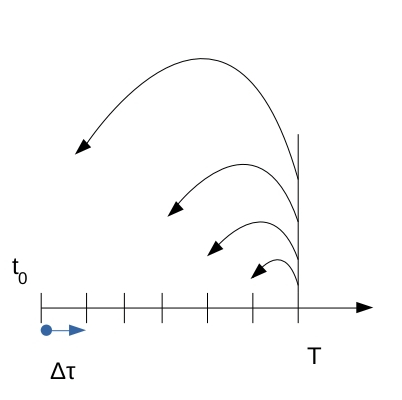}
\caption{Schematic representation of direct backward propagation}
\label{fig_ds_pic}
\end{center}
\end{figure}

\begin{figure*}
\begin{center}
\includegraphics{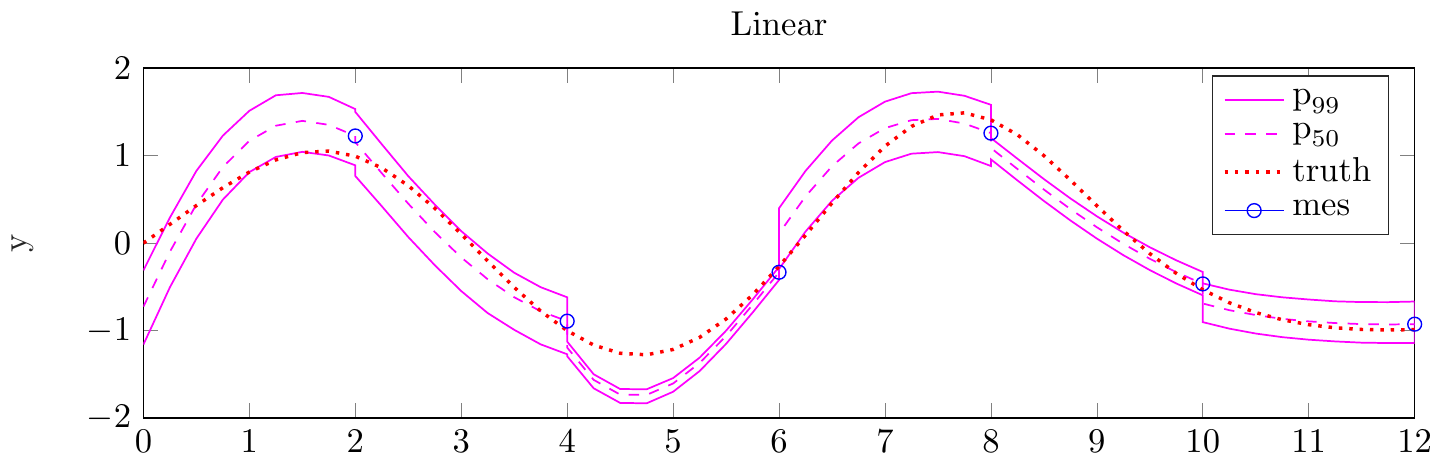}\\
\includegraphics{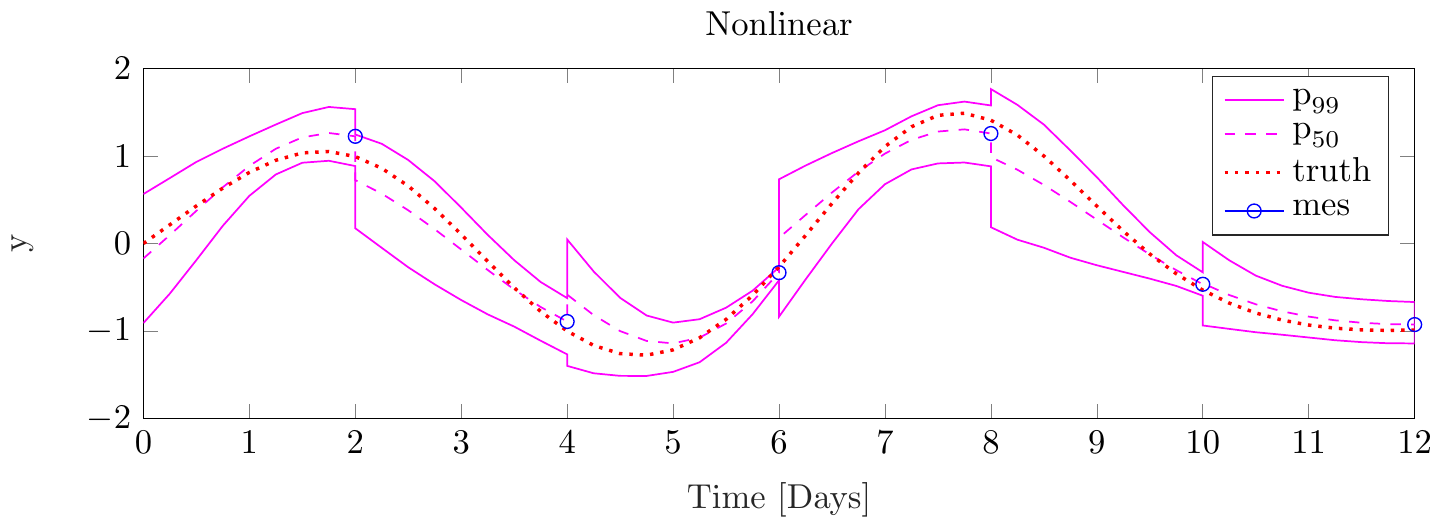}\\
\caption{Linear and nonlinear smoothing of the second component of the Lorenz 1984 system}
\label{fig_linear_nonlinear_update}
\end{center}
\end{figure*}

\begin{figure*}
\begin{center}
 \includegraphics{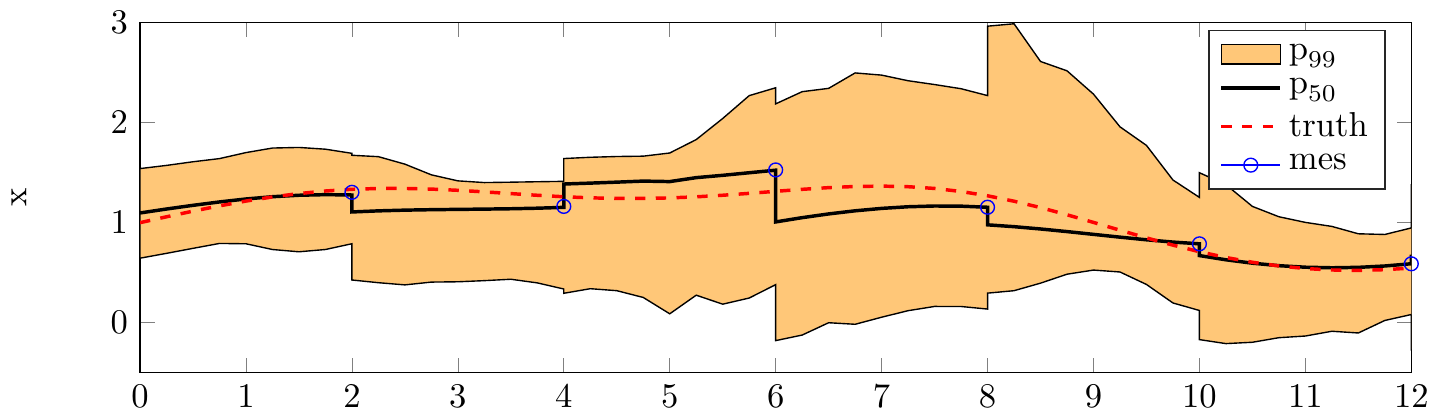}\\
\includegraphics{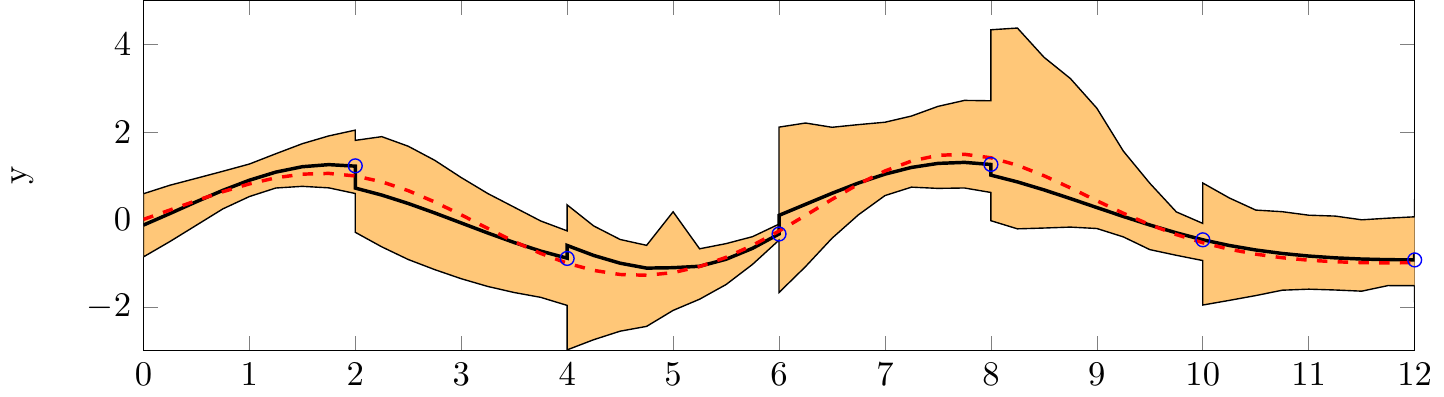}\\
 \includegraphics{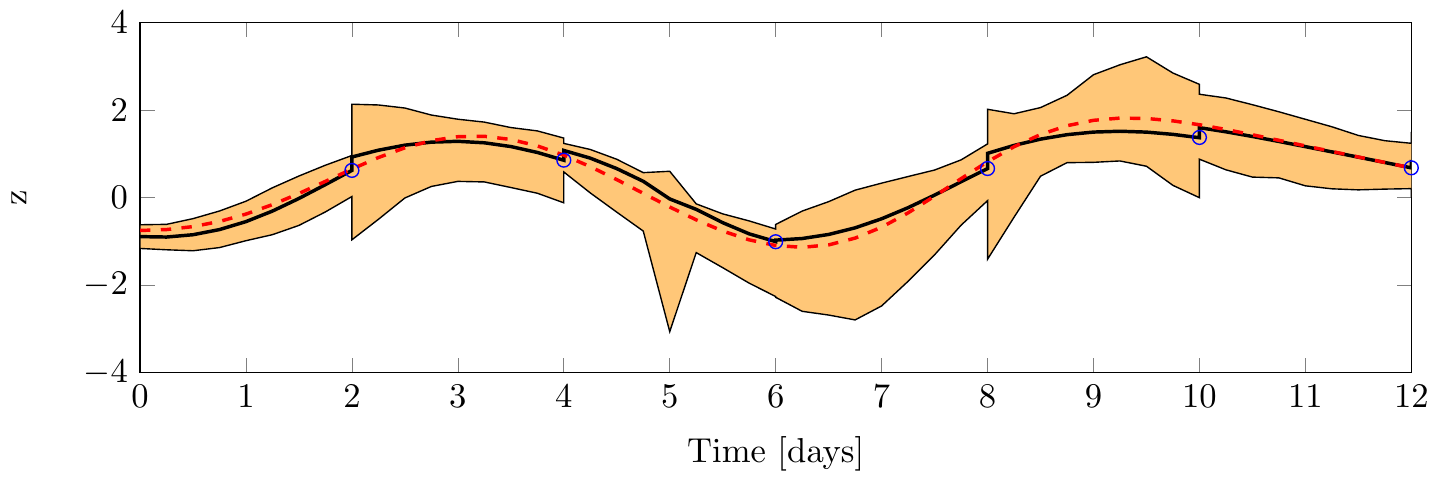}
\caption{The backward estimation of the Lorenz 1984 state with an updating window size of two days}
\label{fig_filter_non_seq}
\end{center}
\end{figure*}

To estimate the initial condition of the dynamical system given in \refeq{initvp}, one may use the previously designed filter in the following form:
 \begin{equation}
 \label{direct_estimate}
  x_{0,a}^{(i+1)}={x_{0,f}}+K_\ell^{(i)}(y^{mes}-y_{\ell}^{(i)}(x_{0,f})),
 \end{equation}
 in which $x_{0,f}$ is the a priori random variable describing the initial condition at $t_0$, $y^{mes}$ is the measurement 
 at the time $T$ and $y_{\ell}:=\check{H}^{(i)}(x_f-\check{x}^{(i)})+\check{h}^{(i)}+\varepsilon$ is the forecasted
 linearised measurement at $T$ and in iteration $(i)$. In a similar manner one may also estimate any state  between $t_0$ and $T$.
 Considering the identification of all states equidistantly separated by the update time step $\Delta \tau$, the Gauss-Newton-Markov-Kalman filter 
 is schematically described in \refalg{alg1}-\refalg{alg11}, and 
 depicted in \refig{fig_ds_pic}. After initialisation of the prior variable, one approximates the forward map $x_f \mapsto y_f$ by the linearised operator estimated 
 either in a classical way, see \refeq{mmse_inverse_iter}, or in a Bayesian manner, see \refeq{bayes_mapp}. Once the linearised measurement is found, one may estimate the inverse map 
 linearly again in two different manners: by projection or by Bayes's rule. Once both maps are estimated one may assimilate the state 
 using the measurement data, and hence update the linearisation point. This method of estimating the state will be called direct smoothing (DS) further on. 
 A numerical example is shown in
\refig{fig_linear_nonlinear_update}. Here, the smoothing algorithm
with a window size of two days is used to estimate the second component of the Lorenz 1984 system (for
model details see the Appendix) given 
noisy full state measurement data, see \refalg{alg11}. Clearly, the linear update observed in the upper plot
fails to properly estimate the state in any other time moment than the time of the measurement itself. 
On the other hand, the nonlinear filtering counterpart taking the iterative form as described
previously produces satisfying results, see the
lower plot in \refig{fig_linear_nonlinear_update}. This also holds for all three Lorenz components as depicted in \refig{fig_filter_non_seq}.

\begin{figure}[htpb]
\begin{center}
\includegraphics{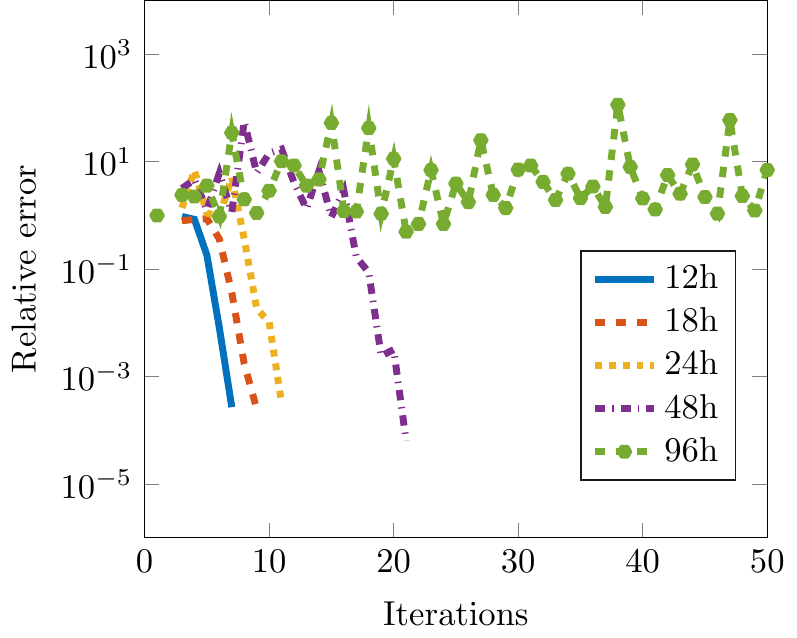}
\caption{Convergence of posterior estimate of the initial condition $\vek{x}_0$ w.r.t.~time at which the measurement data arrive}
\label{fig_firsta}
\end{center}
\end{figure}

 In general the Gauss-Newton procedure is known to be convergent when the residuals are assumed to be small. 
  However, if the state $x_n$ is estimated given $y_k$, in which $k$ is many times
 larger than $n$ (e.g.~estimation of the initial state after long time integration), 
 and/or the system 
 is highly nonlinear, the direct estimation can be a problem. Fig.~(\ref{fig_firsta}) depicts an example of filter divergence when estimating 
  the initial condition of the Lorenz 1984 system given the state measured after 96 hours. To overcome this, the large ``update step'', 
  i.e.~the time interval $[t_n,t_k]$,
 is split into smaller update steps defined by pseudo-time moments 
$t_n\leq \tau_\ell \leq t_k, \tau_\ell=t_n+\ell \Delta \tau $ via $\Delta \tau=c\Delta t$ stepping 
in which $\Delta t$ is the time discretisation step, and $1\leq c\in \mathbb{N}$. In this way one divergent Gauss-Newton iteration is 
substituted by several convergent ones, and the direct estimation is substituted by an incremental one.  

%

\begin{algorithm*}
   \caption{\textbf{Direct smoothing (DS)}: Bayesian-Gauss-Newton-Markov-Kalman filter, backpropagation}
   \label{alg1}
    \begin{algorithmic}[1]
      \Function{BGNMK}{$x_{0,f},y^{mes},@integ, \Delta t, \Delta \tau,t_0,T,\varepsilon$}\Comment{Where $x_{0,f}$ - 
      prior on initial value, $t_0$ beginning of time interval, $T$ end of time interval,
      $y^{mes}$ - measurements, $integ$ - forward function (model) handle, $\Delta t$ - time integration step, $\Delta \tau$ - 
      update step, $\varepsilon$ - measurement error}
      \State
         \For{$tt=t_0:\Delta \tau:T$} \Comment{Update the assimilation time}
          \State  \textbf{Set prior}
         \State $\quad x_{f}=integ(x_{0,f},\Delta t,t_0,tt)$\Comment{integrate ODE system from $t_0$ to $tt$ by time step $\Delta t$}
         \State  \textbf{Update}
         \State $\quad x_a=\textrm{GNMK}(x_{f},y^{mes},@integ, \Delta t,tt,T,\varepsilon)$
      \EndFor
%
       \EndFunction

\end{algorithmic}
\end{algorithm*}

\begin{algorithm*}
   \caption{\textbf{Direct smoothing (DS)}: Bayesian-Gauss-Newton-Markov-Kalman filter, backpropagation}
   \label{alg11}
    \begin{algorithmic}[1]
      \Function{GNMK}{$x_{f},y^{mes},@integ, \Delta t, \Delta \tau,tt,T,\varepsilon$}\Comment{Where $x_{0,f}$ - 
      prior on initial value, $t_0$ beginning of time interval, $T$ end of time interval,
      $y^{mes}$ - measurements, $integ$ - forward function (model) handle, $\Delta t$ - time integration step, $\Delta \tau$ - 
      update step, $\varepsilon$ - measurement error}
      \State
          \State  \textbf{Set linearisation point}
        \State $\quad \mathring{x}^{(0)}=\mathbb{E}(x_f), \quad x_a^{(0)}=x_f$
          \State \textbf{Set} $i=0$, $\textrm{err}=2\cdot\textrm{tol}$, $\textrm{maxiter}=100$
        \While{$i\leq \textrm{maxiter}$ \& $\textrm{err}\leq \textrm{tol}$} 
          \State  \textbf{Predict measurement}
         \State $\quad z_{a}^{(0)}=integ(x_{a}^{(i)},\Delta t,tt,T)$\Comment{integrate ODE system from $tt$ to $T$}
          \State \textbf{Approximate forward map $x_a^{(i)} \mapsto z_a^{(i)}:=Y(x_a^{(i)})$ by}
          \State $\quad$ $\varphi_y(x_a^{(i)})=\mathring{H}^{(i)}(x_a^{(i)}-\mathring{x}^{(i)})+\mathring{h}^{(i)}$
            \State \textbf{Estimate forward map coefficients $\beta:=(\mathring{H}^{(i)},\mathring{h}^{(i)})$ by}
        \State $\quad$ - projection:
        \State $\quad\quad\quad$   $\mathring{H}^{(i)}=C_{z_{a}^{(i)},x_{a}^{(i)}}C_{x_{a}^{(i)}}^\dagger,\quad 
        \mathring{h}^{(i)}=\mathbb{E}(z_{a}^{(i)})-\mathring{H}^{(i)}(x_{a}^{(i)}-
        \mathring{x}^{(i)}),$
        \State $\quad$ - or by Bayes's rule 
        \State $\quad\quad\quad$ given data $d^{sim}=(x_a(\omega_j)^{(i)},z_a(\omega_j)^{(i)}), j=1,...,N$ (see Section \ref{spaopt})
        \State $\quad\quad\quad$ update $\pi_{\beta|d^{sim}}(\beta|d^{sim})\propto\pi_{d^{sim}|\beta}
 (d^{sim}|\beta)
 \pi_{\beta}(\beta)$
        \State \textbf{Linearise predicted measurement}
        \State $\quad y_{\ell}^{(i)}({x_{f}})=\mathring{H}^{(i)}(x_{f}-\mathring{x}^{(i)})+\mathring{h}^{(i)}+\varepsilon$
        \State \textbf{Approximate inverse map $y_\ell^{(i)} \mapsto x_f$ by}
          \State $\quad$ $\varphi(y_\ell^{(i)})=K^{(i)}y_\ell^{(i)}+b^{(i)}$
       \State \textbf{Estimate inverse map coefficients $w:=(K^{(i)},b^{(i)})$ by}
        \State $\quad$ - projection:
        \State $\quad\quad\quad$   $K^{(i)}=C_{x_{f},y_{\ell}^{(i)}}C_{y_{\ell}^{(i)}}^\dagger, \quad b=\mathbb{E}(x_f)-K^{(i)}\mathbb{E}(y_\ell^{(i)})$ 
        \State $\quad$ - or by Bayes's rule 
        \State $\quad\quad\quad$ given data $d^{sim}=(x_f(\omega_j),y_{\ell}^{(i)}(\omega_{(j)}), j=1,...,N$ (see Section \ref{spaopt})
        \State $\quad\quad\quad$ update $\pi_{w|d^{sim}}(w|d^{sim})\propto\pi_{d^{sim}|w}
 (d^{sim}|w)
 \pi_{w}(w)$
       \State \textbf{Update state}
       \State $\quad$ $ x_{a}^{(i+1)}={x_{f}}+K^{(i)}(y^{mes}-y_{\ell}^{(i)}),\quad $
       \State \textbf{Update linearisation point}
       \State $\quad i=i+1;$
       \State $\quad \mathring{x}^{(i)}=\mathbb{E}( x_{a}^{(i)})$
       \State \textbf{Convergence criterion} \Comment{e.g. mean based}
       \State $\quad \textrm{err}=\|\mathbb{E}(x_{a}^{(i)})-\mathbb{E}(x_{a}^{(i-1)})\|\cdot \|\mathbb{E}(x_{a}^{(i-1)})\|^{-1}$
        
      \EndWhile
%
       \EndFunction

\end{algorithmic}
\end{algorithm*}

\begin{figure}
\begin{center}
\includegraphics{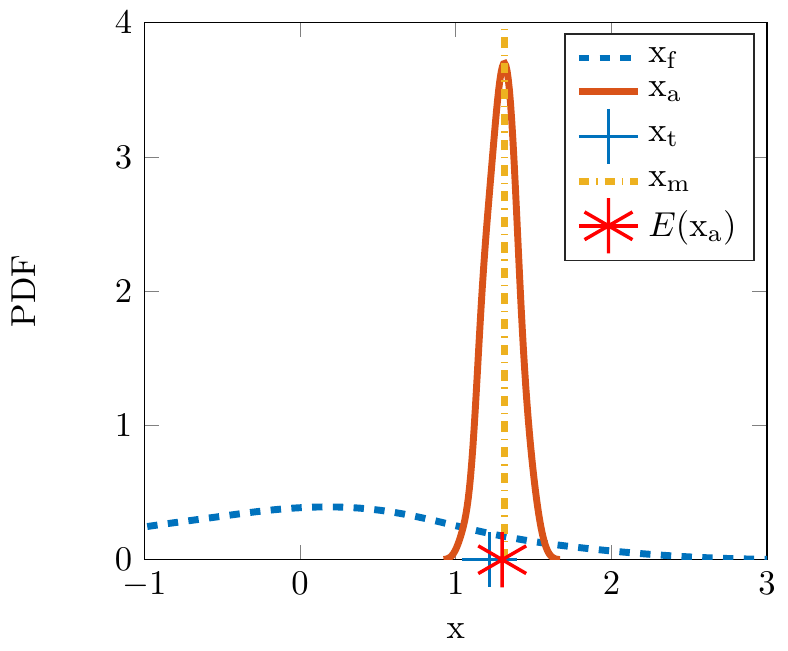}
\caption{Update of the current state $x$ at $t=6$ days}
\label{fig_first_state_filter}
\end{center}
\end{figure}

The initial value estimation via a pseudo-time stepping Gauss-Newton procedure can be done in different ways. Here, two 
variants are considered: the mean-based and the random variable-based smoothing. Both start with 
filtering of the current state $x_k$ given the measurement data $y_k^{mes}$ at $t_k$ via
\begin{equation}
\label{first_stage}
  x_{k,a}^{(i+1)}={x_{k,f}}+K_k^{(i)}(y_k^{mes}-y_{kh}^{(i)}(x_{k,f}))
 \end{equation}
in which ${x_{k,f}}$ is the prior knowledge on the current state, and $y_{kh}^{(i)}(x_{k,f})$ is the 
measurement prediction. As
$y_{kh}^{(i)}(x_{k,f})$ is linear in the state ${x_{k,f}}$, the iterative filter in 
\refeq{first_stage} consists of only one iteration. Fig.~(\ref{fig_first_state_filter})
shows the posterior probability density function of $x_a$ of the current state $x$ after six days of integration
given the perturbed full measurement data $x_m=x_t+\hat{\varepsilon}$
and the measurement noise with $C_\varepsilon=(0.1x_t)^2I$. 

\begin{figure}
\begin{center}
\includegraphics[width=0.35\textwidth]{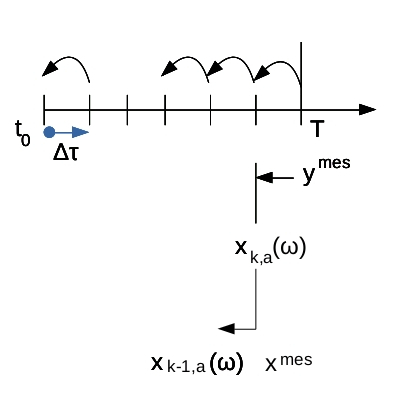}
\caption{The schematic representation of pseudo-backward propagation}
\label{fig_ps_pic}
\end{center}
\end{figure}

Once converged, the a posteriori state $ x_{k,a}$ is adopted as a pseudo-measurement 
for the preceding state $x_{k-\Delta \tau}$ at the time $t_k-\Delta \tau$. However, 
this could be done in at least two different ways: i) by assuming 
that the posterior mean is a pseudo-measurement and the posterior covariance 
is the measurement/modelling error describing our confidence 
in the ``measured'' value, or ii) by assuming that $ x_{k,a}$ is an uncertain ``perfect'' measurement, see \refig{fig_ps_pic}.

\subsection{Gaussian based pseudo-measurement}

Instead of evaluating the initial condition in \refeq{direct_estimate} 
directly one may use the ``smoothing'' procedure in which the intermediate states are estimated before the desired one, see \refig{fig_ps_pic}. 
In other words,
the first unknown state $x_k$ at the measurement time $t_k$ is estimated via \refeq{first_stage},
 whereas the preceding state $x_{k-\Delta \tau}$ at the time $t_k-\Delta \tau$ is further evaluated 
 given the Gaussian approximation 
 ${x}_{k,a}^g\sim\mathcal{N}(\bar{x}_{k,a},C_{{x}_{k,a}})$ of the convergent $x_{k,a}$ such that
 \begin{equation}
 \label{smoothing_filter12}
   x_{k-\Delta \tau,a}^{(i+1)}={x_{k-\Delta \tau,f}}+K_{k-\Delta \tau}^{(i)}({x}_{k,a}^g-y_{kh}^{(i)}({x_{k-\Delta \tau,f}}))
 \end{equation}
 holds. Here, $x_{k-\Delta \tau,f}$ is the apriori assumption on the state at the time $t_{k-\Delta \tau}$,
 $y_{kh}^{(i)}({x_{k-\Delta \tau,f}})$ is the linearised measurement operator 
 (i.e.~the linearised forward map $x_{k-\Delta \tau,f}\mapsto x_{k,f}$) around the point $\mathring{x}^{(i)}$ in iteration $(i)$:
 \begin{equation}
 \label{mes1}
  y_{kh}^{(i)}({x_{k-\Delta \tau,f}})=\mathring{H}^{(i)}(x_{k-\Delta \tau,f}-\mathring{x}^{(i)})+\mathring{h}^{(i)}.
 \end{equation}
 The map coefficients $\mathring{H}^{(i)},\mathring{h}^{(i)}$ are estimated either by the projection algorithm or by Bayesian 
 update similar to those depicted in \refalg{alg1}-\refalg{alg11}, whereas the linearisation point is chosen as  
 \begin{equation}
 \label{point_x_k_a}
  \mathring{x}^{(i+1)}=\mathbb{E}(x_{k-\Delta \tau,a}^{(i+1)}).
 \end{equation}
 Decoupling $x_{k,a}^g$ into the mean $\bar{x}_{k,a}$ and perturbation
 $\varepsilon_{k,f}\sim \mathcal{N}(0,C_{{x}_{k,a}})$ parts, one may rewrite \refeq{smoothing_filter12} to 
 \begin{equation}
 \label{smoothing_filter1}
   x_{k-\Delta \tau,a}^{(i+1)}={x_{k-\Delta \tau,f}}+K_{k-\Delta \tau}^{(i)}(\bar{x}_{k,a}-(y_{kh}^{(i)}+\varepsilon_{k,f})),
 \end{equation}
 thanks to the symmetry of the Gaussian distribution representing $\varepsilon_{k,f}$. In this manner 
 \refeq{smoothing_filter1} 
 can be understood as the state estimation given deterministic measurement $\bar{x}_{k,a}$ at the time $t_k$. Hence, 
 the algorithm of pseudo-time stepping is only a slight extension of the one presented in \refalg{alg1}. The new procedure requires 
 estimation of the current state, after which the original filter is called, see \refalg{alg2}.

\begin{algorithm*}
   \caption{\textbf{Pseudo-smoothing I (PS)}: incremental BGNMK (iBGNMK) with Gaussian approximation}
   \label{alg2}
    \begin{algorithmic}[1]
      \Function{iGNMK}{$x_{0,f},y^{mes},@integ, \Delta t, \Delta \tau,t_0,T,\varepsilon$}\Comment{Where $x_{0,f}$ - 
      prior on initial value, $t_0$ beginning of time interval, $T$ end of updating interval,
      $y^{mes}$ - measurement at $T$, $integ$ - forward function handle, $\Delta t$ - time integration step, $\Delta \tau$ - 
      update step, $\varepsilon$ - measurement error}
      \State \textbf{Predict current state at $T$}
      \State  $\quad x_{f}=integ(x_{0,f},\Delta t,t_0,T)$\Comment{integrate ODE system from $t_{0}$ to $T$ by time step $\Delta t$}
      \State \textbf{Predict measurement}
      \State $\quad y_f=I_x(x_f)+\varepsilon$ \Comment{$I_x$ is the indicator operator in case
      $\textrm{dim}(x_f)>\textrm{dim}(y^{mes})$}
      \State \textbf{Update current state at $T$ given $y^{mes}$}
       \State $\quad x_a=x_f+C_{x_f,y_f}C_{y_f,y_f}^{-1}(y^{mes}-y_f)$
        \For{$tt=T-\Delta \tau:-\Delta \tau:t_0$}
          \State \textbf{Set pseudo-measurement}
        \State $\quad x_a^g=\textrm{Gaussian}(x_a)$
         \State \textbf{Decompose pseudo-measurement}:
         \State $\quad$ to the mean value $\bar{x}_{a}=\mathbb{E}(x_{a}^g)$
         \State $\quad$ and the fluctuation term $\varepsilon_f:=x_{a}^g-\bar{x}_{a}$
         \State \textbf{Set preceding prior} 
          \State $\quad x_{f}=integ(x_{0,f},\Delta t,t_0,tt)$ \Comment{integrate ODE system from $t_{0}$ to $tt$ by time step $\Delta t$}
         \State \textbf{Update preceding state}
         \State $\quad x_a=GNMK(x_{f},\bar{x}_{a},@integ, \Delta t, \Delta \tau,tt,tt+\Delta \tau, \mathring{x},\varepsilon_f)$
      \EndFor
%
       \EndFunction

\end{algorithmic}
\end{algorithm*}
%

Rewriting Eqs.~(\ref{first_stage})-(\ref{smoothing_filter}) for all preceding states, one obtains the general form of a smoothing iterative filter:
\begin{equation}
 \label{smoothing_filter}
   x_{\ell-1,a}^{(i+1)}={x_{\ell-1,f}}+K_\ell^{(i)}(\bar{x}_{\ell,a}-(y_{\ell,h}^{(i)}({x_{\ell,f}})+\varepsilon_{\ell})),
 \end{equation}
for all $\ell=k,k-\Delta \tau,..,k-n\Delta \tau$. The last formula further can be 
 generalised by taking into account all estimated states from the time moment $t_{k}$ to the current time $t_n$
 as measurements, similarly to
 the classical smoothing algorithm. 
 
\begin{figure*}
\begin{center}
\includegraphics{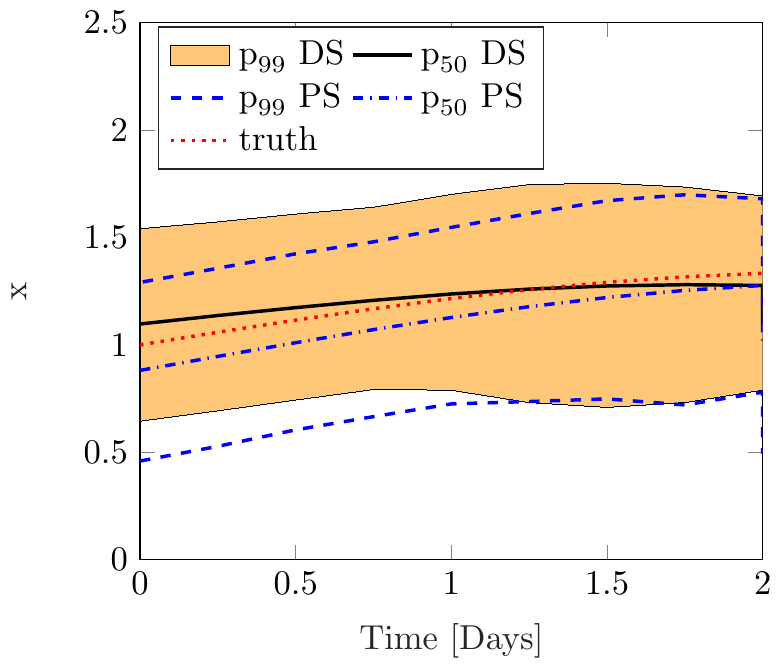}
\caption{Bias propagation over time for the first Lorenz 1984 component. DS is the direct simulation estimate given in \refeq{direct_estimate} 
and PS is the pseudo-estimate given in \refeq{smoothing_filter1} with $\Delta \tau=6h$. $p_{n}$ denotes 
$n\%$ quantile.}
\label{fig_second}
\end{center}
\end{figure*}

 Unfortunately, the estimate in \refeq{smoothing_filter} is biased 
 due to nonlinearity of the time-dependent problem.  If not corrected, the bias becomes propagated through the 
 model with each new update as shown in Fig.~(\ref{fig_second}) on the example of the first Lorenz 1984 component. 
 The 
 mean value deteriorates from the measured one with each update such that the deviation becomes larger with the reduction 
 of the update step size $\Delta \tau$ 
  in contrast to expectations.


%

 \begin{figure*}[htpb]
 \begin{center}
 \includegraphics{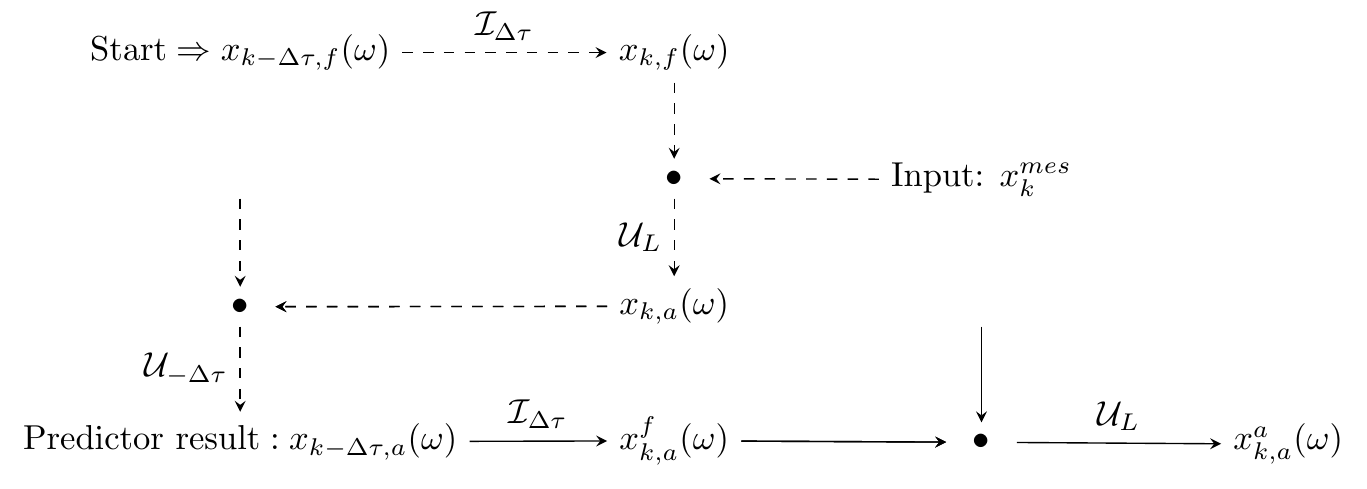}
\end{center}
\caption{The scheme of bias correction}
\label{bias_scheme}
\end{figure*}

\begin{figure*}
\begin{center}
\includegraphics{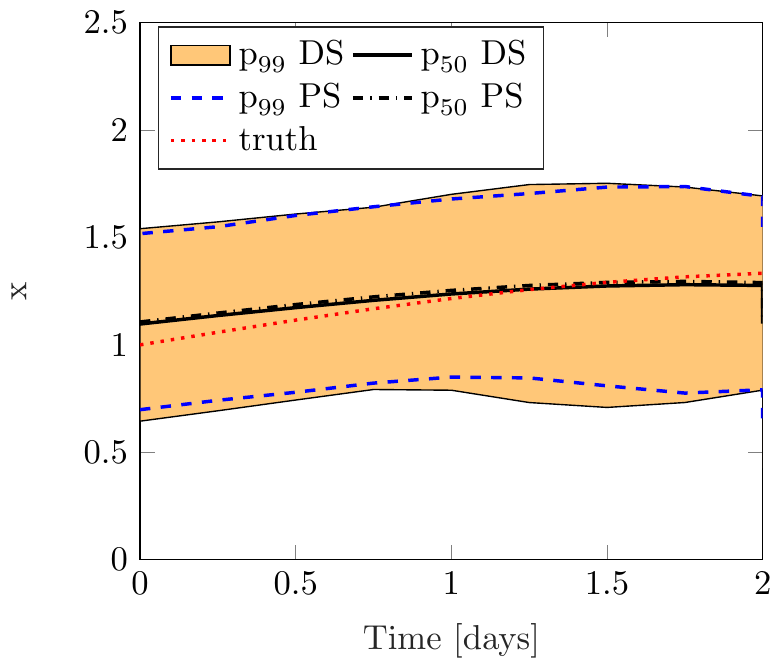}
\includegraphics{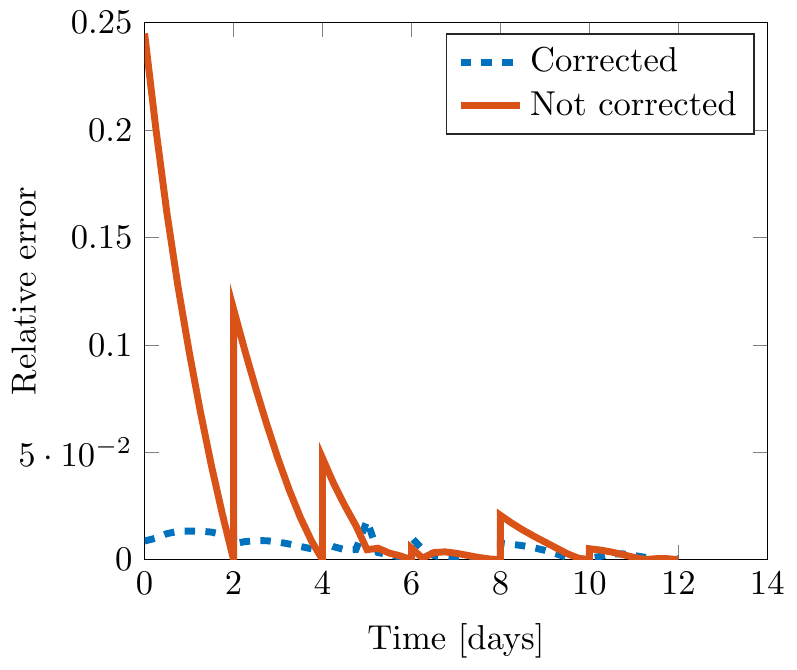}
\caption{The bias propagation in time (left) and the bias correction in time (right)}
\label{fig_unbiased}
\end{center}
\end{figure*}

\begin{figure*}
\begin{center}
\includegraphics{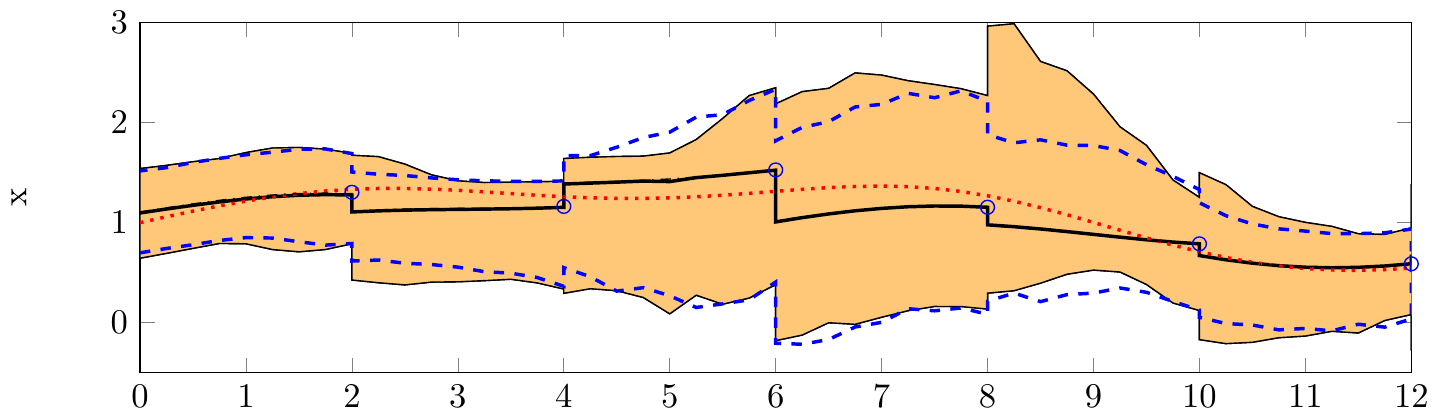}\\
\includegraphics{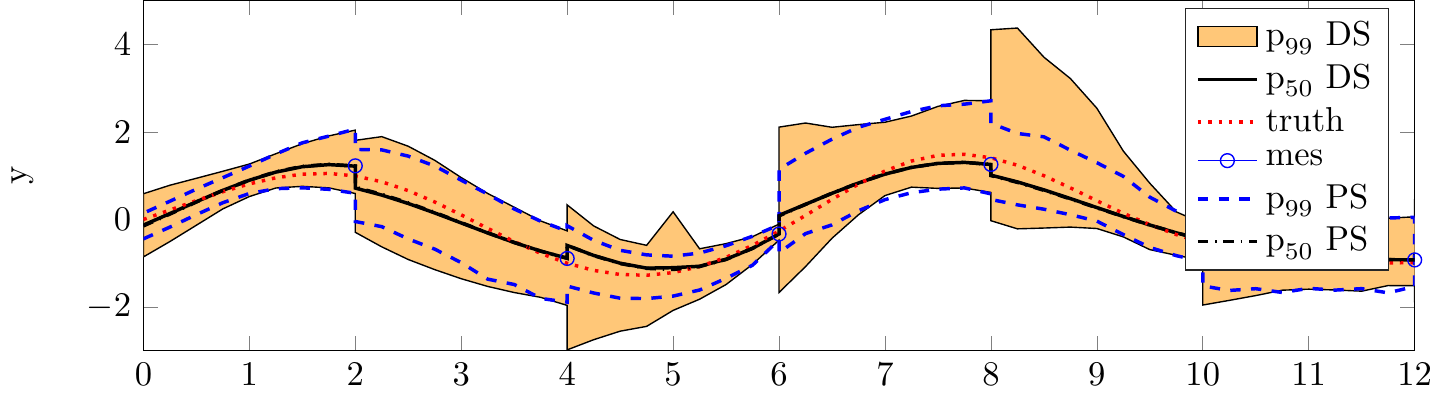}\\
\includegraphics{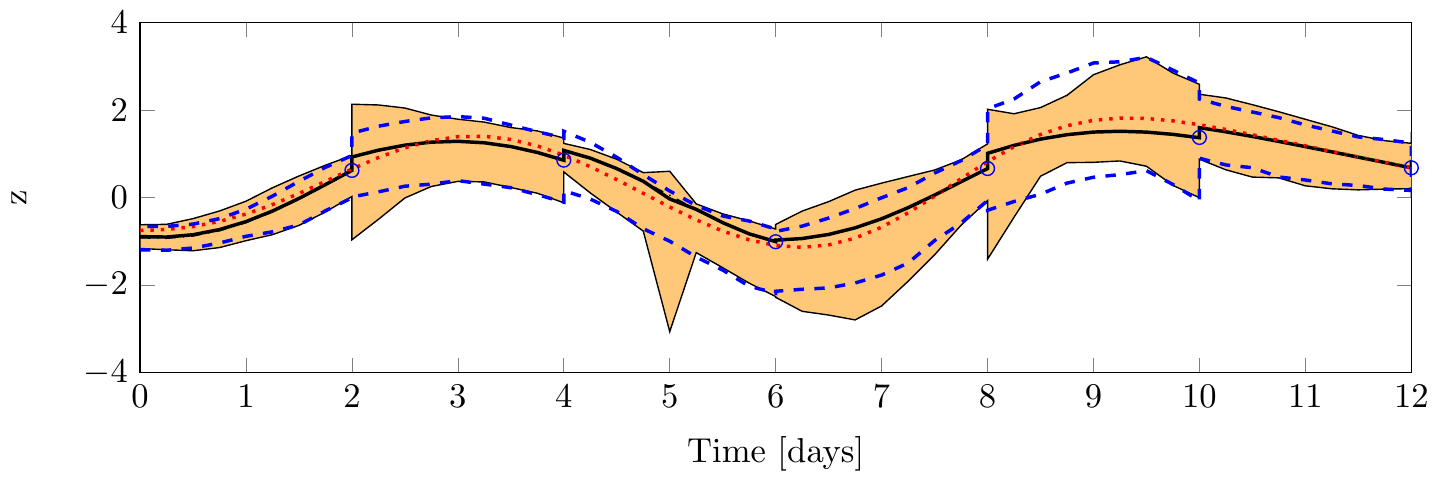}
\caption{The mean based corrected estimate of the Lorenz 1984 state every two days backwards. The prior is obtained starting from the initial condition 
and the measurement has the coefficient of the variance equal to $10\%$.}
\label{fig_corrected_det}
\end{center}
\end{figure*}


 The posterior $x_{k,a}$ in \refeq{smoothing_filter} has the mean $\bar{x}_{k,a}$
 that differs from the true posterior mean $\bar{x}_{k,a}^{true}$ according to the error 
 \begin{equation}
 \label{eq:correction_term}
  \epsilon_k=\bar{x}_{k,a}^{true}-\bar{x}_{k,a},
 \end{equation}
 which further becomes propagated in time with the state integration/assimilation. 
 Hence, \refeq{smoothing_filter1} (and similarly \refeq{smoothing_filter}) have to be corrected for the amount given in \refeq{eq:correction_term}. 
 
 The correction scheme is schematically depicted in \refig{bias_scheme} and is 
  of the predictor-corrector type. The predictor phase starts with
 \begin{itemize}
  \item the prior assumption on the state $x_{k-\Delta \tau,f}$ at the time $t_{k-\Delta \tau}$ with $\Delta \tau$
  being the backpropagation increment.
  \item The state $x_{k-\Delta \tau,f}$ is integrated forward ($\mathcal{I}_{\Delta \tau}$ in \refig{bias_scheme} denotes 
  the integration operator over time interval $\Delta \tau$ from $t_{k-\Delta \tau}$ to $t_k$) to obtain the 
  current prior state $x_{k,f}$ at the time $t_k$.
  \item The current state  $x_{k,f}$ is further assimilated with the measurement data $x_k^{mes}$ 
  in a linear direct GMK manner (in \refig{bias_scheme} denoted by $\mathcal{U}_L$)  to obtain the 
  posterior $x_{k,a}$. $x_k^{mes}$ may 
  represent the real data only for the state that is being measured, otherwise these
  are pseudo-measurement data. For example, if we update in the time interval $[t_0,T]$ given measurement 
  data at the time $T$, then the measurement $x_k^{mes}$ at $t_k=T$ is the real measurement $y^{mes}$. 
  Otherwise, if $t_k<T$ our measurement at $t_k$ is the posterior estimate obtained by incremental backpropagation 
  of the posterior at $t_k+\Delta \tau$. 
  \item The assimilated current state $x_{k,a}$ is then used as
 a pseudo-measurement for the assimilation of $x_{k-\Delta t,f}$ state via iterative GMK (see \refeq{smoothing_filter1}), in \refig{bias_scheme}
 denoted by 
 backward update operator $\mathcal{U}_{-\Delta \tau}$.
 \end{itemize}
With this the corrector phase starts by 
\begin{itemize}
 \item integrating forward the estimate $x_{k-\Delta \tau,a}$ via $\mathcal{I}_{\Delta \tau}$
 to obtain the prior on the current state
 $x_{k,f}^{a}$ at $t_k$ given posterior $x_{k-\Delta \tau,a}$ at $t_{k-\Delta \tau}$.
 \item Furthermore, 
 the newly obtained estimate $x_{k,a}^f$ is used as a prior
 for a second turn of updating the current state at $t_k$ given measurement $x_k^{mes}$. The update is performed using linear direct GMK rule 
to obtain $x_{k,a}^a$.
\item The difference between the prior $x_{k,a}^f$ and posterior $x_{k,a}^a$ 
estimates then defines the correction error. This is the corrector phase.
 The process is further repeated for $x_{k-2\Delta \tau,f}$ given the
 measurement $x_{k-\Delta \tau}^{mes}$ adopted as the corrected version of $x_{k-\Delta \tau,a}$.
\end{itemize}

 To estimate the correction error, let the converged posterior estimate in \refeq{smoothing_filter1} 
be denoted by $x_{k-\Delta \tau,a}$ (beginning of the corrector phase in \refig{bias_scheme}) such that 
 \begin{equation}
  x_{k-\Delta \tau,a}=x_{k-\Delta \tau,a}^{true}+e_{k-\Delta \tau}
 \end{equation}
holds, in which $e_{k-\Delta \tau}$ denotes the bias error at the time $t_{k-\Delta \tau}$. Propagating 
the a posteriori estimate $x_{k-\Delta \tau,a}$ by time step $\Delta \tau$ forward \footnote{this may include several
time discretisation steps $\Delta t$}, one obtains the 
forecast estimate $x_{k,a}^f$ at $t_k$ such that 
\begin{eqnarray}
  x_{k,a}^f&=&\mathring{H}_k({x}_{k-\Delta \tau,a}-\mathring{x}_k)+\mathring{h}_k\\
  &=& \mathring{H}_k(x_{k-\Delta \tau,a}^{true}+e_{k-\Delta \tau}-\mathring{x}_k)+\mathring{h}_k\nonumber\\
  &=& \mathring{H}_k(x_{k-\Delta \tau,a}^{true}-\mathring{x}_k)+\mathring{h}_k+\mathring{H}_ke_{k-\Delta \tau}\nonumber\\
  &=& x_{k,a}^{f,true}+\mathring{H}_k e_{k-\Delta \tau}
 \end{eqnarray}
 holds. Here, $\mathring{H}_k$ and $\mathring{h}_k$ are converged parameters of the forward map, and 
 $x_{k,a}^{f,true}$ denotes the forecast of the exact a posteriori estimate. 
 The analysis step at time moment $t_k$ is then given by
 \begin{eqnarray}\label{eqq:two}
  x_{k,a}^a&=&x_{k,a}^f+K({x}_{k}^{mes}-x_{k,a}^f-\varepsilon_{k,f})\\
  &=& x_{k,a}^{f,true}+\mathring{H}_k e_{k-\Delta \tau}+\nonumber\\
  && K({x}_{k}^{mes}-x_{k,a}^{f,true}-\mathring{H}_k e_{k-\Delta \tau}-\varepsilon_{k,f})\nonumber\\
  &=&  x_{k,a}^{f,true}+K({x}_{k}^{mes}-x_{k,a}^{f,true})\nonumber\\
  && +\mathring{H}_k (I-K)e_{k-\Delta \tau}-K\varepsilon_{k,f}\nonumber\\
  &=& x_{k,a}^{a,true} +\mathring{H}_k (I-K)e_{k-\Delta \tau}-K\varepsilon_{k,f}. \nonumber
 \end{eqnarray}
 Here, $ x_{k,a}^{a,true}$ is the assimilated value of $x_{k,a}^{f,true}$.
 By subtracting the previous two equations 
 \begin{eqnarray}\label{eqq:one}
   x_{k,a}^f-x_{k,a}^a&=&x_{k,a}^{f,true}-x_{k,a}^{a,true}+\\
  && \mathring{H}_k e_{k-\Delta \tau}-\mathring{H}_k (I-K)e_{k-\Delta \tau}\nonumber\\
  &&+K\varepsilon_{k,f}\nonumber
 \end{eqnarray}
 and taking the mathematical expectation one obtains
 \begin{eqnarray}\label{eq:unerr}
   \mathbb{E}(x_{k,a}^f-x_{k,a}^a)&=&\mathbb{E}(x_{k,a}^{f,true}-x_{k,a}^{a,true})+\\
   && \mathring{H}_k \bar{e}_{k-\Delta \tau}-\mathring{H}_k (I-K)\bar{e}_{k-\Delta \tau}.\nonumber
 \end{eqnarray}
 Furthermore, 
 \begin{eqnarray}\label{eq:unerr1}
   \mathbb{E}(x_{k,a}^f-x_{k,a}^a)&=&\mathbb{E}(x_{k,a}^{f,true}-x_{k,a}^{true})\nonumber\\
   &&+\mathbb{E}(x_{k,a}^{true}-x_{k,a}^{a,true})\\
   && +\mathring{H}_k K\bar{e}_{k-\Delta \tau}\nonumber
 \end{eqnarray}
 in which
 \begin{eqnarray}\label{eqqunbreq}
  \mathbb{E}(x_{k,a}^{f,true}-x_{k,a}^{true})&\mbeq&0\\
   \mathbb{E}(x_{k,a}^{a,true}-x_{k,a}^{true})&\mbeq&0 \nonumber
 \end{eqnarray}
 due to unbiased requirement. This further 
gives
\begin{eqnarray}
   \mathbb{E}(x_{k,a}^f-x_{k,a}^a)&=&\mathring{H}_k K\bar{e}_{k-\Delta \tau}.\nonumber
 \end{eqnarray}
Hence, the mean bias error for the assimilated state $x_{k-\Delta \tau,a}$ at the end of the predictor phase 
reads:
 \begin{eqnarray}
 \label{estimated_error}
   \bar{e}_{k-\Delta \tau}&=&(\mathring{H}_k K)^{-1} \mathbb{E}(x_{k,a}^f-x_{k,a}^a).
 \end{eqnarray}
 In a similar manner one may correct the variance 
 of the posterior by considering the second moment in \refeq{eq:unerr}.

 Introducing the estimated error in \refeq{estimated_error} to the update in \refeq{smoothing_filter1} one obtains the unbiased solution as shown in
 Fig.~(\ref{fig_unbiased}a) for the update of the first Lorenz 1984 component. The total correction 
 over a period of 12 days is shown in Fig.~(\ref{fig_unbiased}b), in which are depicted the relative errors of the biased 
 and unbiased pseudo-estimated states 
 compared to the direct estimated state following \refeq{direct_estimate}. As one may notice the error 
 is decreasing for several orders of magnitudes when the correction is introduced.

 \subsection{Random-variable based pseudo-measurement}
 
  The previous estimation did not take into consideration the full uncertainty in the pseudo-measurement. 
 Hence, the estimate does not have correct variance as only the Gaussian approximation of the measurement is considered.
   This can be seen in \refig{fig_corrected_det} in which the corrected pseudo-estimate is compared to the direct one.
   
 However, by taking the current aposteriori estimate $x_{k,a}$ at the time $t_k$--- 
 obtained by assimilating the measurement data $y^{mes}$ at $t_k$ via linear GMK filter--- 
 as uncertain non-Gaussian pseudo-measurement,
 the estimation of the preceding state $x_{k-\Delta \tau}$ in a backpropagation manner ($t_k\rightarrow t_{k-\Delta \tau}$)
 becomes 
 stochastic as the measurement is a random variable. Following 
 this, one may further state
 \begin{equation}
 \label{smoothing_filter1_gen}
   x_{k-\Delta \tau,a}^{(i+1)}={x_{k-\Delta \tau,f}}+K_{k-\Delta \tau}^{(i)}({x}_{k,a}-
   y_{kh}^{(i)}({x_{k-\Delta \tau,f}})),
 \end{equation}
 similarly to \refeq{smoothing_filter1}. However, in contrast to \refeq{smoothing_filter1} 
 the pseudo-measurement ${x}_{k,a}$ is taken in its full form, and not only as a Gaussian approximation. 
 This further means that 
 $y_{kh}^{(i)}({x_{k-\Delta \tau,f}})$ 
 is a ``perfect'' linearised 
 version of the time-discretised model in \refeq{initvp} around point $\mathring{x}_k^{(i)}$
 \begin{equation}
  y_{kh}^{(i)}({x_{k-\Delta \tau,f}})=\mathring{H}_k^{(i)}(x_{k-\Delta \tau,f}-
  \mathring{x}_k^{(i)})+\mathring{h}_k^{(i)}+\epsilon_k,
 \end{equation}
 and similarly $K_{k-\Delta \tau}^{(i)}$ is the ``perfect'' Kalman gain given as
 \begin{equation}\label{kalman_rv}
  K_{k-\Delta \tau}^{(i)}=C_{x_{k-\Delta \tau,f},y_{kh}^{(i)}}
  C_{y_{kh}^{(i)}}^\dagger.
 \end{equation}
Notice that $\epsilon_k$ represents the modelling/ discretisation error, the estimate of which is  
further described in Section \ref{spaopt}. 

 The posterior estimate in \refeq{smoothing_filter1_gen} has different second order statistics than those specified by
 the ``classical'' Kalman filter in the previous section.
To simplify the notation let $x_f:={x_{k-\Delta \tau,f}}$, $x_a:={x_{k-\Delta \tau,a}}$, $y_f:=y_{kh}^{(i)}({x_{k-\Delta \tau,f}})$ and $y_m:=x_{k,a}$, then 
 the mean value of the converged posterior reads
\begin{eqnarray}
   \bar{x}_{a}&=&{\bar{x}_{f}}+K_{k-\Delta \tau}(\bar{y}_m- \bar{y}_{f}),
 \end{eqnarray}
 whereas the covariance follows from 
 \begin{eqnarray}\label{eqdercov}
  C_{{x}_{a}}&=&{C_{{x}_{f}}}+K_{k-\Delta \tau}C_{y_m}K_{k-\Delta \tau}^T\nonumber\\
  &&+K_{k-\Delta \tau}C_{y_f}K_{k-\Delta \tau}^T\nonumber\\
  &&-2K_{k-\Delta \tau}C_{y_f,y_m}K_{k-\Delta \tau}^T\nonumber\\
 && -C_{x_f,y_f}K_{k-\Delta \tau}^T\nonumber\\
&&-K_{k-\Delta \tau}  C_{x_f,y_f}^T\nonumber\\
&&+K_{k-\Delta \tau}C_{x_f,y_m}\nonumber\\
&&+C_{x_f,y_m}K_{k-\Delta \tau}^T.
 \end{eqnarray}
 In the previous equations the index $^{(i)}$ is avoided, as the last 
 two equations are written for
 $i=i_{conv}$ in which $i_{conv}$ is the number of iterations of the converged estimate. 
%
 In \refeq{eqdercov} note that
 \begin{eqnarray}
   &&-2K_{k-\Delta \tau}C_{y_f,y_m}K_{k-\Delta \tau}^T\nonumber\\
   &&=-2K_{k-\Delta \tau}\mathring{H}_k^{(i)}C_{x_f,y_m}K_{k-\Delta \tau}^T\nonumber\\
   &&=-2C_{x_f,y_m}K_{k-\Delta \tau}^T
 \end{eqnarray}
as the Kalman gain is optimal \footnote{In numerical computations $K_{k-\Delta \tau}\mathring{H}_k^{(i)}\approx I$}, i.e.~$K_{k-\Delta \tau}\mathring{H}_k^{(i)}=I$.
Having that $K_{k-\Delta \tau}=C_{x_f,y_f}C_{y_f}^\dagger$ and after 
substituting the last equation in \refeq{eqdercov} one obtains 
\begin{eqnarray}
  C_{{x}_{k-\Delta \tau,a}}&=&C_{x_f}+C_{x_{f},y_f}C_{y_f}^\dagger
  C_{y_m}(C_{y_f}^\dagger)^TC_{x_f,y_f}^T\nonumber\\
  &&+C_{x_{f},y_f}C_{y_f}^\dagger C_{y_f}(C_{y_f}^\dagger)^TC_{x_{f},y_f}^T\nonumber\\
 && -C_{x_f,y_f}(C_{y_f}^\dagger)^TC_{x_f,y_f}^T\nonumber\\
&&-C_{x_{f},y_f}C_{y_f}^\dagger C_{x_f,y_f}^T\\
&=& C_{x_f}+C_{x_{f},y_f}C_{y_f}^\dagger\nonumber\\
&&(C_{y_m}-C_{y_f})(C_{y_f}^\dagger)^TC_{x_f,y_f}^T\nonumber
 \end{eqnarray}
%
which in the original notation reads
 \begin{eqnarray}
 \label{rv_cov_final}
  C_{{x}_{k-\Delta \tau,a}}&=&{C_{{x}_{k-\Delta \tau,f}}}+
  {C_{{x}_{k-\Delta \tau,f}}}\mathring{H}^TC_{y_{kh}^{(i)}}^\dagger \nonumber \\
  && (C_{x_{k,a}}-C_{y_{kh}^{(i)}})
  {C_{y_{kh}^{(i)}}^\dagger}\mathring{H}^T{C_{{x}_{k-\Delta \tau,f}}}^T.
 \end{eqnarray}

\begin{algorithm*}
   \caption{\textbf{Pseudo-smoothing II (PS)}: incremental RV-based GNMK (irvGNMK)}
    \begin{algorithmic}[1]
      \Function{irvGNMK}{$x_{0,f},y^{mes},@integ, \Delta t, \Delta \tau,t_0,T, \mathring{x},\varepsilon$}\Comment{Where $x_{0,f}$ - 
      prior on initial value, $t_0$ beginning of time interval, $T$ end of updating interval,
      $y^{mes}$ - measurement at $T$, $integ$ - forward function handle, $\Delta t$ - time integration step, $\Delta \tau$ - 
      update step, $\varepsilon$ - measurement error}
      \State \textbf{Predict current state at $T$}
      \State  $\quad x_{f}=integ(x_{0,f},\Delta t,t_0,T)$\Comment{integrate ODE system from $t_{0}$ to $T$ by time step $\Delta t$}
      \State \textbf{Predict measurement}
      \State $\quad y_f=I_x(x_f)+\varepsilon$ \Comment{$I_x$ is the indicator operator in case
      $\textrm{dim}(x_f)>\textrm{dim}(y^{mes})$}
      \State \textbf{Update current state at $T$ given $y^{mes}$}
       \State $\quad x_a=x_f+C_{x_f,y_f}C_{y_f,y_f}^{-1}(y^{mes}-y_f)$
        \For{$tt=T-\Delta \tau:-\Delta \tau:t_0$}
          \State \textbf{Set pseudo-measurement}
        \State $\quad x_a^m=x_a$
        \State \textbf{Set preceding prior}
        \State  $\quad x_{f}=integ(x_{0,f},\Delta t,t_0,T)$\Comment{integrate ODE system from $t_{0}$ to $tt$ by time step $\Delta t$}
         \State \textbf{Update preceding state}
         \State $\quad x_a=orvGNMK(x_{f},x_a^{m},@integ, \Delta t, \Delta \tau,tt,T,\varepsilon)$
      \EndFor
%
       \EndFunction

\end{algorithmic}
\end{algorithm*}

\begin{algorithm*}
   \caption{\textbf{Pseudo-smoothing (PS) II}: incremental BGNMK filter, backpropagation}
   \label{alg31}
    \begin{algorithmic}[1]
      \Function{orvGNMK}{$x_{f},x_a^{m},@integ, \Delta t, \Delta \tau,tt,T,\varepsilon$}\Comment{Where $x_{0,f}$ - 
      prior on initial value, $t_0$ beginning of time interval, $T$ end of time interval,
      $y^{mes}$ - measurements, $integ$ - forward function (model) handle, $\Delta t$ - time integration step, $\Delta \tau$ - 
      update step, $\varepsilon$ - measurement error}
      \State
          \State  \textbf{Set linearisation point}
        \State $\quad \mathring{x}^{(0)}=\mathbb{E}(x_f), \quad x_a^{(0)}=x_f$
          \State \textbf{Set} $i=0$, $\textrm{err}=2\cdot\textrm{tol}$, $\textrm{maxiter}=100$
        \While{$i\leq \textrm{maxiter}$ \& $\textrm{err}\leq \textrm{tol}$} 
          \State  \textbf{Predict measurement}
         \State $\quad z_{a}^{(0)}=integ(x_{a}^{(i)},\Delta t,tt,T)$\Comment{integrate ODE system from $tt$ to $T$}
          \State \textbf{Approximate forward map $x_a^{(i)} \mapsto z_a^{(i)}:=Y(x_a^{(i)})$ by}
          \State $\quad$ $\varphi_y(x_a^{(i)})=\mathring{H}^{(i)}(x_a^{(i)}-\mathring{x}^{(i)})+\mathring{h}^{(i)}$
            \State \textbf{Estimate forward map coefficients $\beta:=(\mathring{H}^{(i)},\mathring{h}^{(i)})$ by}
        \State $\quad$ - projection:
        \State $\quad\quad\quad$   $\mathring{H}^{(i)}=C_{z_{a}^{(i)},x_{a}^{(i)}}C_{x_{a}^{(i)}}^\dagger,\quad 
        \mathring{h}^{(i)}=\mathbb{E}(z_{a}^{(i)})-\mathring{H}^{(i)}(x_{a}^{(i)}-
        \mathring{x}^{(i)}),$
        \State $\quad$ - or by Bayes's rule 
        \State $\quad\quad\quad$ given data $d^{sim}=(x_a(\omega_j)^{(i)},z_a(\omega_j)^{(i)}), j=1,...,N$ (see Section \ref{spaopt})
        \State $\quad\quad\quad$ update $\pi_{\beta|d^{sim}}(\beta|d^{sim})\propto\pi_{d^{sim}|\beta}
 (d^{sim}|\beta)
 \pi_{\beta}(\beta)$
        \State \textbf{Linearise predicted measurement}
        \State $\quad y_{\ell}^{(i)}({x_{f}})=\mathring{H}^{(i)}(x_{f}-\mathring{x}^{(i)})+\mathring{h}^{(i)}+\varepsilon$
        \State \textbf{Approximate inverse map $y_\ell^{(i)} \mapsto x_f$ by}
          \State $\quad$ $\varphi(y_\ell^{(i)})=K^{(i)}y_\ell^{(i)}+b^{(i)}$
       \State \textbf{Estimate inverse map coefficients $w:=(K^{(i)},b^{(i)})$ by}
        \State $\quad$ - projection:
        \State $\quad\quad\quad$   $K^{(i)}=C_{x_{f},y_{\ell}^{(i)}}C_{y_{\ell}^{(i)}}^\dagger, \quad b=\mathbb{E}(x_f)-K^{(i)}\mathbb{E}(y_\ell^{(i)})$ 
        \State $\quad$ - or by Bayes's rule 
        \State $\quad\quad\quad$ given data $d^{sim}=(x_f(\omega_j),y_{\ell}^{(i)}(\omega_{(j)}), j=1,...,N$ (see Section \ref{spaopt})
        \State $\quad\quad\quad$ update $\pi_{w|d^{sim}}(w|d^{sim})\propto\pi_{d^{sim}|w}
 (d^{sim}|w)
 \pi_{w}(w)$
       \State \textbf{Update state}
       \State $\quad$ $ x_{a}^{(i+1)}={x_{f}}+K^{(i)}(x_a^{m}-y_{\ell}^{(i)}),\quad $
       \State \textbf{Update linearisation point}
       \State $\quad i=i+1;$
       \State $\quad \mathring{x}^{(i)}=\mathbb{E}( x_{a}^{(i)})$
       \State \textbf{Convergence criterion} \Comment{e.g. mean based}
       \State $\quad \textrm{err}=\|\mathbb{E}(x_{a}^{(i)})-\mathbb{E}(x_{a}^{(i-1)})\|\cdot \|\mathbb{E}(x_{a}^{(i-1)})\|^{-1}$
        
      \EndWhile
%
       \EndFunction

\end{algorithmic}
\end{algorithm*}

\begin{figure*}
\begin{center}
\includegraphics{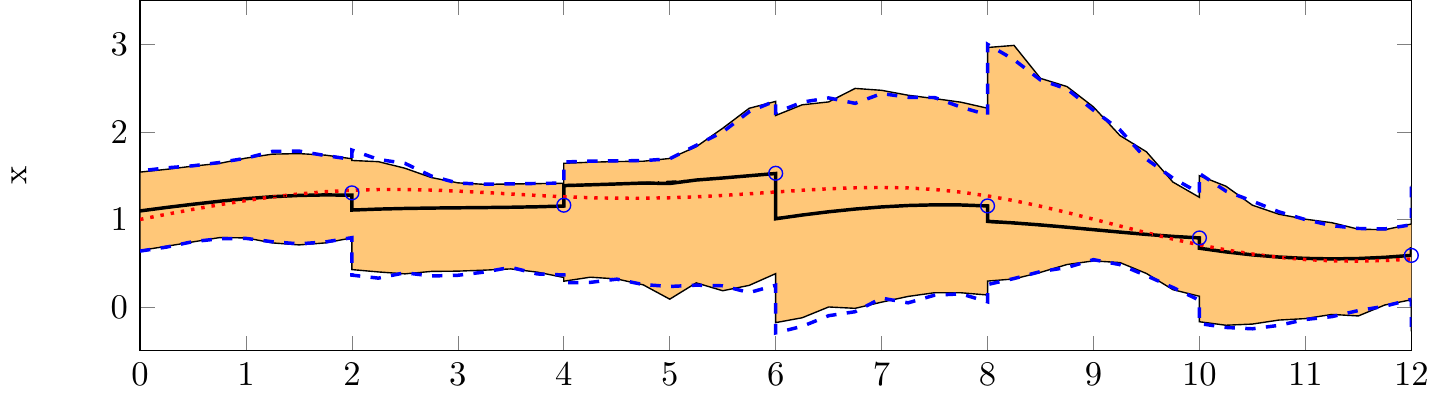}\\
\includegraphics{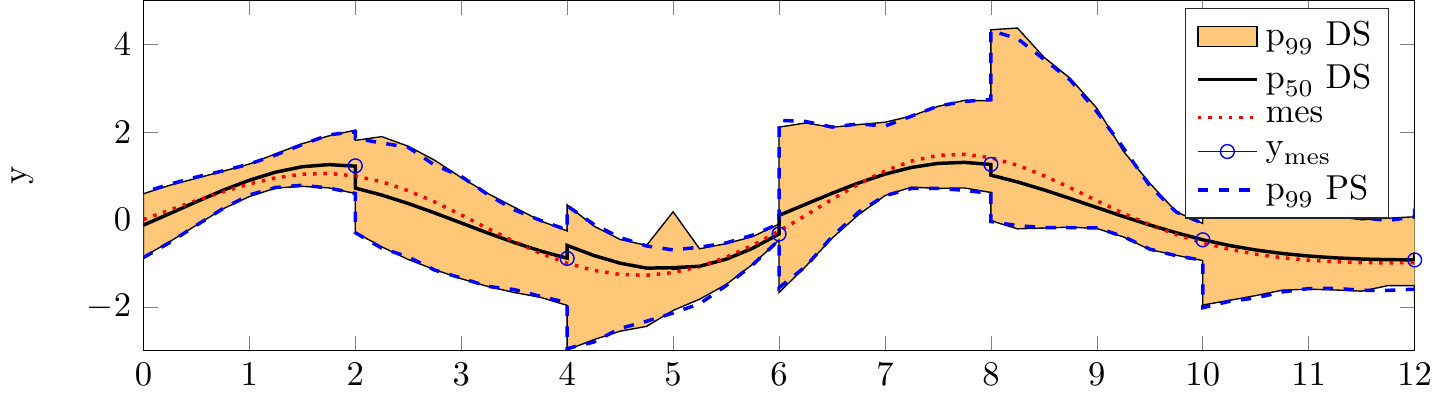}\\
\includegraphics{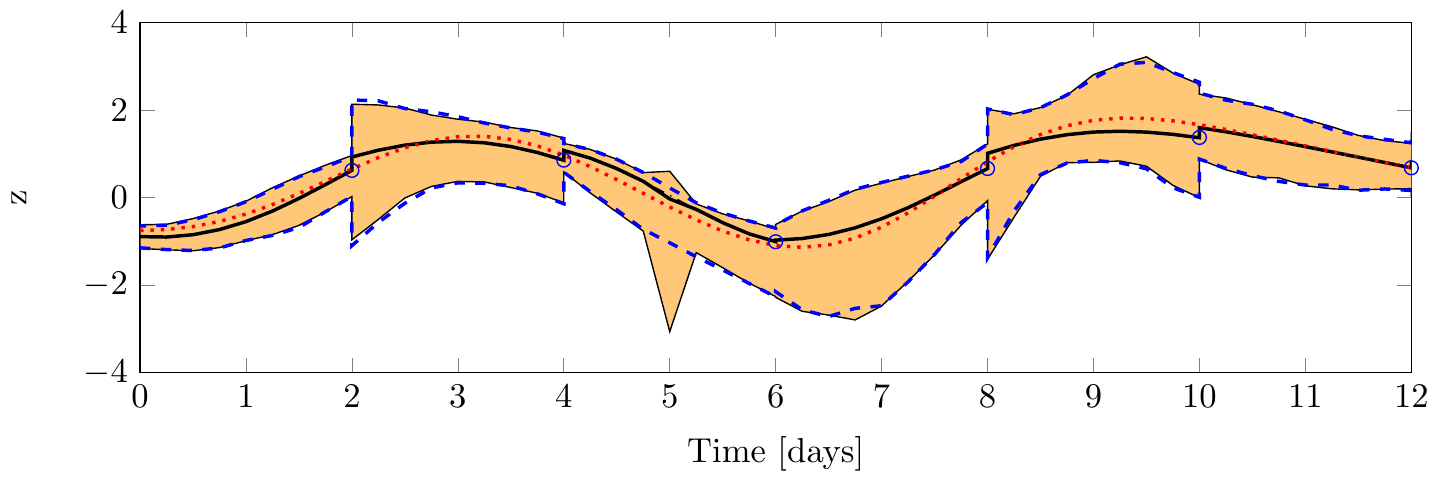}
\caption{The pseudo-estimation of the Lorenz 1984 state every two days backwards using the random variable algorithm. The pseudo-update is
made every 6 hours.}
\label{fig_filter_rv_ja}
\end{center}
\end{figure*}

\begin{figure}
\begin{center}
\includegraphics{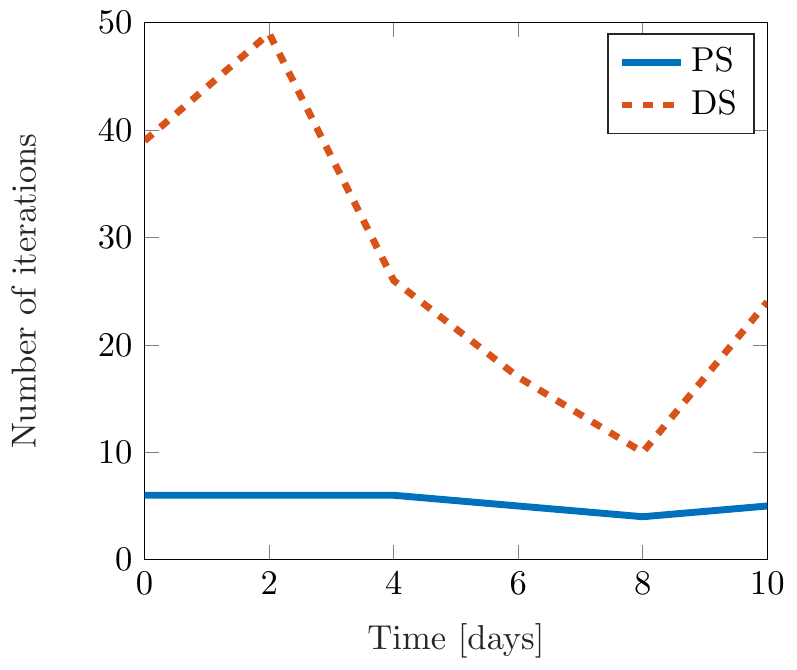}
\caption{Number of iterations neccessary to achieve posterior accuracy of 1e-3 in the mean}
\label{fig_filter_rv_num_iter}
\end{center}
\end{figure}

Using the estimation in \refeq{smoothing_filter1_gen} one obtains the correct estimate of the posterior variance as obtained by the direct simulation, see
\refig{fig_filter_rv_ja} for the comparison of the update obtained by direct iteration (DS) and the pseudo (PS) one. Note that the pseudo-updating is here performed 
every 6 hours. The same estimate is also depicted earlier in \refig{fig_linear_nonlinear_update}, in which the iterative pseudo-estimation is 
compared to the linear pseudo-estimation every 6 hours. The pseudo-nonlinear posterior estimate converges faster than the direct one, see \refig{fig_filter_rv_num_iter}
for comparison of the number of iterations neccessary to achieve the relative error in the posterior mean of magnitude 1e-3. Usually the posterior converges
very fastly 
after two or three iterations up to tolerance on the first decimal. However, this number raises up to ten iterations if the accuracy
is up to 1e-3 in all three components.
On the other hand, a direct iteration of the initial condition requires up to 50 iterations for the same accuracy. This behaviour 
also depends on the discretisation 
of the previously described filters which will be discussed later. 

\begin{figure*}
\begin{center}
\includegraphics{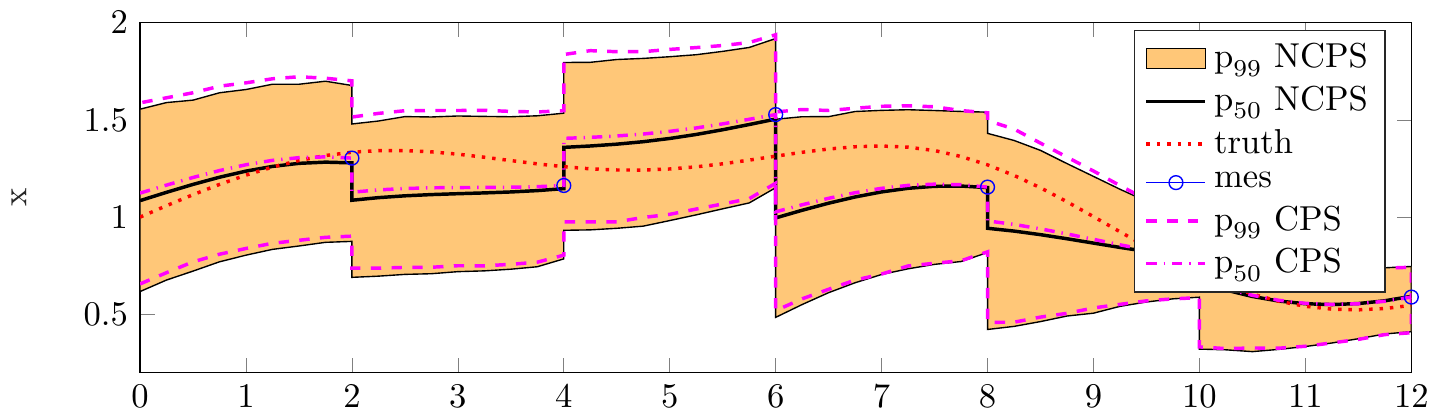}
\caption{The random variable pseudo-update of the Lorenz 1984 first component with (CPS) and without (NSP) correction.}
\label{fig_filter_rv_corr_non}
\end{center}
\end{figure*}

The random variable updating does not introduce bias into the estimation, and hence the bias correction introduced earlier
does not change much the posterior 
estimate, see \refig{fig_filter_rv_corr_non}. A small difference between the corrected and original estimates exists
due to numerical integration of discretisation errors.

\section{Iterative polynomial chaos filter}\label{itpce}

The advantage of the filtering approach as presented in \refeq{newton_gauss_markov1}, \refeq{smoothing_filter1} and \refeq{smoothing_filter1_gen} 
compared to the other Bayesian 
numerical procedures lies in the simplicity of the posterior variable estimation. Once the random variables appearing in 
\refeq{newton_gauss_markov1} are approximated using the standard Galerkin functional approximation tools in their minimal form, the filtering procedure 
reduces to the purely algebraic method for estimating the posterior variable. However, in high-dimensional problems, or when using commercial softwares, 
sometimes it is not possible to use spectral, but pseudo-spectral approximations. Therefore, here the focus is put on the 
discretisation of random variables 
in a data-driven sparse functional approximation form. In this light the optimal approximations of the state variable,
their numerical evaluations using the minimal number of sample points, as well as an efficient estimation of forward and inverse maps, i.e. the Jacobian 
of linearised forward maps, as well as Kalman gain are discussed here. 

\subsection{Random variable discretisations}

For the purpose of discretisation, the random variables appearing in \refeq{newton_gauss_markov1} 
can be expressed in terms of some known simpler kind of random variables, as previously studied by the author
and colleagues in a purely linear setting, see \cite{BRrosic12}. This can be achieved by introducing a truncated polynomial 
chaos approximation of the state variable
\begin{equation}
\label{pcex1}
x(\omega)\approx \hat{\vek{x}}(\omega)=\sum_{\alpha \in \mathcal{J}_x} \vek{x}^{(\alpha)}\Psi_\alpha(\vek{\vartheta}(\omega)),
\end{equation}
in which $\Psi_\alpha$ are multi-variate polynomials in 
random variables $\vek{\vartheta}$ as arguments. The random variables
$\vek{\vartheta}$ represent the parameterisation
of the prior uncertainties in the initial conditions or even model parameters. They are usually taken as independent,
uncorrelated random variables of some simpler kind
such as for example normal or uniform random variables corresponding to the Askey scheme as discussed in \cite{Xiu2010}. In a
similar manner, one may approximate the 
predicted error
\begin{equation}
\varepsilon(\omega)\approx \hat{\vek{\varepsilon}}(\omega)=\sum_{\alpha \in \mathcal{J}_{\varepsilon}} 
\vek{\varepsilon}^{(\alpha)}\Psi_\alpha(\vek{\eta}(\omega))
\end{equation}
in which $\vek{\eta}(\omega)$ and $\vek{\vartheta}(\omega)$ are assumed to be independent and uncorrelated.
Collecting all random variables of consideration, the global discretisation of the state reads
\begin{equation}
\label{pcex}
x(\omega)\approx \hat{\vek{x}}(\omega)=\sum_{\alpha \in \mathcal{J}_{\Psi}} \vek{x}^{(\alpha)}\Psi_\alpha(\vek{\xi}(\omega)),
\end{equation}
 in which $\vek{\xi}(\omega):=(\vek{\vartheta},\vek{\eta})$
 
\begin{figure*}
\begin{center}
\includegraphics{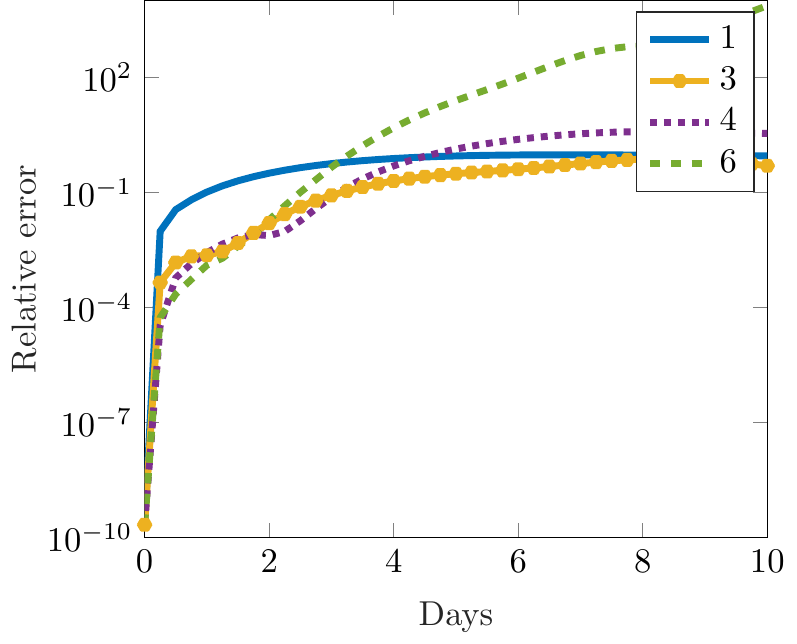}
\includegraphics{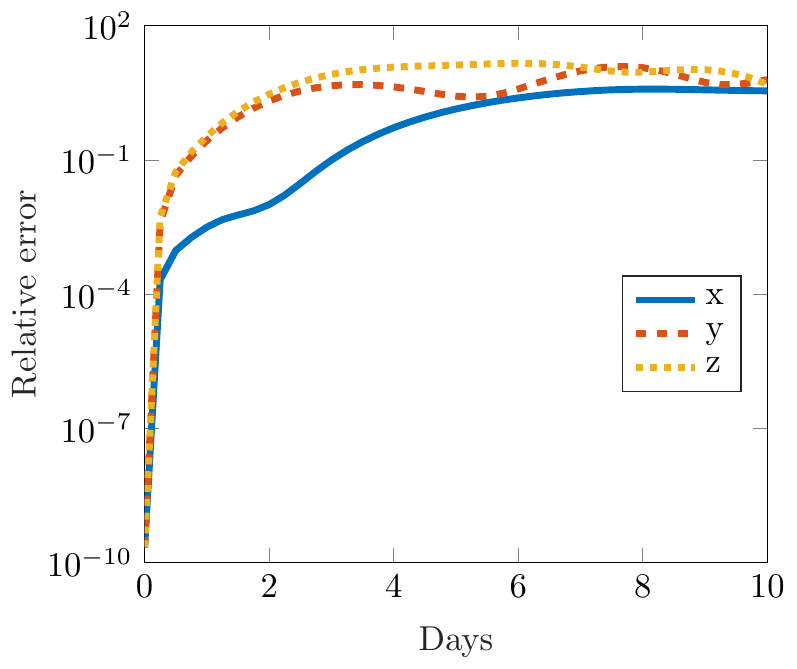}
\caption{Relative error of  a) the first state PCE w.r.t. to the polynomial order for 100 randomly chosen samples b) the state PCE for $p=4$ and
100 randomly chosen points}
\label{fig_pce_approximations}
\end{center}
\end{figure*}

 \begin{figure}
\begin{center}
\includegraphics{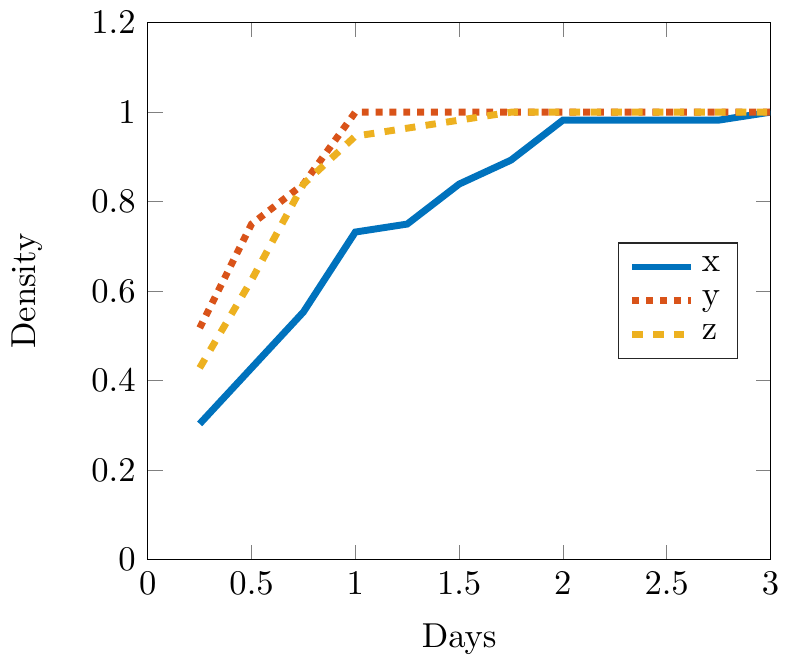}
\caption{State sparsity in time for fixed polynomial order $p=4$}
\label{fig_spar_state}
\end{center}
\end{figure}

When dealing with time-dependent systems, the 
approximation as given previously is not optimal when the time integration
of the nonlinear system before the update is too long. In such a case the state becomes highly non-Gaussian and 
requires high-order polynomial chaos approximations. \refig{fig_pce_approximations} shows the decrease of the state PCE accuracy in
time, and its improvement with the increase of the polynomial order. Similarly, the non-Gaussianity increases the number of 
sampling points neccessary for the estimation of PCE coefficients as the sparsity of
the solution decreases, see \refig{fig_spar_state}. 

To resolve this problem, the idea is to change the basis in \refeq{pcex} to
\begin{equation}
\label{regression1}
 \hat{\vek{x}}_{k}(\omega)=\sum_{\alpha \in \mathcal{J}_{\Phi}} \vek{x}_{k}^{(\alpha)}\varPhi_\alpha(\vek{\zeta}(\omega))=\vek{\Phi}_k\vek{v}_k,
\end{equation}
in which the random variable $\vek{\zeta}(\omega)$
follows the distribution of the last known state $\vek{x}_{k-n}(\omega)$ for which the lower 
order approximation in \refeq{pcex} is still suitable, and $\varPhi_\alpha(\vek{\zeta}(\omega))$ are the basis functions chosen either 
as orthogonal via a modified Gram-Schmidt process, or non-orthogonal ones as polynomial maps
of the last known state. 

The basis transformation starts with the definition of new random variables $\vek{\zeta}(\omega)$
driven by the evolution law in \refeq{initvp} such that
\begin{equation}
\vek{\zeta}(\omega)=g(\vek{\xi}(\omega))
\end{equation}
holds, in which $g(\vek{\xi}(\omega))$ describes a nonlinear transformation of the initial random variables $\vek{\xi}(\omega)$ over
some predefined period of time. Let $t_{k-n}$ be the last time moment in which the classical PCE basis can be used to approximate
the state $\vek{x}_{k-n}$. Then, given a small number
$N$ of model trajectories $(\vek{x}_{k-n}(\vek{\xi}(\omega_i))_{i=1}^N$ for $(\vek{\xi}(\omega_i))_{i=1}^N$ one may 
estimate the state coefficients $\vek{x}_{k-n}^{(\alpha)}$ in the original basis $\varPsi_\alpha(\vek{\xi}(\omega))$. 
Since $\vek{x}_{k-n}(\omega)$ is fully defined, one may take $\vek{\zeta}(\omega):=\vek{x}_{k-n}(\omega)$. By arranging $\vek{\zeta}(\omega)$
into multivariate polynomial form, we may define the new basis $\varPhi_\alpha(\vek{\zeta}(\omega))$ using 
the modified Gram-Schmidt (MGS) orthogonalisation
process, for more details please see \cite{Gerritsma:2010}.  
%
%
%
In such a case the new state $\vek{x}_k(\omega)$ at time $t_k$ can be estimated given a small
number of trajectories $(\vek{x}_k(\vek{\xi}(\omega_i))_{i=1}^N$ and their corresponding basis functions $\vek{\Phi}(\vek{\zeta}(\omega_i))$.
Having
\begin{eqnarray}
\label{regressiona}
 \hat{\vek{x}}_k(\omega_i)&=&\sum_{\alpha \in \mathcal{J}_\Phi} \vek{x}_k^{(\alpha)}\varPhi_\alpha(\vek{\zeta}(\omega_i))\nonumber\\
 &=&
 \sum_{\alpha \in \mathcal{J}_\Phi}
 \vek{x}_k^{(\alpha)}\varPhi_\alpha(g(\vek{\xi}(\omega_i)))
\end{eqnarray}
one may estimate the coefficients $\vek{x}_k^{(\alpha)}$ via Bayesian regression as described in Section \ref{spaopt}. Here,
$\mathcal{J}_\Phi$ is a new multi-index set defined by a polynomial order that is lower than the corresponding Hermite one.
This procedure further allows the evaluation of a large number of samples of $\vek{x}_k$ as the large
number of samples of $\vek{\zeta}$ resp. $\vek{\xi}$ is known, and hence one may repeat the process to estimate the next unknown state in time $t_{k+1}$. 

\begin{figure}
\begin{center}
\includegraphics{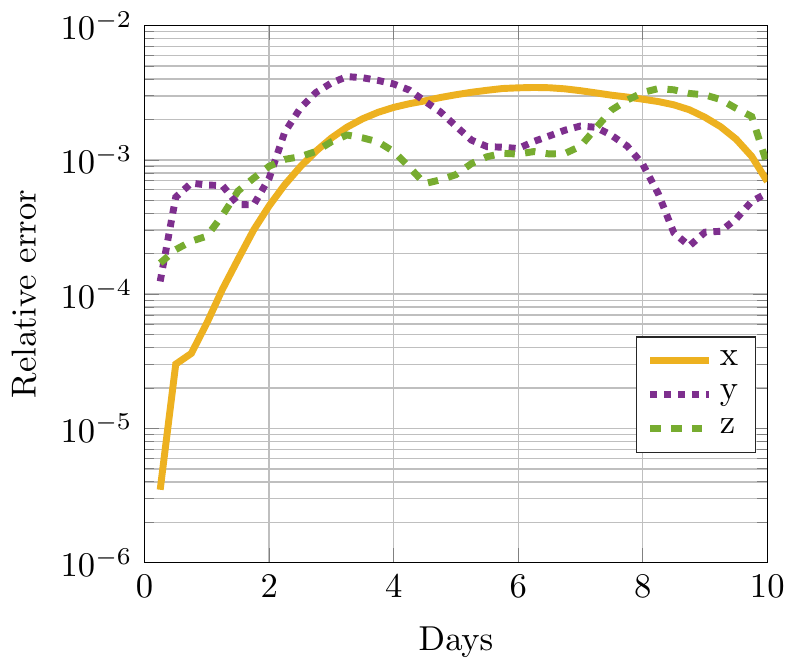}
\caption{Accuracy of the MGS basis in time for all three Lorenz states}
\label{mgs_state_accuracy}
\end{center}
\end{figure}

\begin{figure}
\begin{center}
\includegraphics{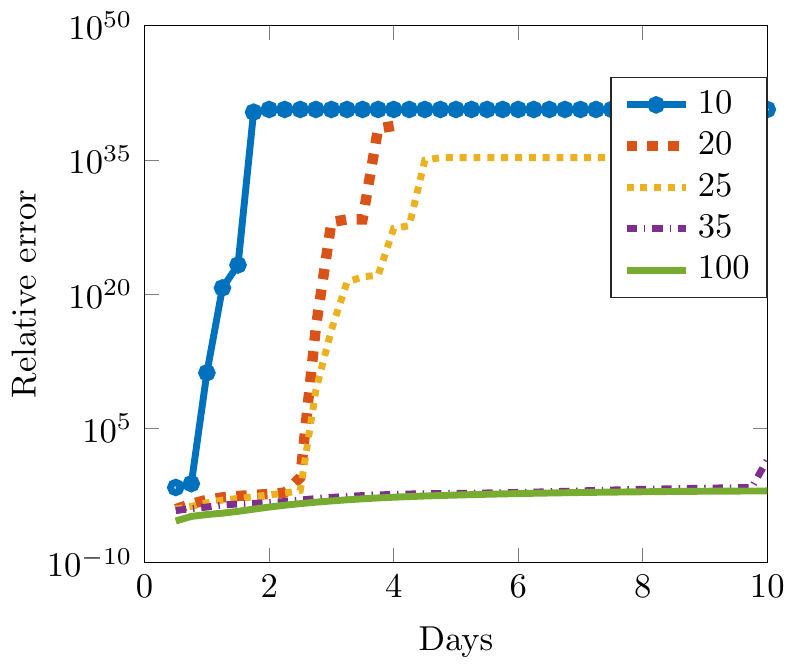}
\caption{Accuracy of the MGS solution w.r.t. the number of randomly chosen samples for the first Lorenz state for $p=3$}
\label{mgs_state_accuracy_mc_samples}
\end{center}
\end{figure}

\begin{figure}
\begin{center}
\includegraphics{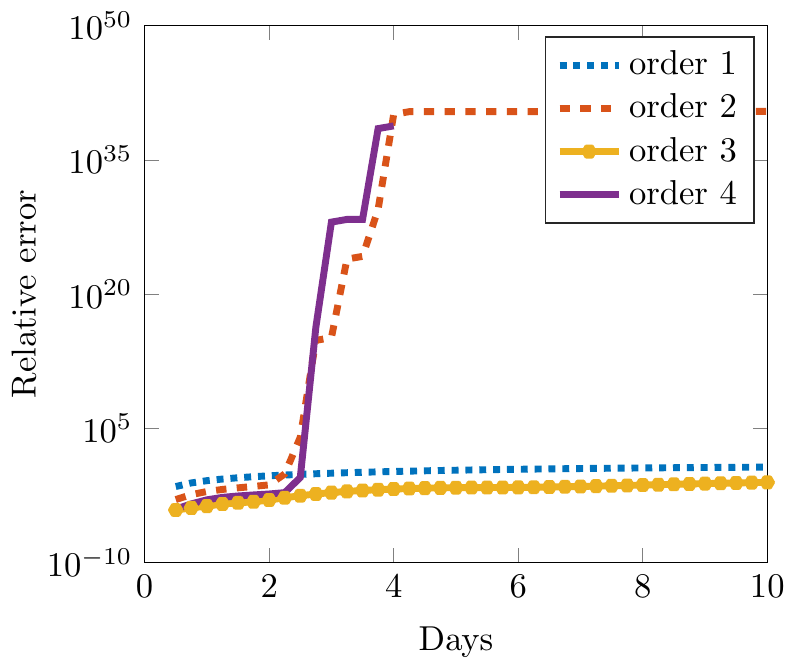}
\caption{Accuracy of the MGS solution w.r.t. the polynomial order}
\label{mgs_state_accuracy_order}
\end{center}
\end{figure}

\begin{figure}
\begin{center}
\includegraphics{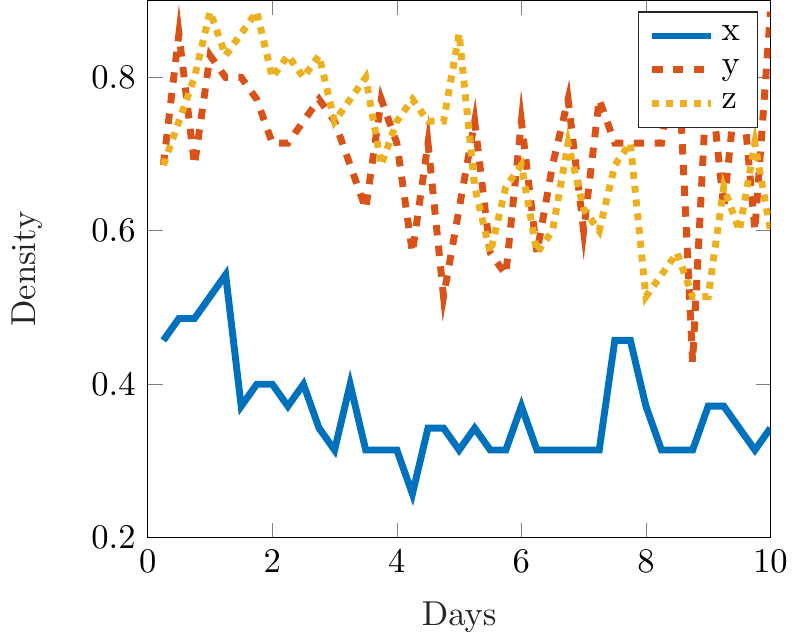}
\caption{Sparsity of the MGS Lorenz state in time}
\label{fig_mgs_sparsity}
\end{center}
\end{figure}

\refig{mgs_state_accuracy} shows the accuracy of the MGS for the polynomial order $p=3$ 
and $100$ randomly chosen samples w.r.t. the solution obtained from $10^6$ Monte Carlo runs. In comparison to the 
Bayesian regression on classical PCE depicted in \refig{fig_pce_approximations} one may note that the accuracy 
of the MGS solution improves 
by an order of magnitude for the same number of samples. The dependence of the MGS solution on the number of samples 
and
the polynomial order can be seen in 
\refig{mgs_state_accuracy_mc_samples} and \refig{mgs_state_accuracy_order}, respectively. As expected, the accuracy improves with the 
sample number. Similar holds for the polynomial order. Finally, the sparsity of the newly obtained approximation is shown in 
\refig{fig_mgs_sparsity}, where it is observed that the first state is much sparser than the other two.

The basis estimation via Gram-Schmidt orthogonalisation can be computationally demanding. Thus, 
a much more efficient solution is to consider the non-orthogonal basis. The simplest choice is 
to observe the
current state $\vek{x}_{k}(\omega)$ as a nonlinear map of the last known one $\vek{x}_{k-n}(\omega)$, i.e.~
\begin{equation}\label{eq:np}
 \vek{x}_{k}(\omega)=\sum_{\alpha \in \mathcal{J}_{\varUpsilon}} \vek{x}_{k}^{(\alpha)} \varUpsilon_\alpha(\vek{x}_{k-n}(\omega))
\end{equation}
in which the coefficients $\vek{x}_{k}^{(\alpha)}$ are obtained via regression described in Section \ref{spaopt}. 
Here, $\varUpsilon_\alpha(\vek{x}_{k-n}(\omega))$ are taken to be the non-orthogonal multivariate polynomials defined as:
\begin{eqnarray}
\varUpsilon_\alpha(\vek{x}_{k-n}(\omega))&=&\vek{x}_{k-n}^{(\alpha)}\nonumber\\
&=&x_{k-n}^{\alpha_1}y_{k-n}^{\alpha_2}z_{k-n}^{\alpha_3}
\end{eqnarray}
with $(\alpha)$ being the multi-index set similarly defined to the one that describes the classical PCE. 

\begin{figure*}
\begin{center}
\begin{minipage}{\textwidth}
 \begin{minipage}{0.5\textwidth}
 \includegraphics{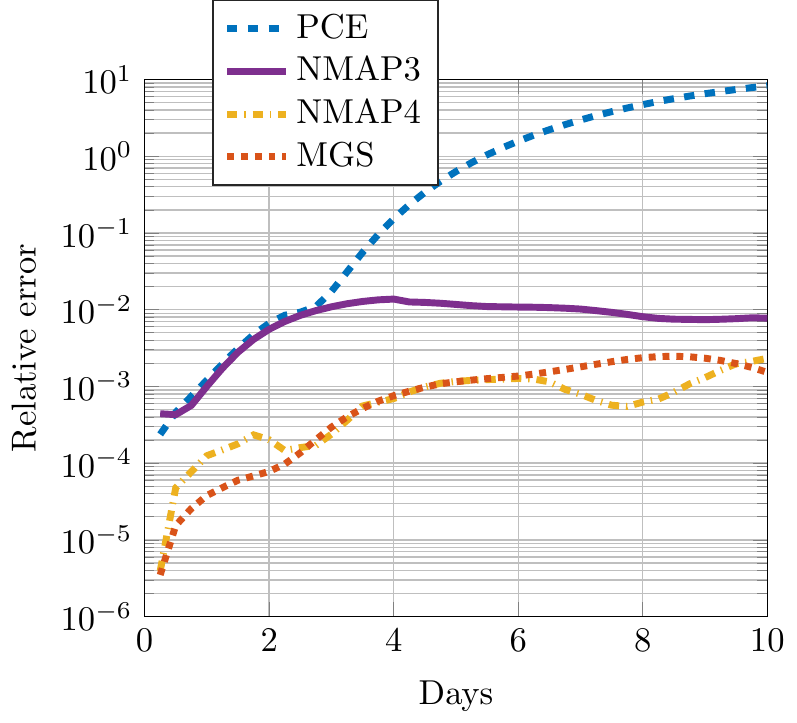}\\
 a) Comparison of accuracy
 \end{minipage}
 \begin{minipage}{0.5\textwidth}
 \vskip 30pt
 \includegraphics{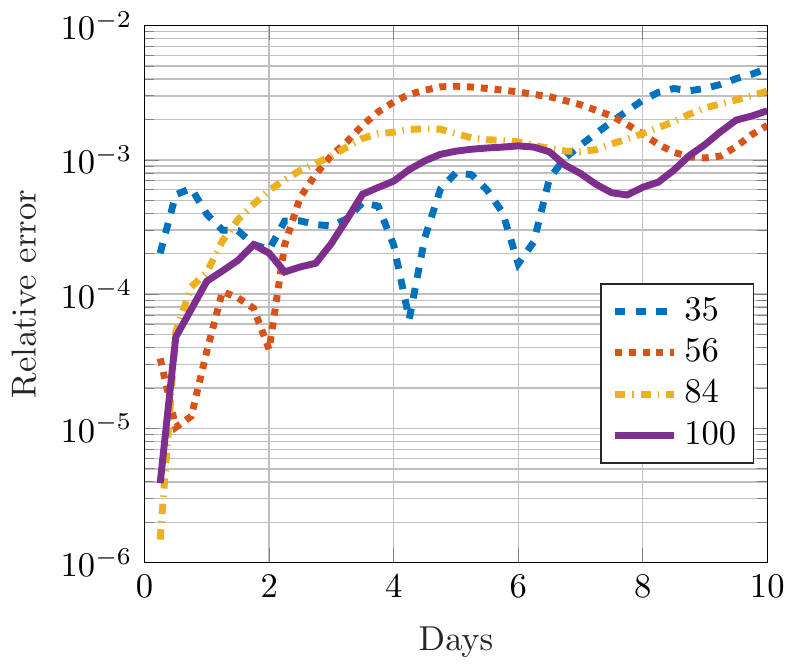}\\
 b) Accuracy w.r.t. to the number of training points used in regression
 \end{minipage}
 \begin{minipage}{\textwidth}
 \vskip 30pt
 \begin{center}
 \includegraphics{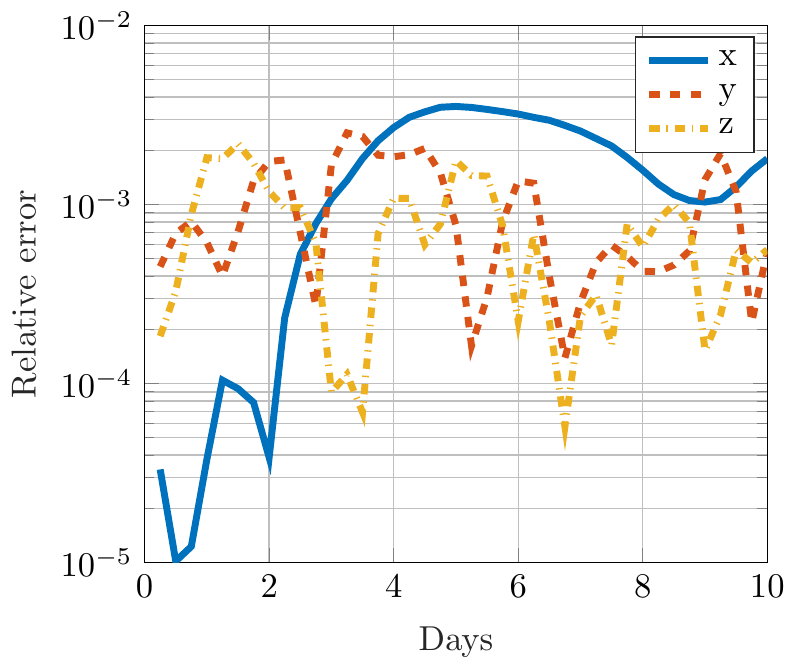}\\
 c) State accuracy
 \end{center}
 \end{minipage}
\end{minipage}
%
\caption{Non-orthogonal approximation of solution}
\label{fig_nmap_prop}
\end{center}
\end{figure*}

\refig{fig_nmap_prop}a) shows the accuracy of the third and fourth order nonlinear map (NMAP) approximations (of order 3
(NMAP3) resp. order 4 NMAP4) 
compared to the solution obtained by regressing on the fixed Hermite polynomial basis of fourth order (PCE), and 
the MGS solution of fourth order. While the PCE solution is not accurate enough, both the MGS and 
the nonlinear map solutions give similar 
results for the same order of approximation. In the beginning lower order nonlinear map solution (NMAP3) matches the solution 
obtained by fixed regression. 
In contrast 
to the PCE solution, the error stabilises over time and does not over-estimate the Lorenz state. Furthermore, the accuracy
of the NMAP4 solution 
is tested on different data set sizes 
in \refig{fig_nmap_prop}b). The experiment shows that even a low number of samples (56 samples)
can be used to achieve the desired accuracy, see \refig{fig_nmap_prop}c).

By using approximations in \refeq{regressiona} and in \refeq{eq:np} one may use a small number 
of solution trajectories of $\vek{x}_{k-n}(\omega)$ to estimate a large number of samples
$\vek{x}_{k}(\omega)$. The approximations are made adaptive such that the last known basis is used 
in a current time, and the Kullback-Leibler divergence is used to estimate the error compared to the validation set. 
If the error is bigger than tolerance then the basis is adaptively modified. 

Even though
both of the previous approximations are significantly better than the original basis, they are not suitable 
to be used in the filtering process due to correlated arguments, and in the latter case also due to the non-orthogonality. 
Therefore, to compute the time dependent polynomial chaos approximations, the previous approximations at the update time 
are transformed such that the non-Gaussian correlated random variables $\zeta(\omega)$ are mapped to 
uncorrelated Gaussian ones via nonlinear transformation. 
The main idea of the transformation process is to map the state variable 
$\vek{x}_{k-n}(\omega_x)$
by an isoprobabilistic map to a Gaussian random variable $\theta(\omega_\theta), \omega_ \theta \in \Omega_\theta$, i.e.~
\begin{equation}
\label{here_for_trans}
 T:\quad \vek{x}_{k-n}(\omega_x) \mapsto \theta(\omega_\theta)
\end{equation}
such that the approximations 
rewrite to the PCE with multivariate Hermite orthogonal basis $\Psi_\alpha(\vek{\theta}(\omega_\theta))$:
\begin{eqnarray}
\label{refrvestim_ipce_trans}
  \breve{\vek{x}}_{k}(\omega_\theta)=\sum_{\alpha \in \mathcal{K}} x_{k}^{(\alpha)}\Psi_\alpha(\vek{\theta}(\omega_\theta))
\end{eqnarray}
characterised by much lower cardinality than the one in \refeq{regressiona} or \refeq{eq:np}.

Due to simplicity reasons, the transformation in \refeq{here_for_trans} is assumed to be of the Nataf-type, 
which shows good performance for this kind of problem. The other more general type of transformations
are the current state of the research and will be discussed in another paper. 

The Nataf transform is a composition of maps $T=T_1\circ T_2$ in which the first one $T_1$
maps the vector of non-Gaussian random variables $\vek{\zeta}$ with marginal cumulative distributions $F_{\vek{\zeta}}$
to the vector of correlated standard Gaussian 
variables $\kappa$ via inverse cumulative distribution of the standard normal $\varPhi_{\mathcal{N}}^{-1}$:
\begin{equation}\label{eq:ss}
 T_1:\quad \vek{\zeta}(\omega_\zeta) \rightarrow \vek{\kappa}(\omega_\theta) :=\varPhi^{-1} (F_{\vek{\zeta}}(\vek{\zeta})),
\end{equation}
whereas the second one $T_2$ maps correlated random variables into uncorrelated ones
\begin{equation}
 T_2:\quad \vek{\kappa}(\omega_\theta) \rightarrow \vek{\theta}(\omega_\theta) =\vek{C}_{\kappa}^{-1/2}\vek{\kappa}(\omega_\theta).
\end{equation}
Here, the factor $\vek{C}_{\kappa}^{-1/2}$ is evaluated using the Cholesky decomposition. 
To perform the step in \refeq{eq:ss}, one requires knowledge on the cumulative distribution function (cdf) $F_{\vek{\zeta}}$. As this information 
is not accessible, 
but only instances of the random variable $(\vek{\zeta}(\omega_j))_{j=1}^M$ are known, one may use
the kernel density estimator as the one presented in 
\cite{botev2010} to obtain $F_{\vek{\zeta}}$. In addition, $F_{\vek{\zeta}}$ is interpolated in a Bayesian manner (see Section 8)
using the polynomial of order 3.

Hence, for the further process of assimilation
one may rewrite \refeq{eq:np} to the orthogonal polynomial chaos expansion expressed 
in terms of newly estimated standard random variables:
\begin{eqnarray}
\label{refrvestim_ipce}
   \hat{\vek{x}}_{k-\Delta \tau,a}^{(i+1)}={\hat{\vek{x}}_{k-\Delta \tau,f}}+
   \hat{\vek{K}}_{k-\Delta \tau}^{(i)}(\hat{\vek{x}}_{k,a}-\hat{\vek{y}}_{kh}^{(i)})
\end{eqnarray}
in which 
\begin{equation}
\label{pcepredm}
y_{kh}(\omega)\approx \hat{\vek{y}}_{kh}(\omega)=\sum_{\alpha \in \mathcal{K}} 
\vek{y}_{k\ell}^{(\alpha)}\Psi_\alpha(\vek{\theta}(\omega))+\hat{\varepsilon}(\omega)
\end{equation}
i.e.
\begin{equation}
 \hat{\vek{y}}_{kh}(\omega)=\hat{\mathring{H}}(\hat{x}_{k-\Delta \tau,f}-\mathring{x})+\hat{\mathring{h}}+\hat{\varepsilon}.
\end{equation}

\begin{figure*}
\begin{center}
\includegraphics{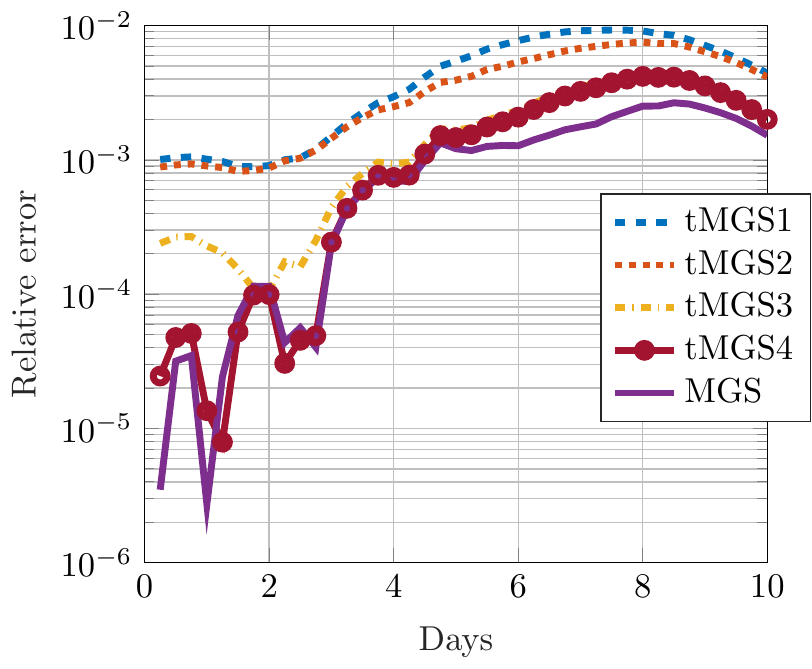}
\includegraphics{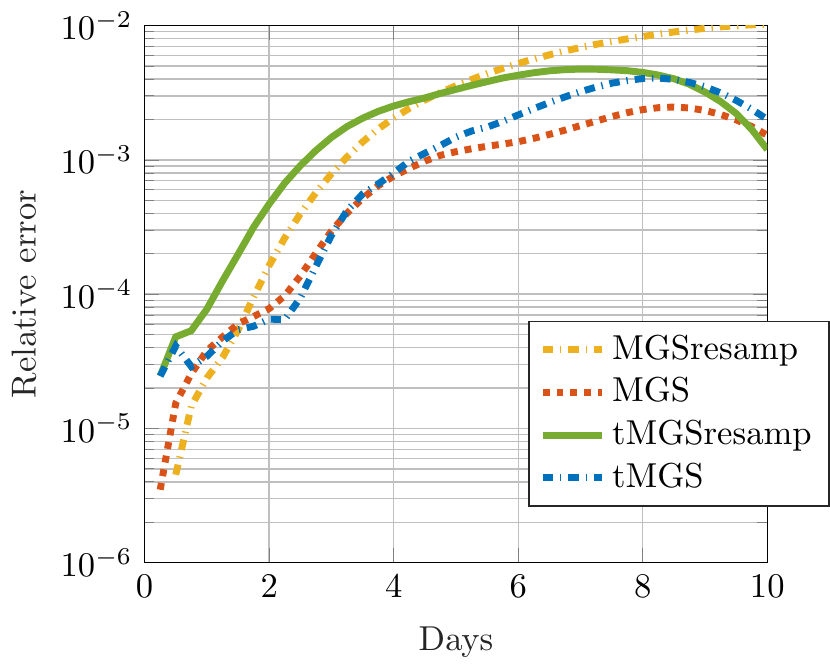}
\caption{a) Accuracy of transformed MGS solution in time w.r.t. to polynomial order b) Comparison of accuracies of different MGS approaches}
\label{fig_comparison_mgs}
\end{center}
\end{figure*}

The accuracy of the transformed solution in a
Gaussian basis (tMGS- transformed modified Gram-Schmidt process)
compared to non-Gaussian ones (denoted by MGS in plot) w.r.t. to the polynomial order is 
shown in \refig{fig_comparison_mgs}a). 
As expected, the Gaussian basis requires higher polynomial order to achieve the same accuracy as the non-Gaussian one.

The comparison of the transformed approach to the classical MGS one is depicted in Fig.~(\ref{fig_comparison_mgs}b). Here, four different types of 
solutions are considered. The solutions denoted by MGS and tMGS (transformed MGS) 
are obtained by integrating original samples of the state in time,
whereas solutions MGSresamp and tMGSresamp
are obtained by sampling 
the polynomial chaos approximations that are further integrated in time. In the latter case the approximation error gets propagated in time, and hence the solution 
is less accurate than the corresponding sampled solution. The reason to investigate the second case lies in the updating procedure. After the update 
of the state
is made one does not have the original state samples coming from sampling the initial condition. 
Instead, one samples the newly obtained polynomial chaos approximation.

%
%

\begin{figure*}
\begin{center}
\includegraphics{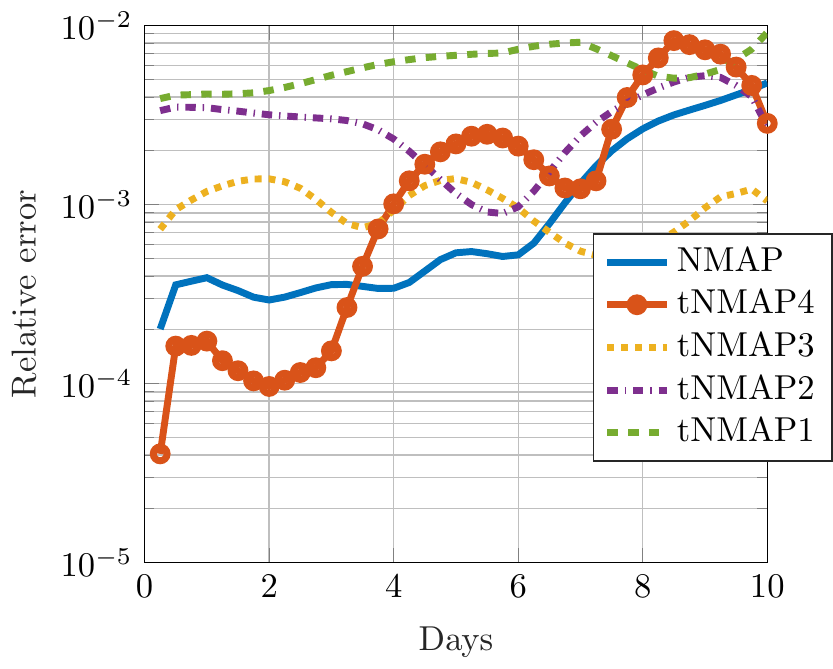}
\includegraphics{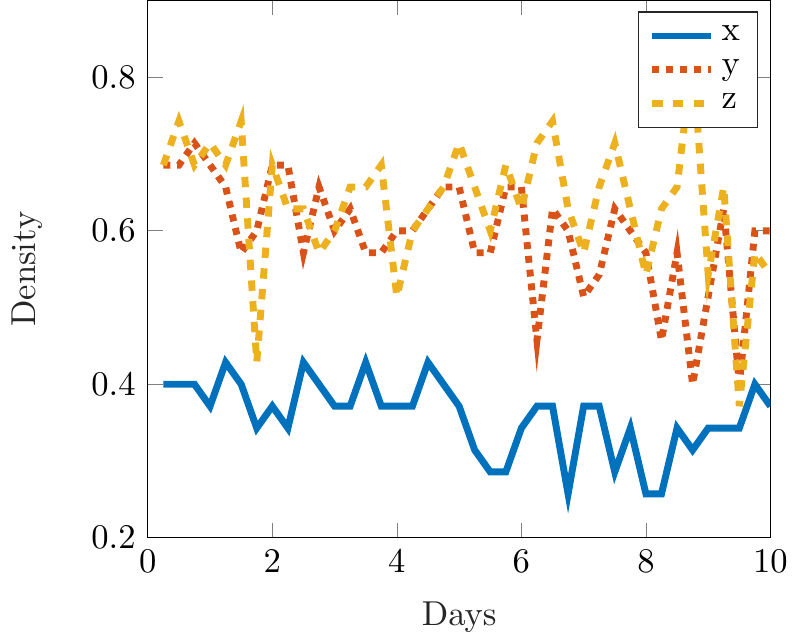}
\caption{State sparsity of the non-orthogonal approximation of solution}
\label{fig_non_orth}
\end{center}
\end{figure*}

The discretised posterior in \refeq{refrvestim_ipce} is described by both the state random variables as well as the
variables describing the measurement noise. 
The number of the latter ones increases with the number of measurements, and hence the cardinality of the posterior PCE grows. 
However, the dimension increase can be avoided by same transformation process as described before.  
%
%
In \refig{fig_non_orth}a) the accuracy of the transformed state for NMAP estimate with respect to the polynomial order is depicted. 
The sparsity of the newly obtained state is depicted in \refig{fig_non_orth}b) and is slightly higher than the one described by the MGS procedure. 

\begin{figure*}
\begin{center}
\includegraphics{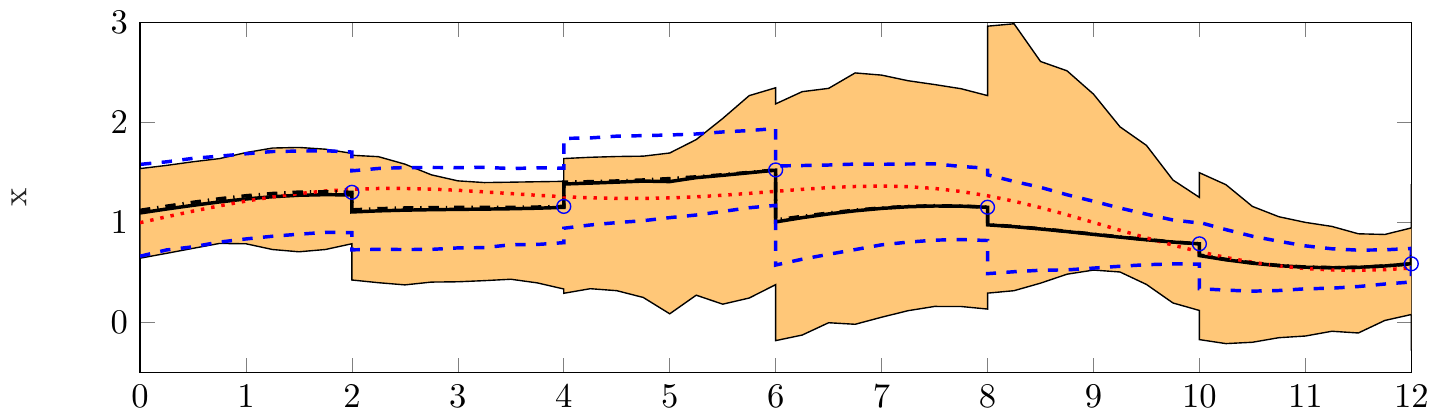}\\
\includegraphics{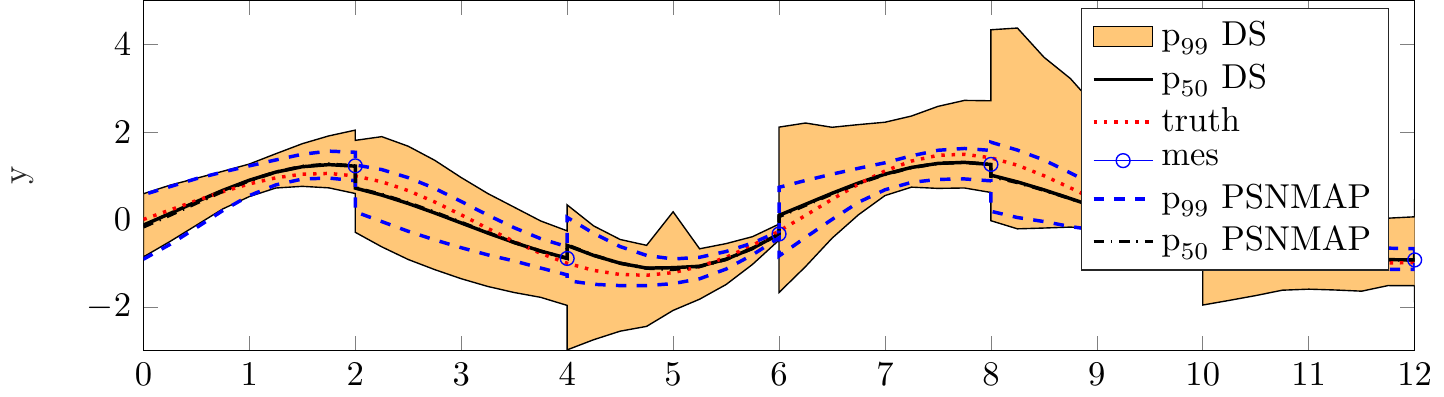}\\
\includegraphics{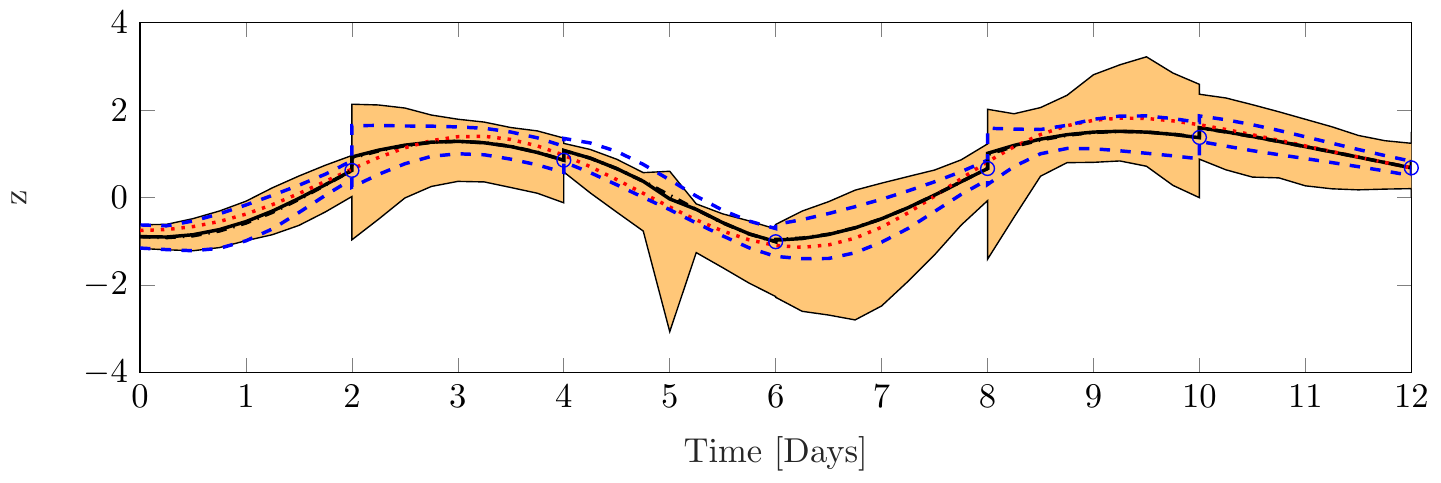}
\caption{Update of the Lorenz 1984 state backwards every two days. The forecast is estimated using the nonlinear map.}
\label{fig_filter_map_seq}
\end{center}
\end{figure*}

\begin{figure*}
\begin{center}
\includegraphics{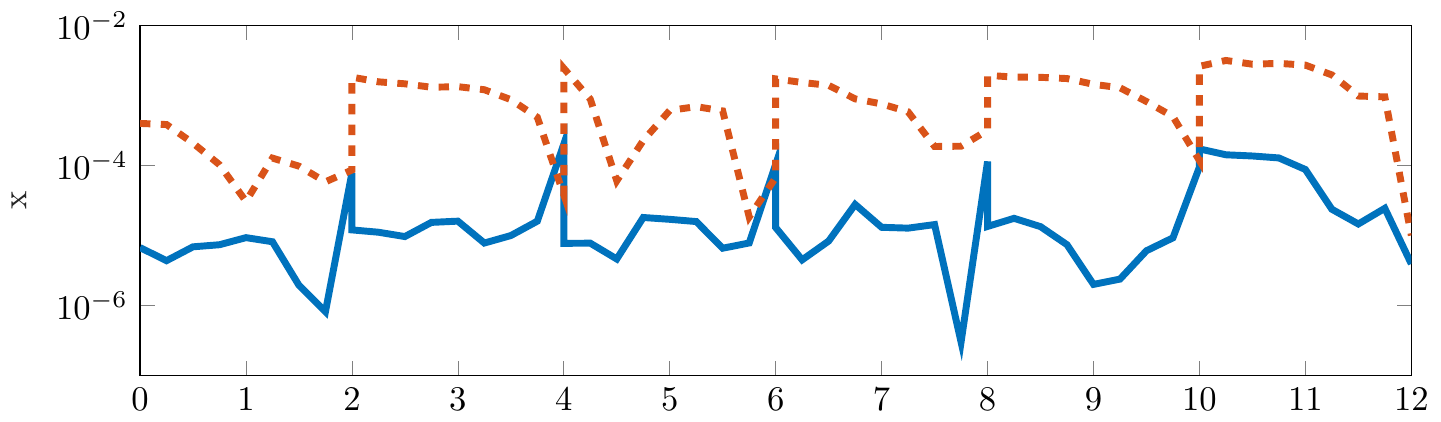}\\
\includegraphics{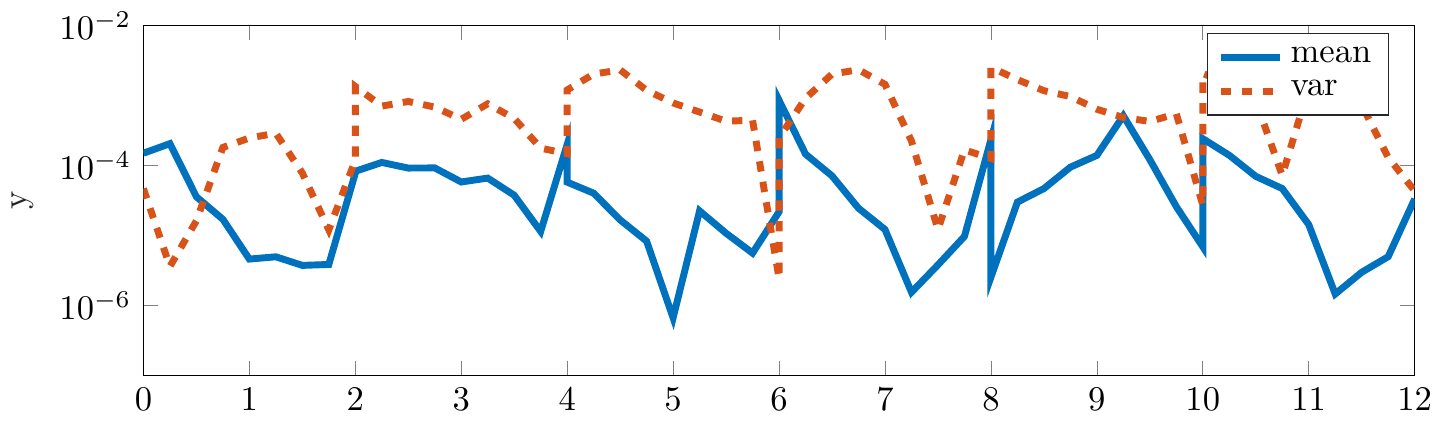}\\
\includegraphics{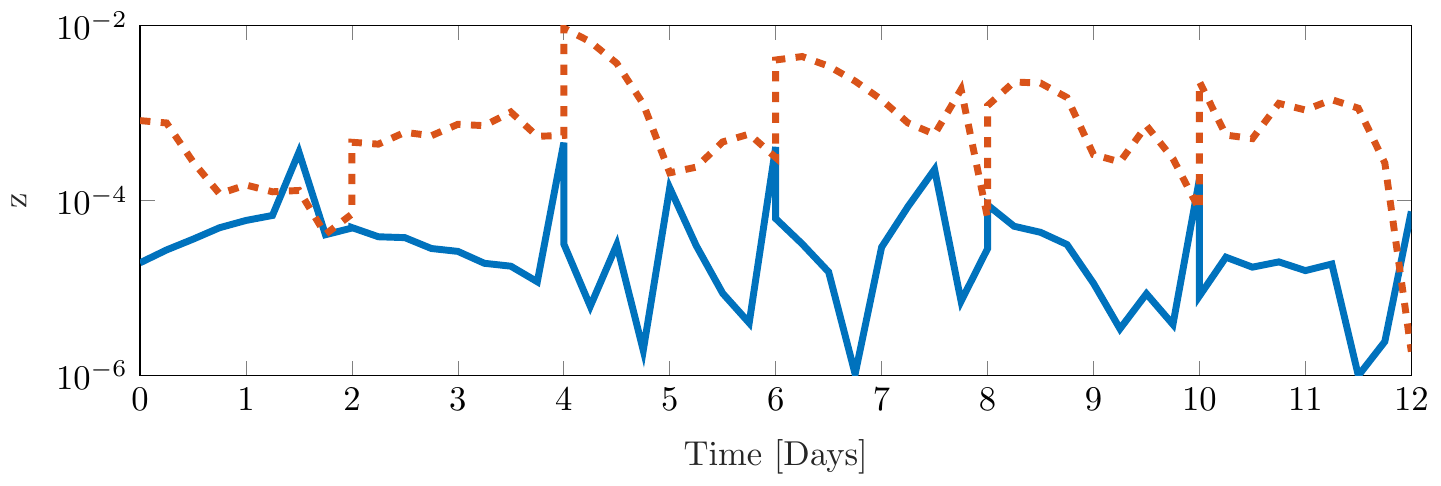}
\caption{Relative error between the posterior mean and variance of the transformed MGS based solution with respect to the transformed nonlinear map one. Both are obtained by using
the square root algorithm}
\label{fig_filter_map_mgs}
\end{center}
\end{figure*}

Finally, the assimilation results can be significantly different than those obtained using the classical PCE. 
In \refig{fig_filter_map_seq} one may see a comparison of the assimilated state using the direct iteration with non-adaptive 
classical polynomial chaos expansion of order 4 (DS)
and the pseudo-state (PS) update (frequency of update is 6h) using the transformed nonlinear map estimate of same order. 
Clearly, the DS estimation leads
to the overestimation of the posterior variance already after one day of estimation, as expected. This is due to the inaccuracy of the state approximations.
On the other hand, the transformed nonlinear map estimate and the
one based on the transformed modified
Gram-Schmidt estimation are giving very close results. \refig{fig_filter_map_mgs} depicts the mean
and variance relative errors between these two solutions.

\begin{figure*}
\begin{center}
\includegraphics{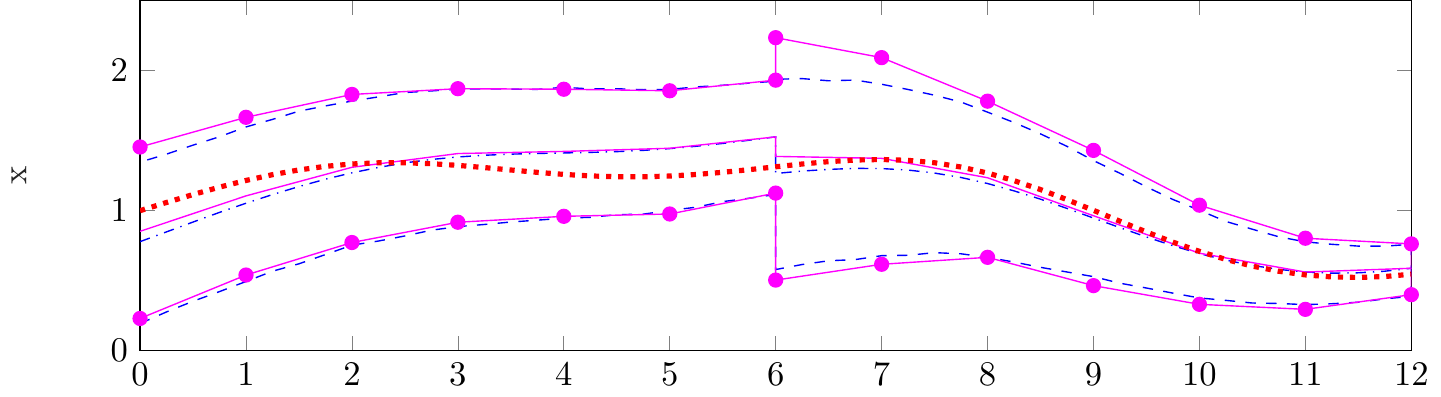}\\
\includegraphics{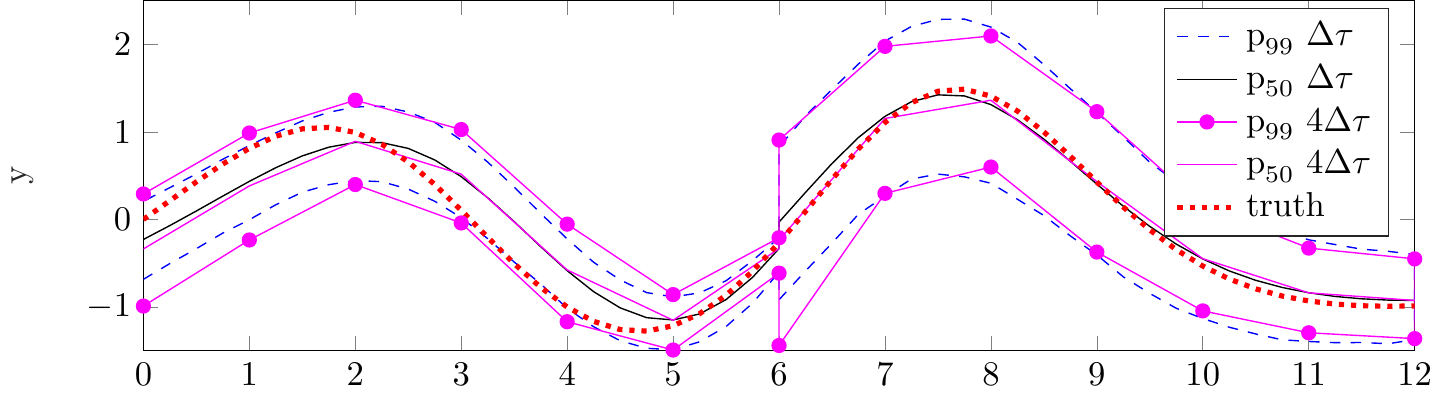}\\
\includegraphics{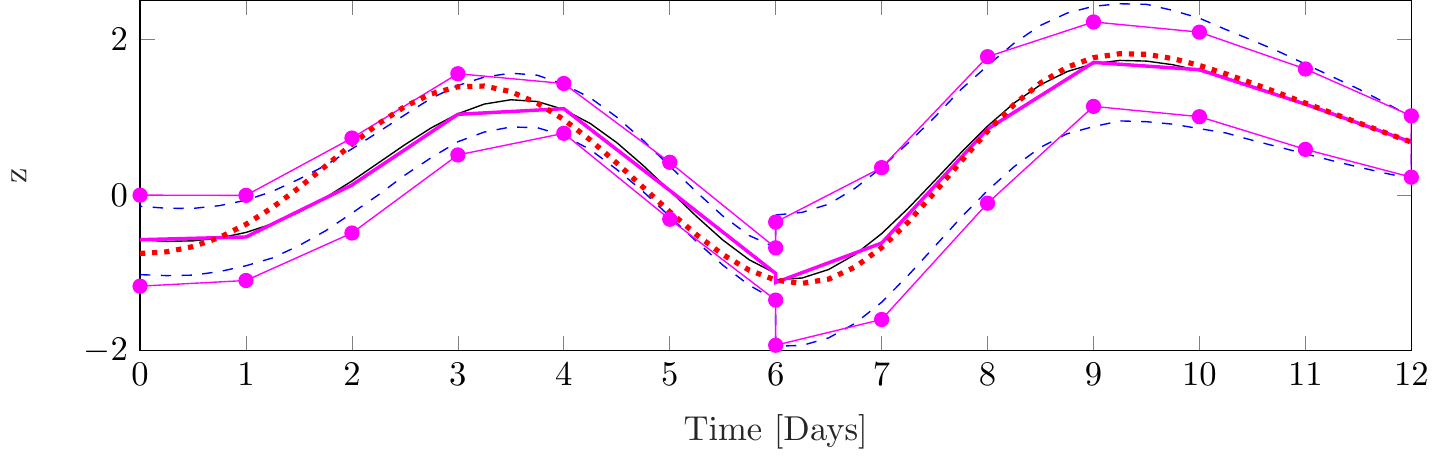}
\caption{The robustness of the square-root update of the full state with respect to the time step size}
\label{fig_filter_all_states_diff_steps}
\end{center}
\end{figure*}

\begin{figure*}
\begin{center}
\includegraphics{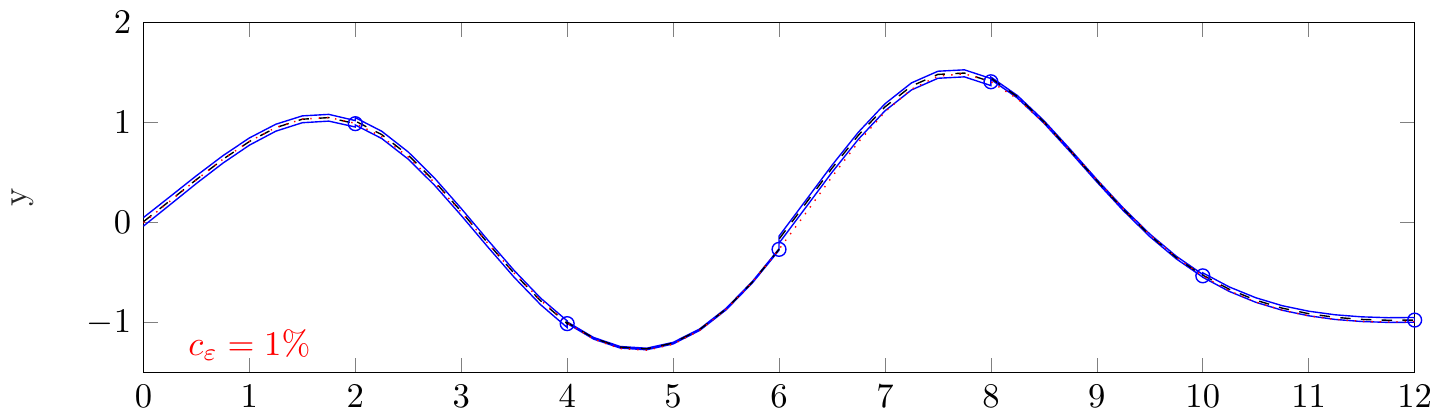}\\
\includegraphics{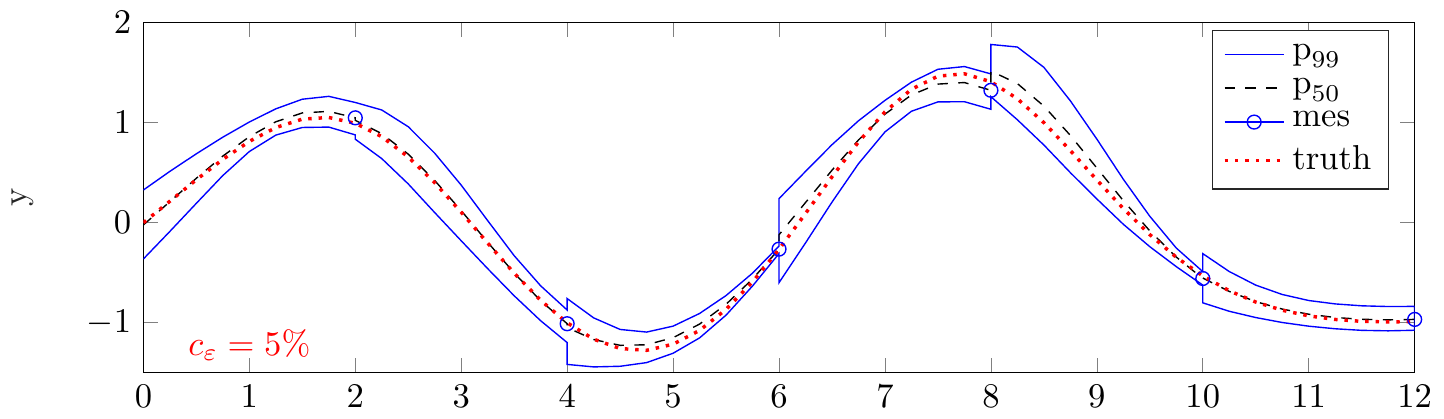}\\
\includegraphics{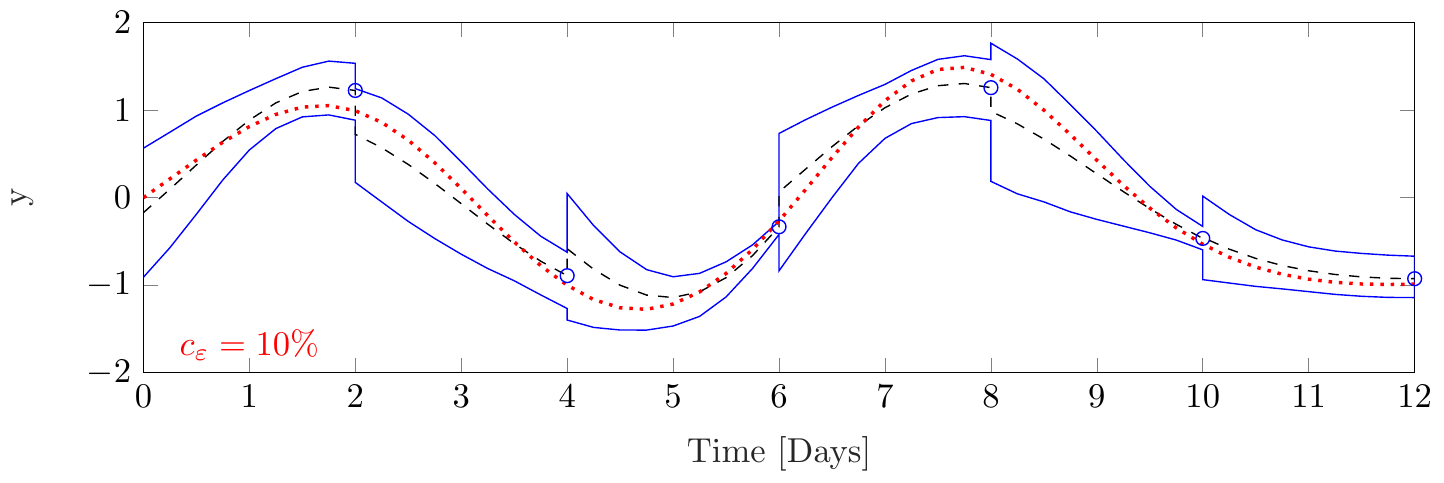}
\caption{Estimation of the Lorenz state with respect to different coefficients of the variation of noise $c_{\varepsilon}$}
\label{fig_filter_noise}
\end{center}
\end{figure*}

The modified basis approach results in a stable posterior variance with respect to the 
pseudo-time step size, see \refig{fig_filter_all_states_diff_steps} in which are depicted 
posterior bounds for two different updating step sizes. 
This, however, does not hold for the classical PCE. 
In addition, the modified updating procedure is robust with respect to the measurement
noise as presented in \refig{fig_filter_noise}. Here, the posterior $99\%$ regions are shown for three different values of the measurement noise
with $c_{\varepsilon}$ being the coefficient of the variation of noise. 

%

%
%

\subsection{Sparse polynomial chaos approximations}\label{spaopt}

The Gauss-Newton-Markov-Kalman filter requires repeated evaluations of the forward problem. 
To reduce the overall computational burden,
 the propagation of the uncertainty through the forward problem can be achieved in 
a data-driven non-intrusive spectral setting.
Given the approximation of the state in a polynomial chaos setting
\begin{equation}
\label{regression}
 \hat{\vek{x}}_f(\omega)=\sum_{\alpha \in \mathcal{J}_\Psi} \vek{x}_f^{(\alpha)}\Psi_\alpha(\vek{\xi}(\omega))=\mat{\Psi}\vek{v},
\end{equation}
the goal is to estimate the unknown coefficients 
$\vek{v}$ given $N$ samples $(\vek{x}_f(\omega_i))_{i=1}^N$, i.e.~
\begin{equation}
 \vek{x}_f(\omega_i)=\sum_{\alpha \in \mathcal{J}_\Psi} \vek{x}_f^{(\alpha)}\Psi_\alpha(\vek{\xi}(\omega_i))
\end{equation}
for $i=1,...,N$.
 In a vector form the previous relation reads
\begin{equation}
 \label{reg}
 \vek{u}=\mat{\Psi} \vek{v}.
\end{equation}
In a case when $N\leq P:=\textrm{card } \mathcal{J}_\Psi$, the system in \refeq{reg}
is undetermined, and requires additional information. As a priori knowledge on the 
current state exists (e.g. for small time step sizes 
the subsequent state is close to the current one), one may model 
the unknown coefficients $\vek{v}$ a priori as random variables
in $L_2(\varOmega_v,\mathcal{F},\mathbb{P};\mathbb{R}^P)$, i.e.
\begin{center}
 $\vek{v}(\omega_v):=[v^{(\alpha)}(\omega_v)]: \varOmega_v \rightarrow \mathbb{R}^P,$
\end{center}
and further use the linear Gauss-Markov-Kalman filtering procedure as previously described
to determine their conditional mean. 
%
%
%
As the coefficients can be both positive and negative, one may assume that the prior $\vek{v}(\omega_v)$ is normally distributed 
\begin{center}
$ \vek{v} \sim \mathcal{N}(\vek{0},\vek{I}) \sim  e^{-{\frac {\|\vek{v} \|^{2}}{2}}} $
\end{center}
resulting in 
\begin{equation}
\vek{v}_a(\omega_v)=\vek{v}_f(\omega_v)+\vek{K}(\vek{u}-\vek{\Psi} \vek{v}_f(\omega_v)).
\end{equation}
Having the Kalman filter on both the updating and forecasting levels, the last equation together with \refeq{smoothing_filter1_gen}
forms the hierarchical structure of the iterative Gauss-Newton-Kalman filter. 
However, such estimation still requires a large number of samples as all coefficients in the polynomial chaos expansion 
are taken into consideration even those close to zero. To promote for sparsity, see \refig{fig_spar_state},
the prior distribution has to be concentrated more around the zero value. This can be achieved by 
taking a Laplace prior
\begin{center}
$ \vek{v} \sim e^{-{\|\vek{v}\|_1}} $
\end{center}
to model the unknown coefficients. As the work with a Laplace prior is computationally difficult, in this paper we use
the corresponding hyperprior instead as advocated in relevance vector machine approach, see \cite{Tipping}. 
The hyperprior is modelled as
$$p(\vek{v}|\varpi)=\prod_{\alpha \in \mathcal{J}} p(\vek{v}^{(\alpha)}|\varpi_\alpha), \quad \vek{v}^{(\alpha)}\sim \mathcal{N}(0,\varpi_\alpha^{-1})$$
with $\varpi_\alpha$ being the precision of each PCE coefficient modelled by a Gamma prior $p(\varpi_\alpha)$.
By marginalising over $\varpi$ one obtains the overall prior
 $$p(\vek{v})=\prod_{\alpha \in \mathcal{J}} \int   p(\vek{v}^{(\alpha)}|\varpi_\alpha) p(\varpi_\alpha) \textrm{d}\varpi_\alpha$$
 which is further simplified by taking the most probable values for $\vek{\varpi}$, i.e.~$\vek{\varpi}_{MP}$.
 

\begin{figure}
\begin{center}
\includegraphics{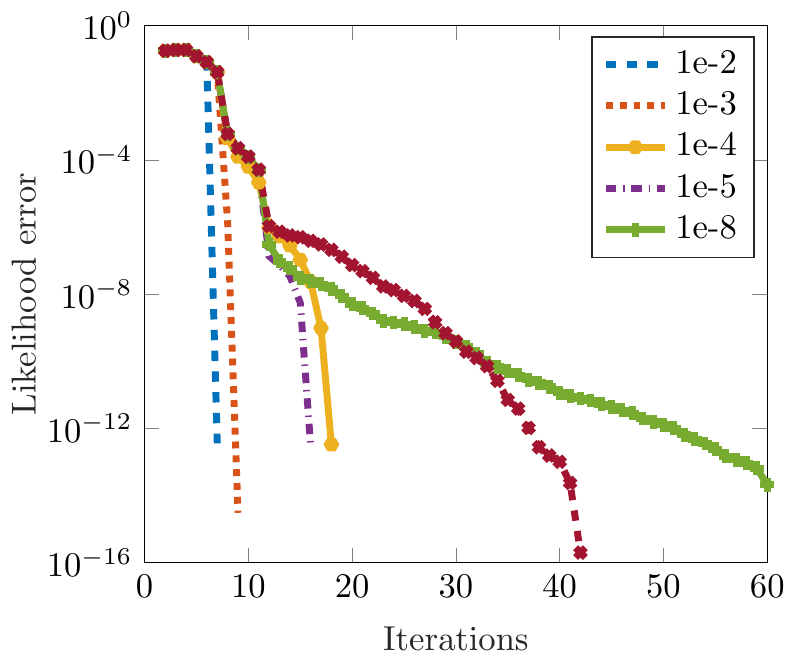}
\caption{Convergence of marginal likelihood error with respect to the priorly assumed 
regression error $\mathcal{N}(0,\sigma^2)$ for the Lorenz 1984 example}
\label{fig_like_err}
\end{center}
\end{figure}

\begin{figure}
\begin{center}
\includegraphics{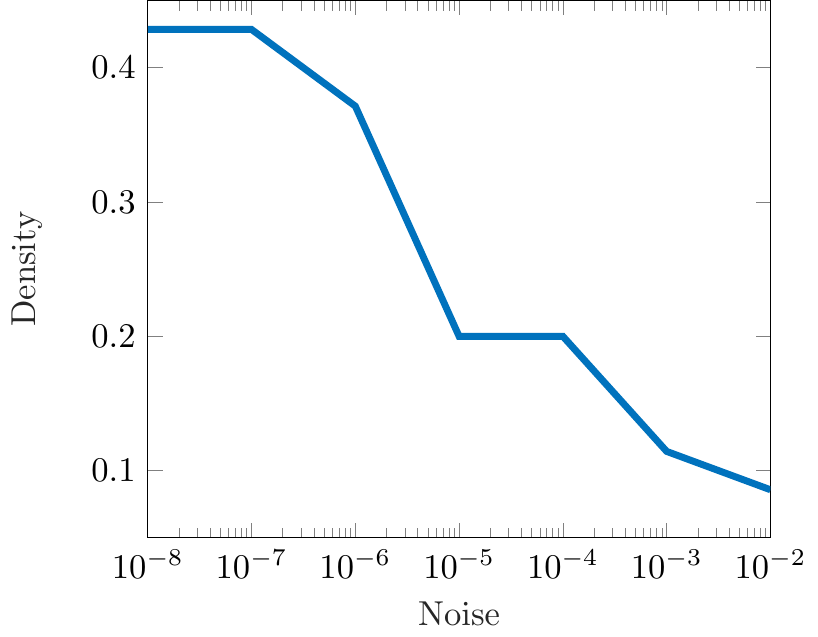}
 \caption{Lorenz 1984 state sparsity with respect to the data noise level}
 \label{fig_sparse_noise}
\end{center}
\end{figure}

\begin{figure}
\begin{center}
\includegraphics{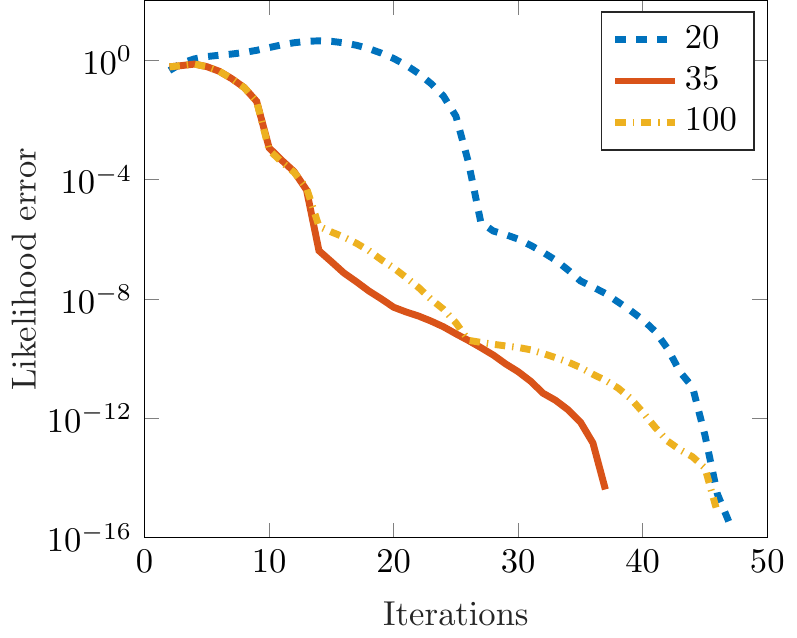}
\caption{Convergence of the marginal likelihood error with 
respect to the number of data points for a fixed regression error prior. The state is the second component of the Lorenz 1984 system.}
\label{fig_like_err_samp}
\end{center}
\end{figure}


 To estimate the coefficients in \refeq{reg} we further use Bayes's rule in a form 
 \begin{equation}\label{bayes_sparse}
p(\vek{v},\vek{\varpi},\vek{\sigma}^2|\vek{u})=\frac{p(\vek{u}|\vek{v},\vek{\varpi},\vek{\sigma}^2)}{p(\vek{u})} p(\vek{v},\vek{\varpi},\vek{\sigma}^2)
\end{equation}
 in which the coefficients $\vek{v}$, the precision $\vek{\varpi}$ and the regression error $\vek{\sigma}^2$ are assumed to be uncertain.
For computational reasons the posterior is further factorised into
$$p(\vek{v},\vek{\varpi},\vek{\sigma}^2|\vek{u})=p(\vek{v}|\vek{u},\vek{\varpi},\vek{\sigma}^2)p(\vek{\varpi},\vek{\sigma}^2|\vek{u})$$
in which the first factoring term is the convolution of normals
$p(\vek{v}|\vek{u},\vek{\varpi},\vek{\sigma}^2) \sim \mathcal{N}(\vek{\mu},\vek{\Sigma})$, whereas
the second factoring term $p(\vek{\varpi},\vek{\sigma}^2|\vek{u})$ cannot be computed analytically,
and thus is approximated by delta function $p(\vek{\varpi},\vek{\sigma}^2|\vek{u})\approx \delta (\vek{\varpi}_{MP},\vek{\sigma}_{MP})$.
The estimate $(\vek{\varpi}_{MP},\vek{\sigma}_{MP})$ is obtained from
$$p(\vek{\varpi},\vek{\sigma}^2|\vek{u})\propto p(\vek{u}|\vek{\varpi},\vek{\sigma}^2)p(\vek{\varpi})p(\vek{\sigma}^2)$$
by maximising the evidence (marginal likelihood)
$$p(\vek{u}|\vek{\varpi},\vek{\sigma}^2)=\int p(\vek{u}|\vek{v},\vek{\sigma}^2)p(\vek{v}|\vek{\varpi}) \textrm{d}\vek{v}$$
taking the form
\begin{eqnarray}
 p(\vek{u}|\vek{\varpi},\vek{\sigma}^2)&=&(2\pi)^{-P/2}(\vek{R})^{-1/2}\nonumber\\
 &&\textrm{exp}\left(-\frac{1}{2}\vek{u}^T(\vek{R})^{-1}\vek{u}\right)
\end{eqnarray}
in which $\vek{R}=B^{-1}+\vek{\varPsi} A^{-1}\vek{\varPsi}^T$ with $B=\vek{\sigma}^{-2}I_P$, and $A=\textrm{diag}(\varpi_\alpha)_{\alpha \in \mathcal{J}}$. By optimality criteria
$$\frac{\partial p(\vek{u}|\vek{\varpi},\vek{\sigma}^2)}{\partial \vek{\varpi}}=\vek{0}$$
and
$$\frac{\partial p(\vek{u}|\vek{\varpi},\vek{\sigma}^2)}{\partial \vek{\sigma}^2}=\vek{0}$$
one may iteratively obtain the values for $\vek{\varpi}$ and $\vek{\sigma}^2$. The number of iterations neccessary to achieve the desired accuracy 
depends on the value of the measurement noise if not marginalised, see \refig{fig_like_err}. If the signal is clean, the convergence is faster and vice versa.
Likewise, the sparsity increases with the increase of the noise magnitude, see \refig{fig_sparse_noise}. For a higher noise magnitude more 
polynomial chaos terms can be 
considered as zero, and vice versa.
In addition, the convergence also depends on the size of the data set, see \refig{fig_like_err_samp} on the example of a 
randomly chosen (i.e.~Monte Carlo)
data set. The convergence is hence faster when more samples are available.

\subsection{Sparse optimal maps}
Once the functional approximation of the forecasted state is computed, the discretisation of the coefficients of 
the forward (e.g. Jacobian $\hat{\mathring{\vek{H}}}$) and inverse maps
(i.e. Kalman gain $\hat{\vek{K}}_{k-\Delta \tau}^{(i)}$) in \refeq{refrvestim_ipce} is the only remaining operation before having 
full discretisation of the posterior variable. This can 
be simply achieved by using the direct projection method in which Jacobian and Kalman gains 
are computed directly by employing \refeq{chechkhapp} and \refeq{kalman_gain_seq}, respectively, 
and the formula for the evaluation of the respective covariance matrices:
\begin{eqnarray}
\label{cov_pce}
  \vek{C}_{q,w}&=&\mathbb{E}((\hat{\vek{q}}-\bar{\vek{q}})\otimes (\hat{\vek{w}}-\bar{\vek{w}}))\\
  &=&\sum_{\alpha,\beta \in \mathcal{K}} \mathbb{E}(\Psi_\alpha \Psi_\beta) 
 \vek{{q}}^{(\alpha)}\otimes {\vek{w}}^{(\beta)}-\bar{\vek{q}} \otimes \bar{\vek{w}}. \nonumber
\end{eqnarray}
The last relation can be further rewritten in a matrix form as
\begin{equation}
\label{cov_pcea}
 \vek{C}_{q,w}=\tilde{\vek{Q}}_f\vek{\vDelta}\tilde{\vek{W}}_f^T
\end{equation}
in which $(\vek{\vDelta})_{\alpha \beta}=\mathbb{E}(\Psi_\alpha \Psi_\beta)=\textrm{diag}(\alpha!)$ and
$\tilde{\vek{Q}}$ is equal to ${\vek{Q}}:=(...,\vek{x}^{(\alpha)},...)^T$ without the mean part. 
Similar holds for $\vek{W}$.

However, in case of  high dimensional problems this approach can be expensive. Therefore,  
the estimation of a linearised map in \refeq{optmaplin}
can be rephrased in a similar setting as described in the previous section. 
Given samples of the a posteriori estimate of the state $\vek{x}_{na}^{(i)}(\omega_j)$
in $i$-th iteration, one may evaluate the samples of the measurement forecast 
$\vek{u}(\omega)=[Y(\vek{x}_{na}^{(i)}(\omega_j))]_{j=1}^N$ such that 
\begin{equation}
\label{regression_measur}
 {\vek{u}}(\omega)=\sum_{\alpha \in \mathcal{K}} \vek{u}^{(\alpha)}\Psi_\alpha(\vek{\xi}(\omega))=\vek{\Psi}\vek{v}
\end{equation}
holds. Hence, Bayesian regression as introduced earlier can be used for the estimation of unknown sparse coefficients $\vek{v}$. In this regard
\begin{equation}
Y_{k}(\vek{x}_{na}^{(i)})=\mathring{\vek{H}}^{(i)}(\vek{x}_{na}^{(i)}-\check{\vek{x}}^{(i)})+\mathring{\vek{h}}_k+\vek{\epsilon}_k
\end{equation}
holds in which $\mathring{\vek{H}}^{(i)},\mathring{\vek{h}}_k$ and $\epsilon$ are unknown, and are to be estimated from underdetermined data. 
However, in contrast to the problem in the previous section, here one aims at estimating the matrix parameter type. 
To reduce the estimation to the same form as in \refeq{regression_measur}, one may vectorise the previous equation to
\begin{equation}
\vek{u}_x=\vek{X}\vek{q}_x+\vek{\epsilon}_x
\end{equation}
in which
\begin{eqnarray}\label{state_x_app}
 \vek{X}&=&[\vek{1}\quad (\vek{x}_{na}^{(i)}(\omega_j)-\check{\vek{x}}^{(i)})^T]_{j=1}^N \in \mathbb{R}^{N\times (d+1)}\nonumber\\
 \vek{q}_x&=&[\vek{h}_x;(\check{H}^{(i)}(1,:))^T]
\end{eqnarray}
and $\vek{\epsilon}_x=\epsilon_x\otimes \vek{e}$. Here, $\vek{x}_{na}(\omega_j)\in \mathbb{R}^d$ is the state sample, 
$\vek{e}=[1,0,...,0]^T$ and $\epsilon_x$ is the approximation error of the first state. Similarly, one may write
\begin{eqnarray*}\label{other_states}
\vek{u}_y&=&\vek{X}\vek{q}_y+\vek{\epsilon}_y\nonumber\\
\vek{u}_z&=&\vek{X}\vek{q}_z+\vek{\epsilon}_z.
\end{eqnarray*}
In these forms Eqs.~(\ref{state_x_app})-~(\ref{other_states}) can be also solved in a sparse Bayesian setting. 

\begin{figure*}
\begin{center}
\includegraphics{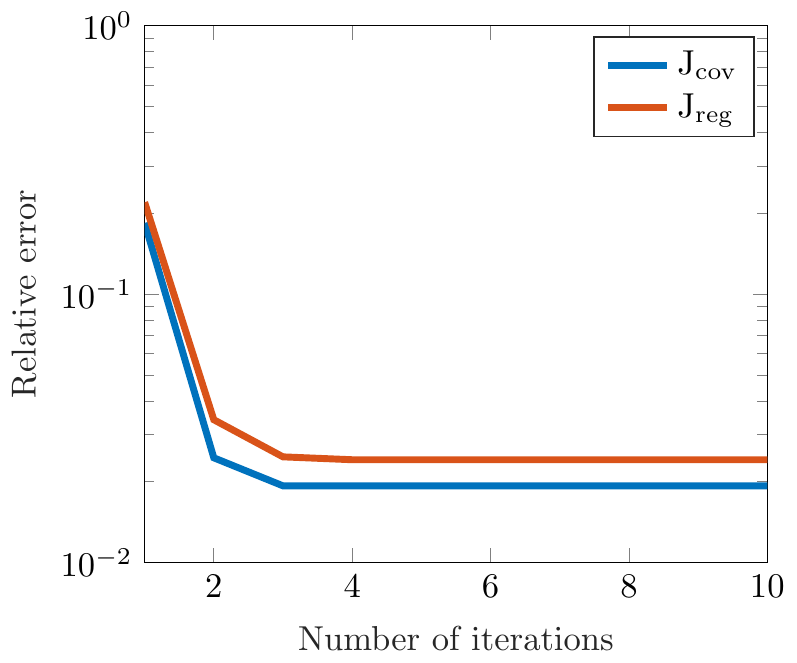}
\includegraphics{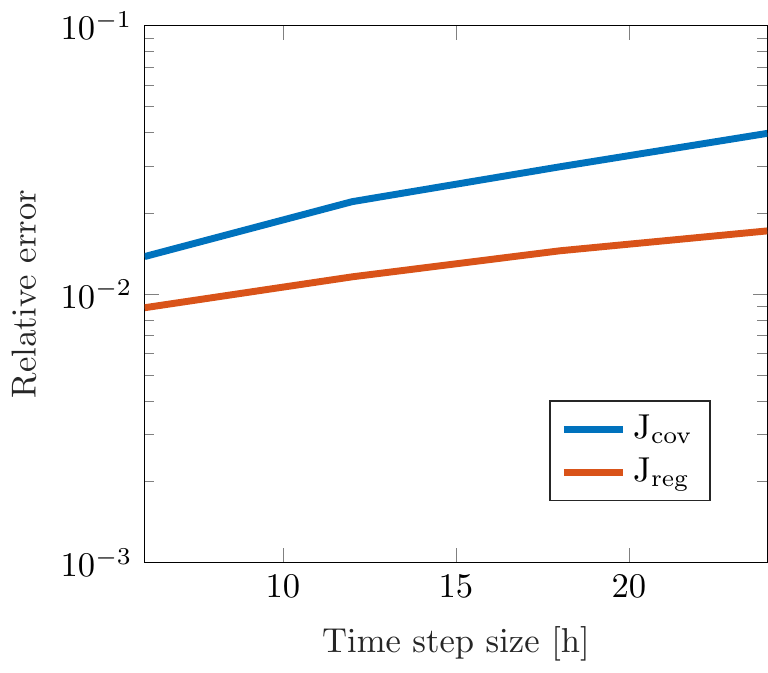}
\caption{Accuracy of the approximated Jacobian compared to the exact Jacobian}
\label{fig_jacobian_approx}
\end{center}
\end{figure*}

The Jacobian estimated in this manner is slightly better than the estimate obtained using \refeq{chechkhapp} as can be seen in \refig{fig_jacobian_approx}.
Here, the relative errors of regression ($J_{reg}$) and covariance ($J_{cov}$) type of Jacobians compared to the analytical value of Jacobian are depicted. 
Both Jacobians 
converge very fast, already after 3 iterations, whereas their accuracy deteriorates with the increase of the length of pseudo-update step as expected. 


Besides promoting sparsity in the polynomial chaos approximations and the Jacobian, one may also use 
the Bayesian method to estimate the Kalman gain
by solving the linear system
\begin{equation}
 \vek{x}_{k-\Delta \tau,f}(\omega_j)=\vek{K}_{k-\Delta \tau}^{(i)}\vek{y}_{kf}^{(i)}(\omega_j)+\vek{b}^{(i)}+\vek{\epsilon}_K
\end{equation}
given the set of sample points $(\vek{x}_{nf}^{(i)}(\omega_j),\vek{y}_{kf}^{(i)}(\omega_j))$. Collecting samples of each of the states into vectors 
$\vek{u}_1,\vek{u}_2$ and $\vek{u}_3$ respectively for $x,y$ and $z$ one may rewrite the previous equation as
\begin{equation}\label{kalman_regression}
 \vek{u}_m=\vek{W}\vek{k}_m+\vek{\epsilon}_m, \quad m=1,..3
\end{equation}
in which $\vek{W}=[\vek{1}\quad (\vek{y}_{kf}^{(i)}(\omega_j))^T] \in \mathbb{R}^{N\times (d+1)}$ 
and $\vek{k}_m=[\vek{b}^{(i)};(\vek{K}_\ell^{(i)}(m,:))^T]\in \mathbb{R}^d$. The unknown coefficients  $\vek{k}_m$ 
can be then evaluated by using the Bayes's rule. 

\section{Conclusion}\label{concl}

We have developed the iterative incremental predictor-corrector Gauss-Newton-Markov-Kalman smoothing algorithm for
the non-Gaussian state estimation given 
noisy measurements. The method is 
based on the nonlinear local approximation of the conditional expectation, and is mathematically generalised to take into account 
possible measurement uncertainty. The resulting update equation is discretised by using the time-adaptive polynomial chaos expansion  
in terms of the standard normal random variables, the number of which matches the state dimension.
These are obtained by isoprobabilistic transformation 
of the non-Gaussian posterior random variable expressed in terms of generalised polynomials
of the last known state. The adjustment of the basis functions is achieved 
via modified Gram-Schmidt as well as nonlinear mapping algorithm such that the desired updating accuracy 
does not change when the measurement frequency is too low.  The resulting Kalman-type update formula
for the PCE coefficients can be efficiently computed solely within the PCE. As it does not
rely on sampling, the method is robust, fast and exact.

As compared to Monte Carlo, the method is not directly affected by sampling error. However, the method accuracy involves
regression error, the truncation error of polynomial approximations (PCE, approximation of optimal map and approximation of inverse map)
and errors characterising the transformation of the non-Gaussian random variables. 
The polynomial approximations here are all evaluated in a data learning setting via Bayes's rule given randomly chosen samples. This 
may lead to potential over-estimation of some of the polynomial coefficients. 
However, note that the PCE approximations can be easily exchanged with a fully deterministic Galerkin
algorithm for the state estimation obtained given the variational form of the stochastic ordinary differential
equations as previously studied
in \cite{Pajonk2011c}.  

The updating procedure has been applied to a low-dimensional state estimation problem of the chaotic
Lorenz-84 system. It is shown that the method is robust and able to estimate the initial state of the Lorenz-84 system 
even when the updating step is large and the measurement noise is high. The extension of the presented method to more 
realistic applications is currently ongoing research. As the numerical complexity of the method increases with the state dimension,
the future plan is to consider low-rank techniques as well as to implement more efficient adaptive sampling strategies.
This would then allow the use of quadratic approximations in the iterative form. Finally, the proposed method is based on the 
approximation of the conditional expectation of the state and not its higher moments. The further step is to also include the conditional 
expectation of the second moment into the updating process as well.

\textbf{Acknowledgment} The author greatly appreciates partial financial funding by the German Science Foundation (Deutsche Forschungsgemeinschaft, DFG) as part of
priority programs GRK 2075, SPP 1886 and SPP 1748. 
%
%
%
%


%
%

\section{Appendix: The Lorenz 1984 system}\label{ss:example}

For the numerical evaluation of the estimation method described in the previous, here we consider the well-known ``Lorenz-84'' model
\cite{Lorenz1984,Lorenz2005}. It is described by a set of three state variables
$\OPstatevector=(x,y,z)^T$. Here $x$ represents a symmetric, globally
averaged westerly wind current, whereas $y$ and $z$ represent the cosine and
sine phases of a chain of superposed large-scale eddies transporting heat
polewards. The state evolution 
is described by the following set of
ordinary differential equations (ODEs):
\begin{eqnarray}
\frac{dx}{dt}&=&-ax-y^2-z^2+aF_1 \nonumber\\
\frac{dy}{dt}&=&-y+xy-bxz+F_2 \label{eqn:lorenz84det}\\
\frac{dz}{dt}&=&-z+xz+bxy, \nonumber
\end{eqnarray}
in which $F_1$ and $F_2$ represent known thermal forcings, and $a$ and $b$ are
fixed constants. 

In the numerical experiment considered in the paper the initial condition of the
``unknown truth''~is $( 1 . 0 , 0 . 0 , − 0 . 75 )$, the thermal forcings
are set to $F_1 = 8$ and $F_2 = 1$, whereas the parameters are set to
$a =0.25$ and $b = 4$. Given the initial values, the previous system is integrated forward in
time using an adaptive embedded Runge-Kutta (RK) scheme of orders 4 and 5.

As the Lorenz-84 model shows chaotic behaviour and is very sensitive to the initial
conditions, we model them as independent Gaussian random variables:
\begin{eqnarray}
x_0(\omega)&\sim& \mathcal{N}(x_0,\sigma_1)\nonumber\\
y_0(\omega)&\sim& \mathcal{N}(y_0,\sigma_2)\label{eq:iclorenz84}\\
z_0(\omega)&\sim& \mathcal{N}(z_0,\sigma_3)\nonumber
\end{eqnarray}
with $x_0 = y_0 = z_0 = 0$ and standard deviations $\sigma_1=\sigma_2=\sigma_3=1$.

%

\bibliographystyle{plain}
\bibliography{main2pdf}

\begin{thebibliography}{10}

\bibitem{banerjee}
A.~Banerjee, X.~Guo, and H.~Wang.
\newblock On the optimality of conditional expectation as a bregman predictor.
\newblock {\em {IEEE} Trans. Information Theory}, 51(7):2664--2669, 2005.

\bibitem{bell}
B.~M. Bell.
\newblock The iterated {Kalman Smoother as a Gauss-Newton Method}.
\newblock {\em SIAM Journal on Optimization}, 4(3):626--636, 1994.

\bibitem{ABobrowski}
A.~Bobrowski.
\newblock {\em Functional analysis for probability and stochastic processes: an
  introduction}.
\newblock Cambridge University Press, Cambridge, Cambridge, UK, 2005.

\bibitem{botev2010}
Z.~I. Botev, J.~F. Grotowski, and D.~P. Kroese.
\newblock Kernel density estimation via diffusion.
\newblock {\em Ann. Statist.}, 38(5):2916--2957, 10 2010.

\bibitem{bregman}
L.M. Bregman.
\newblock The relaxation method of finding the common point of convex sets and
  its application to the solution of problems in convex programming.
\newblock {\em USSR Computational Mathematics and Mathematical Physics},
  7(3):200 -- 217, 1967.

\bibitem{chen2012ensemble}
Y.~Chen and D.~S. Oliver.
\newblock Ensemble randomized maximum likelihood method as an iterative
  ensemble smoother.
\newblock {\em Mathematical Geosciences}, 44(1):1--26, 2012.

\bibitem{reich2016}
N.~Chustagulprom, S.~Reich, and M.~Reinhardt.
\newblock A hybrid ensemble transform particle filter for nonlinear and
  spatially extended dynamical systems.
\newblock {\em SIAM/ASA Journal on Uncertainty Quantification}, 4(1):592--608,
  2016.

\bibitem{einicke2012smoothing}
G.~A. Einicke.
\newblock {\em Smoothing, filtering and prediction: estimating the past,
  present and future}.
\newblock InTech, 2012.

\bibitem{einicke1999robust}
G.~A. Einicke and B.~Langford.
\newblock Robust extended {K}alman filtering.
\newblock {\em IEEE Transactions on Signal Processing}, 47(9):2596--2599, 1999.

\bibitem{gamerman1997}
D.~Gamerman.
\newblock {\em {Markov chain Monte Carlo: stochastic simulation for Bayesian
  inference}}.
\newblock Chapman \& Hall, Boca Raton, USA, 2 edition, May 1997.

\bibitem{Gerritsma:2010}
Marc Gerritsma, Jan-Bart van~der Steen, Peter Vos, and George Karniadakis.
\newblock Time-dependent generalized polynomial chaos.
\newblock {\em J. Comput. Phys.}, 229(22):8333--8363, November 2010.

\bibitem{gratton}
S.~Gratton, A.~S. Lawless, and N.~K. Nichols.
\newblock Approximate {Gauss-Newton} methods for nonlinear least squares
  problems.
\newblock {\em SIAM Journal on Optimization}, 18(1):106--132, 2007.

\bibitem{kai2010robust}
X.~Kai, C.~Wei, and L.~Liu.
\newblock Robust extended {K}alman filtering for nonlinear systems with
  stochastic uncertainties.
\newblock {\em IEEE Transactions on Systems, Man, and Cybernetics-Part A:
  Systems and Humans}, 40(2):399--405, 2010.

\bibitem{kalman}
R.E. Kalman.
\newblock A new approach to linear filtering and prediction problems.
\newblock {\em ASME. J. Basic Eng.}, 82(1):35--45, 1960.

\bibitem{Lorenz1984}
E.~N. Lorenz.
\newblock Irregularity: a fundamental property of the atmosphere.
\newblock {\em Tellus A}, 36(2):98--110, 1984.

\bibitem{Lorenz2005}
Edward~N. Lorenz.
\newblock A look at some details of the growth of initial uncertainties.
\newblock {\em Tellus A}, 57(1):1--11, 2005.

\bibitem{Matthies2016}
H.~G. Matthies, E.~Zander, B.~Rosi\'{c}, and A.~Litvinenko.
\newblock Parameter estimation via conditional expectation: a {B}ayesian
  inversion.
\newblock {\em Advanced Modeling and Simulation in Engineering Sciences},
  3(1):1--21, 2016.

\bibitem{van2001square}
R.~Van~Der Merwe and E.~A. Wan.
\newblock The square-root unscented {K}alman filter for state and
  parameter-estimation.
\newblock In {\em Acoustics, Speech, and Signal Processing, 2001.
  Proceedings.(ICASSP'01)}, volume~6, pages 3461--3464. IEEE, 2001.

\bibitem{Pajonk2011c}
O.~Pajonk, B.~Rosi\'{c}, A.~Litvinenko, and H.~G. Matthies.
\newblock A deterministic filter for non-{G}aussian {B}ayesian estimation.
\newblock {\em Physica D: Nonlinear Phenomena}, 241(7):775--788, 2012.

\bibitem{ristic}
B.~Risti\'{c}, S.~Aurlampalam, and N.~Gordon.
\newblock {\em Beyond the {K}alman filter: particle filters for tracking
  applications}.
\newblock Artech House Publishers, Boston, 2004.

\bibitem{BRrosic12}
B.~Rosi\'{c}, O.~Pajonk, A.~Litvinenko, and H.~G. Matthies.
\newblock Sampling-free linear {B}ayesian update of polynomial chaos
  represenations.
\newblock {\em Journal of Computational Physics}, 231(17):5761--5787, 2012.

\bibitem{sakov2012iterative}
P.~Sakov, D.~Oliver, and L.~Bertino.
\newblock An iterative {EnKF} for strongly nonlinear systems.
\newblock {\em Monthly Weather Review}, 140(6):1988--2004, 2012.

\bibitem{simon2006optimal}
D.~Simon.
\newblock {\em Optimal state estimation: {K}alman, {H} infinity, and nonlinear
  approaches}.
\newblock John Wiley \& Sons, 2006.

\bibitem{Smith93}
A.~F.~M. Smith and G.~O. Roberts.
\newblock Bayesian computation via the {G}ibbs sampler and related {M}arkov
  chain {M}onte {C}arlo methods.
\newblock {\em Journal of the Royal Statistical Society. Series B
  (Methodological)}, 55(1):3--23, 1993.

\bibitem{Bazargan2013}
H.~A. Tchelepi, H.~Bazargan, and M.~A. Christie.
\newblock Efficient {M}arkov chain {M}onte {C}arlo sampling using polynomial
  chaos expansion.
\newblock In {\em Proceedings of the SPE Reservoir Simulation Symposium}, The
  Woodlands, Texas, United States, 2013. online.

\bibitem{Tipping}
M.~E. Tipping.
\newblock Sparse bayesian learning and the relevance vector machine.
\newblock {\em Journal of Machine Learning Research}, 1:211--244, 2001.

\bibitem{wan2000unscented}
E.~A. Wan and R.~Van~Der Merwe.
\newblock The unscented {Kalman} filter for nonlinear estimation.
\newblock In {\em Adaptive Systems for Signal Processing, Communications, and
  Control Symposium 2000. AS-SPCC. The IEEE 2000}, pages 153--158. IEEE, 2000.

\bibitem{wang2016randomized}
K.~Wang, T.~Bui-Thanh, and O.~Ghattas.
\newblock A randomized maximum a posterior method for posterior sampling of
  high dimensional nonlinear {B}ayesian inverse problems.
\newblock {\em arXiv preprint arXiv:1602.03658}, 2016.

\bibitem{Xiu2010}
Dongbin Xiu.
\newblock {\em Numerical Methods for Stochastic Computations: A Spectral Method
  Approach}.
\newblock Princeton University Press, Princeton, NJ, USA, 2010.

\end{thebibliography}

\end{document}